\documentclass[12pt,a4paper]{article}
\baselinestretch
\usepackage[english]{babel}
\usepackage[title]{appendix}
\usepackage{amsbsy, amsmath, amssymb, multirow,epsf}
\usepackage{setspace}
\usepackage{amsthm}
\usepackage{cite}
\usepackage{graphicx}
\usepackage{subcaption}
\usepackage{xcolor}
\usepackage[normalem]{ulem}
\usepackage{enumitem}
\usepackage{tikz}
\usepackage{xcolor}
\usepackage{bbm}
\usetikzlibrary{arrows.meta,positioning,shapes.geometric}

\theoremstyle{plain}
\newtheorem{theorem}{Theorem}[section]
\newtheorem{proposition}[theorem]{Proposition}
\newtheorem{corollary}[theorem]{Corollary}
\newtheorem{lemma}[theorem]{Lemma}

\theoremstyle{definition}
\newtheorem{definition}{Definition}[section]

\newtheorem{assumption}{Assumption}[section]

\theoremstyle{remark}
\newtheorem{remark}{Remark}[section]

\numberwithin{equation}{section}
\usepackage{authblk} 
\usepackage{bm}
\usepackage{threeparttable}
\usepackage{float}
\usepackage{epstopdf}
\usepackage{hyperref}

\usepackage{array}
\usepackage{tabularx}
\usepackage{booktabs}

\newcolumntype{Y}{>{\raggedright\arraybackslash}X}

\DeclareMathOperator{\V@R}{V@R}

\DeclareMathOperator{\TV@R}{TV@R}

\allowdisplaybreaks[4]

\title{Asymptotic Diversification of Value-at-Risk with Heterogeneous Tails}
\author[1]{Yongzhao Chen\footnote{Email: ychen@hsu.edu.hk}}
\author[2]{Ka Chun Cheung\footnote{Email: kccg@hku.hk}}
\author[3]{Hok Kan Ling\footnote{Email: brian.ling@queensu.ca}}
\author[4]{Qihe Tang\footnote{Email: qihe.tang@unsw.edu.au}}
\author[5]{Sheung Chi Phillip Yam\footnote{Email: scpyam@sta.cuhk.edu.hk}}
\affil[1]{Department of Mathematics, Statistics and Insurance\\ The Hang Seng University of Hong Kong}
\affil[2]{Department of Statistics and Actuarial Science\\ The University of Hong Kong}
\affil[3]{Department of Mathematics and Statistics\\ Queen's University}
\affil[4]{School of Risk and Actuarial Studies\\ UNSW Business School}
\affil[5]{Department of Statistics and Data Science\\ The Chinese University of Hong Kong}

\begin{document}
\newcommand{\mP}{\mathbb{P}}
	\maketitle

\begin{abstract}
In this article, we develop a general asymptotic diversification theory for Value-at-Risk (V@R) in
heterogeneous portfolios without imposing global tail equivalence among the marginal risks. When several risks dominating the portfolio tail belong to the Fr\'echet
maximum domain of attraction  with a common tail index $\alpha$, and their
survival copula is Archimedean with a strict generator regularly varying at
zero, V@R is asymptotically: i) superadditive if $0<\alpha<1$; ii) additive if
$\alpha=1$; and iii) subadditive if $\alpha>1$;  at the same time, all lighter-tailed risks
are asymptotically negligible. Moreover, under the same Archimedean dependence assumption, with generator index
$-\gamma<0$, we derive aggregate-to-marginal tail limits for common-index
Fr\'echet marginals along sequences with positive finite limiting marginal tail
ratios. These limits i) decrease strictly in $\gamma$ for
$0<\alpha<1$; ii) are
constant for $\alpha=1$; and iii) increase strictly in $\gamma$ for $\alpha>1$.  Under the same Archimedean dependence assumption as above,
we also show that aggregation preserves the corresponding Fr\'echet maximum domain of
attraction for common-index marginals. For heterogeneous risks in the Gumbel maximum domain of attraction with
infinite right endpoints, we also establish asymptotic subadditivity under
arbitrary dependence and asymptotic additivity under positive diagonal
upper-tail dependence.

\end{abstract}

\section{Introduction}
\label{sec:introduction}
Value-at-Risk (V@R) is one of the most widely used quantile-based
measures of financial and insurance risk. For a loss random variable $X$
with distribution function $F$, its V@R at confidence level $p\in(0,1)$ is
defined by
\[
\V@R_p(X)
:=
F^{\leftarrow}(p)
=
\inf\{x\in\mathbb R:F(x)\geq p\};
\]
see, for example, \cite{ccy25,mcneil2015quantitative}.
Throughout the paper, let $X_1,\ldots,X_d\geq0$ denote nonnegative loss
random variables, where $d\geq2$ is fixed. For each
$i\in\{1,\ldots,d\}$, let $F_i$ denote the distribution function of
$X_i$ and write
\[
\overline F_i(x)
:=
1-F_i(x)
=
\mathbb P(X_i>x)
\]
for its survival function. Denote the aggregate loss by
\[
S_d:=\sum_{i=1}^d X_i.
\]
Whenever the denominator is positive, define the V@R diversification ratio
by
\begin{equation}\label{eq:diversification-ratio}
D_p
:=
\frac{\V@R_p(S_d)}
     {\sum_{i=1}^d\V@R_p(X_i)}.
\end{equation}
Subadditivity corresponds to $D_p\leq1$, meaning that aggregation does not
increase the required capital beyond the sum of the stand-alone capital
requirements, whereas superadditivity corresponds to $D_p\geq1$, meaning
that the aggregate capital requirement is at least as large as the sum of
the stand-alone requirements. The boundary case $D_p=1$ corresponds to
additivity. V@R is neither subadditive nor superadditive in general, and
its aggregation behavior depends jointly on the marginal tails and the
dependence structure; see \cite{Embrechts00} and the references therein.
According to Basel Accords and regulations alike, we usually focus on portfolio tail risks, especially for large $p$, let say 99.5\% or 99.9\%; although $D_p$ may not be uniformly above or below 1 for small $p$, the limiting behavior of $D_p$ as $p\uparrow 1$ should still be instructive for the tail V@R when $p$ is large. Here, we study this interaction in the high-confidence regime $p\uparrow1$ with the practically most relevant and general enough setting as described below.

Archimedean copulae constitute a widely used class of multivariate dependence models,
including the Clayton, Gumbel–Hougaard, Frank, Joe, and Ali–Mikhail–Haq families, which accommodate different dependence strengths and tail behaviors, including lower-tail dependence, upper-tail dependence, and asymptotic independence; also see \cite{mc09}. These features have led to applications across several disciplines, including multivariate survival analysis in biostatistics, loss modelling in actuarial science, and risk modelling in quantitative finance, also see \cite{n06}. Under Archimedean dependence,
\cite{wuthrich03} studied the tail behavior and high-level V@R of sums of
dependent risks, and \cite{alw04} further analyzed the resulting first-order
tail relation:
\[
\mathbb{P}\left(\sum_{i=1}^d X_i>x\right)
\sim
q_d \overline F(x),
\qquad x\to\infty,
\]
where $q_d$ is the aggregate-tail constant measuring the asymptotic size of
the tail probability of the sum relative to that of a single marginal risk.
Related extensions to more general copula and multivariate extreme-value
structures were developed in \cite{alw07,bfg06}. Our approach of how to deal with heavy-tail portfolios relies on various common tools from extreme-value theory and regular variation, and it is also related to the well-received work \cite{enw08}; in particular, in \cite{enw08}, for identically distributed
regularly varying risks under Archimedean survival dependence, the authors
characterized $q_d$ and showed that it is increasing in the regular-variation
index of the generator when the marginal tail index exceeds 1, decreasing
when the marginal tail index is below 1, and constant at the boundary value
1. Through the corresponding quantile asymptotics, this yields the familiar
transition from asymptotic subadditivity for $\alpha>1$ to asymptotic
superadditivity for $0<\alpha<1$. Complementary analyses based on
multivariate regular variation can be found in
\cite{Embrechts09a,djssv13,mainik13}.

A common restriction in explicit aggregate-tail and diversification-ratio
formulas is the use of identically distributed or, more generally,
tail-equivalent marginals, meaning that
\[
\frac{\overline F_i(x)}{\overline F_j(x)}
\rightarrow c_{ij}\in(0,\infty),
\qquad x\to\infty.
\]
This restriction excludes two distinct forms of marginal tail heterogeneity.
First, even if all marginal survival functions are regularly varying with a
common tail index $\alpha$, their slowly varying factors may differ in such a way
that these ratios do not converge. Second, if the marginal tail orders differ,
only a subset of the portfolio may remain relevant in the far tail. Sums of
dependent, non-identically distributed heavy-tailed risks have been studied
by \cite{kortschak09}, including settings in which some marginal tails are
asymptotically negligible, while \cite{cmh13} extended earlier aggregation
results to non-identical but tail-equivalent marginals. Some other works such as \cite{yz23} impose more
general tail-order conditions on the dependence structure, but with identical
regularly varying marginals. These studies provide
first-order or higher-order asymptotic descriptions of the tail probability
of the aggregate under their respective assumptions. They do not, however,
directly yield the V@R classification pursued here when
global tail-equivalence assumptions are not imposed and the marginal tail
ratios need not converge, and this is precisely the problem we aim to tackle in this work.

A different strand of the literature derives asymptotic additivity from extremal dependence.  \cite{hj12,hj12b} introduced and
studied upper tail comonotonicity and obtained asymptotic additivity for
tail-equivalent risks in the Fr\'echet or Gumbel maximum domain of
attraction (MDA). \cite{cltyy18} removed the tail-equivalence
requirement and proved asymptotic additivity for heterogeneous portfolios
whose marginals may belong to any of the three MDAs or have finite right endpoints. Upper tail comonotonicity,
however, represents asymptotically maximal positive dependence. Under extreme-value copulae, \cite{ccyyz23} classified asymptotic
subadditive and superadditive regimes for identically distributed marginals. Aggregation of asymptotically independent
light-tailed risks has also been studied in \cite{nolde14}. Our present study has a similar spirit to \cite{zlyxs23}, which developed a general
framework for asymptotic V@R subadditivity and superadditivity under copulae
admitting a stable tail dependence function (STDF). Their dependence
assumption is more general than the Archimedean survival-dependence
condition considered here; in particular, a survival Archimedean copula
with a regularly varying generator as in our setting admits an STDF.
The two approaches differ more substantially in their treatment of the
marginals. Theorem 3.1 of \cite{zlyxs23} requires all marginal risks to be
ordered in the same direction, through the tail concave order, relative to
regularly varying reference risks with index $-1$. In particular, for
regularly varying marginals, it does not directly cover a portfolio containing
risks on both sides of the critical tail index one. Our approach instead
identifies the asymptotically dominating Fr\'echet risks: if several such
risks share a common tail index $\alpha$, they determine the asymptotic V@R
classification, while lighter-tailed risks can still be included in a negligible manner. Moreover, Proposition 3.5(ii) of \cite{zlyxs23}, which also
reduces a portfolio to a set of dominating risks, assumes asymptotically
proportional tails among those risks, whereas our results allow their
relative tail scales to oscillate without a global limit. 

When several risks dominating the portfolio tail belong to the Fr\'echet
MDA, we model the dependence among the
asymptotically relevant risks through an Archimedean survival copula with
a strict generator $\psi$ that is regularly varying at zero with index
$-\gamma<0$. This condition specifies only the local behavior of the
generator that governs joint upper-tail events. Connections between
Archimedean tail behavior and regular variation of the generator have
been studied in \cite{jw02,cs07,charpentier2009tails}; see also
\cite{mc09} for the general structure of multivariate Archimedean
copulae. The major result in our work for this Fr\'echet case is to establish the V@R classification
by isolating the asymptotically relevant risks rather than imposing a common
comparison condition on the entire portfolio. Suppose that these dominating
risks share a common Fr\'echet tail index $\alpha$. Under the Archimedean
survival-dependence condition above, we prove in Theorem \ref{thm:main} that:
\[
\left\{\begin{aligned}
   &\ \liminf_{p\uparrow1}D_p\geq1,&&\text{ when }  0<\alpha<1;\\
&\  D_p\to1,&&\text{ when } \alpha=1;\\
&\  \limsup_{p\uparrow1}D_p\leq1,
&&\text{ when } \alpha>1.
\end{aligned}\right.
\]
We derive the aggregate-tail asymptotics underlying this dependence structure.
For common-index Fr\'echet marginals, we introduce sequential tail ratios
that describe their relative scales along diverging sequences. Along any
sequence on which all relevant pairwise ratios converge to positive finite
limits, we obtain an explicit integral representation for the limiting ratio
between the tail probability of the aggregate and a marginal tail
probability. This formulation permits slowly varying marginal factors whose
ratios oscillate and have no global limits.

We also establish strict dependence monotonicity of the resulting
aggregate-tail constant. It is i) strictly
decreasing in $\gamma$ when $0<\alpha<1$; ii) independent of $\gamma$ when $\alpha=1$; and iii) strictly increasing in $\gamma$ when
$\alpha>1$. This extends the monotonicity
analysis of \cite{enw08} from identical to non-identically distributed
marginals and establishes strictness for every $\alpha\neq1$ and
$d\geq2$; see Theorem \ref{thm:limsup-liminf-sharp}. Its proposed proof here uses a convex-order comparison and a coordinate-wise
replacement argument rather than the supermodular-order argument used in
\cite{enw08}. The same representation gives a strict comparison with the
comonotonic benchmark for every finite $\gamma$ and yields closure of the
Fr\'echet maximum domain of attraction under aggregation without global
tail equivalence. Finally, if only a fixed subset of the portfolio has the dominant common
Fr\'echet index, the remaining risks are asymptotically negligible and the
same classification is determined by the dominant subportfolio. The
single-dominating-risk case, which is already known, appears as the
boundary situation in which the diversification ratio converges to one
under arbitrary dependence.

Besides, we also study heterogeneous risks in the
Gumbel MDA with infinite right endpoints. We
first show in Theorem \ref{thm:gumbel-subadditivity}, without imposing any dependence assumption, that
\[
\limsup_{p\uparrow1}D_p\leq1.
\]
Thus, V@R is asymptotically subadditive under arbitrary dependence, generalizing the result in \cite{ccyyz23}. If
the diagonal upper-tail dependence coefficient
\[
\lambda_\Delta
:=
\lim_{u\downarrow0}
\frac{\widehat C(u,\ldots,u)}{u}
\]
exists and is strictly positive, then the conclusion strengthens to that in Theorem \ref{thm:gumbel-tail-dependent}, namely:
\[
D_p\to1.
\]
Compared with Proposition 4.1 of \cite{zlyxs23}, this result does not
require log-regular variation, a common log-regular-variation index,
tail equivalence, or continuity of the marginal distributions.
Archimedean survival dependence with a generator regularly varying at
zero automatically satisfies this condition, with
$\lambda_\Delta=d^{-\frac{1}{\gamma}}>0$.  Finally, as in Theorem \ref{thm:gumbel-mixed}, finite-endpoint additions do not
affect the Gumbel subadditivity conclusion, while if all risks have
finite right endpoints, both the V@R of the aggregate and the sum of the
marginal V@Rs converge to the sum of those endpoints under Archimedean survival dependence with a strict generator. In summary, for the sake of convenience, Table \ref{tab:main-results} summarizes the main
asymptotic V@R aggregation results across the marginal and dependence
regimes considered in this paper.

The remainder of the paper is organized as follows.
Section \ref{sec:model-prelim} introduces the Archimedean survival-dependence
framework and recalls the required facts from regular variation and
extreme-value theory. Section \ref{sect:var_aa} presents the main asymptotic
aggregation results. Section \ref{sec:frechet} treats Fr\'echet risks, including
common-index marginals without tail equivalence and portfolios with one or
several dominating Fr\'echet risks. Section \ref{sec:gumbel-infinite} considers
Gumbel risks with infinite right endpoints, including subadditivity under
arbitrary dependence, additivity under positive upper-tail dependence, and related extensions. Section \ref{sec:num} provides numerical illustrations, then
Section \ref{sec:discussion} discusses the results and concludes. Various technical proofs are presented in the appendices.

\section{Model and Preliminaries}
\label{sec:model-prelim}

\subsection{Archimedean Survival Dependence}
\label{sec:model}
We model the dependence among the risks through their survival copula. Specifically, if $\widehat C$ denotes the survival copula of
$(X_1,\ldots,X_d)$, then
\[
\mathbb P(X_1>x_1,\ldots,X_d>x_d)
=
\widehat C\bigl(
\overline F_1(x_1),\ldots,\overline F_d(x_d)
\bigr).
\]
Following the generator convention in \cite{n06}, let
$\psi:[0,1]\to[0,\infty]$ be a continuous and strictly decreasing function
such that $\psi(1)=0$. A generator is called \emph{strict} if
$\psi(0)=\infty$. The corresponding Archimedean copula is
\[
C^\psi(u_1,\ldots,u_d)
=
\psi^{-1}\left(
\sum_{i=1}^d\psi(u_i)
\right),
\qquad
(u_1,\ldots,u_d)\in[0,1]^d.
\]
By applying Sklar's theorem (see Chapter 2 of \cite{n06}) to the survival-copula representation, when
some marginal distributions are discontinuous, the survival copula is uniquely
determined by the joint distribution only on the Cartesian product of the
ranges of the marginal survival functions; its extension to all of $[0,1]^d$
need not be unique \cite[Theorem 2.10.9]{n06}. Throughout this article, a
survival copula $\widehat C$ is therefore regarded as part of the model
specification. All values of $\widehat C$ evaluated at actual marginal survival
probabilities are uniquely determined by the joint distribution. Conditions
involving values outside this Cartesian product, such as its diagonal behavior,
refer to the specific survival-copula extension. Equivalently, for results
stated in terms of the existence of such a copula, it is sufficient that at
least one survival-copula extension satisfy any stated condition.

Throughout this article, it is understood that $\psi$ satisfies the conditions
required for $C^\psi$ to be a $d$-dimensional copula. Under the present
convention, this is equivalent to $\psi^{-1}$ being $d$-monotone on
$[0,\infty)$; see \cite{mc09}. Complete monotonicity of $\psi^{-1}$
corresponds to the stronger requirement that the same generator defines an
Archimedean copula in every dimension.

In the rest of this work, we make the following mild assumption on dependence structure:

\begin{assumption}[Regularly varying Archimedean survival dependence]
\label{AssumeArchi}
The survival copula of $(X_1,\ldots,X_d)$ is Archimedean with a strict
generator $\psi$.
Moreover, $\psi$ is regularly varying at zero with index $-\gamma<0$, i.e.\ 
for every $t>0$,
\[
\frac{\psi(tu)}{\psi(u)}
\rightarrow
t^{-\gamma},
\qquad u\downarrow0.
\]
\end{assumption}
Under Assumption \ref{AssumeArchi},
\[
\widehat C(u,\ldots,u)
=
\psi^{-1}\bigl(d\psi(u)\bigr).
\]
Clearly $\psi^{-1}\in\operatorname{RV}_{-1/\gamma}$ at infinity, and we also have:
\[
\frac{\widehat C(u,\ldots,u)}{u}
=
\frac{
\psi^{-1}\bigl(d\psi(u)\bigr)
}{
\psi^{-1}\bigl(\psi(u)\bigr)
}
\rightarrow
d^{-\frac{1}{\gamma}},
\qquad u\downarrow0.
\]
Thus larger values of $\gamma$ correspond to larger diagonal upper-tail
dependence coefficients, with
\[
d^{-\frac{1}{\gamma}}\uparrow1
\qquad\text{as }\gamma\uparrow\infty.
\]




\subsection{Extreme-Value Domains and Tail Dominance}
\label{sec:evt-tail-dominance}

We next state the extreme-value and tail-comparison facts used throughout
the paper. For a distribution function $F$, denote its right endpoint by
\[
x_F:=\sup\{x\in\mathbb R:F(x)<1\}\leq\infty.
\]
A distribution $F$ is said to belong to the MDA of
an extreme-value distribution $G$, written $F\in\operatorname{MDA}(G)$, if
there exist constants $a_n>0$ and $b_n\in\mathbb R$ such that
\[
F^n(a_nx+b_n)\to G(x),
\]
at every continuity point of $G$. The three possible nondegenerate
extreme-value limits are the Fr\'echet, Gumbel, and Weibull distributions,
denoted respectively by $\Phi_\alpha$, $\Lambda$, and $\Psi_\alpha$. For the purpose of this paper, only a few standard characterizations of
these domains are needed; see \cite{bgt87,cltyy18} for example. We first recall the following:
\[
X\in\operatorname{MDA}(\Phi_\alpha),
\qquad \alpha>0,
\]
if and only if $x_F=\infty$ and
\[
\overline F\in\operatorname{RV}_{-\alpha};
\]
equivalently, for every $t>0$,
\[
\frac{\overline F(tx)}{\overline F(x)}
\rightarrow
t^{-\alpha},
\qquad x\to\infty.
\]
Thus one may write
\[
\overline F(x)=x^{-\alpha}L(x),
\]
where $L$ is slowly varying. For nonnegative risks, $\alpha>1$ implies a
finite mean, whereas $0<\alpha<1$ implies an infinite mean. The boundary
case $\alpha=1$ may have either a finite or an infinite mean depending on
the slowly varying factor; also see \cite{bgt87}.

Besides, if $X\in\operatorname{MDA}(\Lambda)$, then there exist sequences
$a_n>0$ and $b_n\in\mathbb R$ such that
\[
n\overline F(b_n+a_nx)
\rightarrow e^{-x},
\qquad x\in\mathbb R.
\]
In addition, if $X$ has an infinite right endpoint, then its survival
function is rapidly varying; particularly, for every $t>1$,
\[
\frac{\overline F(tx)}{\overline F(x)}
\rightarrow0,
\qquad x\to\infty,
\]
whereas the ratio tends to infinity for $0<t<1$. This distinction between
regularly varying Fr\'echet tails and rapidly varying Gumbel tails will be
used repeatedly below. 

Finally, if $x_F<\infty$, then
\[
\V@R_p(X)\rightarrow x_F,
\qquad p\uparrow1.
\]
Our results for finite-endpoint risks do not require such risks to belong to
a particular MDA. Thus the final finite-endpoint
regime considered below includes, but is not restricted to, risks in the
Weibull MDA and risks in the Gumbel MDA with finite right endpoints.

Regarding the comparison of risks having different upper-tail orders, namely for two
risks $X_i$ and $X_j$ with $x_{F_j}=\infty$, we say that $X_j$
\emph{dominates} $X_i$ in the upper tail if
\[
\frac{\overline F_i(x)}
     {\overline F_j(x)}
\rightarrow0,
\qquad x\to\infty.
\]
For a Fr\'echet risk $X_j\in\operatorname{MDA}(\Phi_\alpha)$, this
dominance relation holds, in particular, whenever $X_i$: i) belongs to
$\operatorname{MDA}(\Phi_\beta)$ with $\beta>\alpha$; ii) belongs to the
Gumbel MDA; or iii) has a finite right endpoint; see Lemma 3.1 of
\cite{cltyy18}.

Since the continuity of the marginal distributions is not assumed, the tail probability evaluated at a high confidence level need not equal $(1-p)$             exactly; the following lemma shows that this discrepancy is asymptotically negligible in the Fr\'echet and Gumbel regimes as below.
\begin{lemma}[Tail probability at a high quantile]
\label{lemma:tail-at-var}
Let $X$ be a random variable with a distribution function $F$ and survival
function $\overline F$, and has an infinite right endpoint. Suppose that either i) $
F\in\operatorname{MDA}(\Phi_\alpha)
$ for some $\alpha>0$; or ii) $F\in\operatorname{MDA}(\Lambda)$.
Then
\[
\frac{
\overline F\bigl(\V@R_p(X)\bigr)
}{
1-p
}
\to1,
\quad p\uparrow1.
\]
\end{lemma}

\begin{proof}
Let $x_p:=\V@R_p(X)=F^{\leftarrow}(p)$. 
Since $F$ has an infinite right endpoint, $x_p\to\infty$ as $p\uparrow1$.
Write $
\overline F(x-):=1-F(x-)=\mathbb P(X\geq x)$.
By the definition of the generalized inverse, $
\overline F(x_p)\leq1-p\leq\overline F(x_p-)$.
It therefore suffices to prove that
\begin{equation}\label{eq:left-right-tail-equivalence}
\frac{\overline F(x-)}{\overline F(x)}
\rightarrow1,
\qquad x\to\infty.
\end{equation}
If $F\in\operatorname{MDA}(\Phi_\alpha)$, then
$\overline F\in\operatorname{RV}_{-\alpha}$, so that for every
$\delta\in(0,1)$ and all sufficiently large $x$,
\[
1
\leq
\frac{\overline F(x-)}{\overline F(x)}
\leq
\frac{\overline F((1-\delta)x)}{\overline F(x)}.
\]
Hence
$
1
\leq
\liminf_{x\to\infty}
\frac{\overline F(x-)}{\overline F(x)}
\leq
\limsup_{x\to\infty}
\frac{\overline F(x-)}{\overline F(x)}
\leq
(1-\delta)^{-\alpha}$;
letting $\delta\downarrow0$ proves
\eqref{eq:left-right-tail-equivalence}.

On the other hand, if $F\in\operatorname{MDA}(\Lambda)$, then by the standard
auxiliary-function characterization of the Gumbel MDA (also see Theorem 3.3.27 of \cite{Embrechts97}), there exists a
positive function $a$ such that, for every $t\in\mathbb R$,
\[
\frac{\overline F(x+a(x)t)}
     {\overline F(x)}
\rightarrow
e^{-t},
\qquad x\to\infty.
\]
For every $\delta>0$,
$
1
\leq
\frac{\overline F(x-)}{\overline F(x)}
\leq
\frac{\overline F(x-\delta a(x))}
     {\overline F(x)}$, and
consequently,
$
1
\leq
\liminf_{x\to\infty}
\frac{\overline F(x-)}{\overline F(x)}
\leq
\limsup_{x\to\infty}
\frac{\overline F(x-)}{\overline F(x)}
\leq
e^\delta$;
letting $\delta\downarrow0$ again proves
\eqref{eq:left-right-tail-equivalence}.
\end{proof}

\section{Asymptotic Aggregation Properties of Value-at-Risk}\label{sect:var_aa}
\subsection{Portfolios with Risks in Fr\'echet MDA}
\label{sec:frechet}

We first consider portfolios containing only risks from the Fr\'echet MDA. The analysis is
divided according to the relative order of tail indices of the Fr\'echet marginal risks.
We start our study with the case in which the relevant marginal risks share a common
Fr\'echet tail index, without imposing tail equivalence. We then consider
portfolios with different tail orders, for which the asymptotic behavior is
determined by one or several dominating Fr\'echet marginal risks with a common tail index.

\subsubsection{A Common Fr\'echet Tail Index}
\label{acrossFre}

We first suppose that all marginal risks belong to the same Fr\'echet
MDA, namely $\Phi_{\alpha}$ as defined in Section \ref{sec:evt-tail-dominance}.

\begin{assumption}[Common Fr\'echet tail index]
\label{Assumefrechet}
For each $i\in\{1,\ldots,d\}$, $
X_i\in\operatorname{MDA}(\Phi_\alpha), \alpha>0$.
\end{assumption}
Assumption \ref{Assumefrechet} does not imply tail equivalence of the
marginal risks. To see this, we
write
$\overline F_i(x)=x^{-\alpha}L_i(x)$
with $L_i$ slowly varying, for $i=1,\dots,d$. For any $i,j$, the ratio
$\frac{\overline F_i(x)}{\overline F_j(x)}=\frac{L_i(x)}{L_j(x)}$ need not converge as
$x\to\infty$. In this connection, we describe the relative
tail behavior along diverging sequences.

\begin{definition}[Sequential tail ratios]
\label{def:tailratio}
Let $\Theta:=\{\theta_k\}_{k\geq1}$ be a sequence with
$\theta_k\to\infty$. For $i,j\in\{1,\ldots,d\}$, define
\[
w_j^i(\Theta)
:=
\lim_{k\to\infty}
\frac{\overline F_i(\theta_k)}
     {\overline F_j(\theta_k)},
\]
whenever the limit exists in $[0,\infty]$.
\end{definition}

The quantity $w_j^i(\Theta)$ describes the asymptotic relative magnitude of the tail of
$X_i$ to that of $X_j$ along the sequence $\Theta$. In particular,
$w_j^i(\Theta)\in(0,\infty)$ means that the two tails are asymptotically
proportional along the sequence, even though their ratio need not converge
globally as $x\to\infty$.

For a fixed $j\in\{1,\ldots,d\}$, define
\begin{equation}\label{eq:capital_H}
\mathcal H_j(\mathbf{x};\Theta)
:=
\left(
\sum_{i=1}^d
x_i^{-\gamma\alpha}
\bigl(w_j^i(\Theta)\bigr)^{-\gamma}
\right)^{-\frac{1}{\gamma}}, \text{ for
 }\mathbf{x}=(x_1,\ldots,x_d)\in \mathbb{R}_+^d:=(0,\infty)^d,
\end{equation}
whose derivative admits the following expression by direct calculation:
\begin{equation}\label{eq:small_h}
\begin{aligned}
    h_j(\mathbf{x};\Theta)
&:=
\frac{\partial^d}
{\partial x_1\cdots\partial x_d}
\mathcal H_j(\mathbf{x};\Theta)\\
&=
\alpha^d
\prod_{r=0}^{d-1}(1+r\gamma)
\left\{
\prod_{i=1}^d
\bigl(w_j^i(\Theta)\bigr)^{-\gamma}
x_i^{-\gamma\alpha-1}
\right\}
\left(
\sum_{i=1}^d
x_i^{-\gamma\alpha}
\bigl(w_j^i(\Theta)\bigr)^{-\gamma}
\right)^{-\frac{1}{\gamma}-d}.
\end{aligned}
\end{equation}

\begin{theorem}[Sequential aggregate-tail asymptotics]\label{thm:Qexist}
Under Assumptions \ref{AssumeArchi} and \ref{Assumefrechet},
also suppose that $\Theta=\{\theta_k\}_{k=1}^{\infty}$ satisfies
$
w_j^i(\Theta)\in \mathbb{R}_+$ for any $i,j\in\{1,\ldots,d\}$.
Then, for every $j\in\{1,\ldots,d\}$, the limit
\[
Q_j(\gamma,\alpha;\Theta)
:=
\lim_{k\to\infty}
\frac{
\mathbb{P}\left(\sum_{i=1}^d X_i>\theta_k\right)
}{
\overline F_j(\theta_k)
}
\]
exists and takes the form:
\begin{equation}\label{qjint}
Q_j(\gamma,\alpha;\Theta)
=
\int_{\mathcal A}
h_j(\mathbf{x};\Theta) d\mathbf{x},
\qquad
\mathcal A
:=
\left\{
\mathbf{x}\in\mathbb{R}_+^d:
\sum_{i=1}^d\frac{1}{x_i}>1
\right\}.
\end{equation}
Moreover, $1\leq Q_j(\gamma,\alpha;\Theta)<\infty$.
\end{theorem}

\begin{remark}
Our present Theorem \ref{thm:Qexist} extends Proposition 2.2 of
\cite{enw08} from identically distributed marginals to non-identically
distributed Fr\'echet marginals yet with a common tail index. It also replaces global tail
equivalence, namely
$\frac{\overline F_i(x)}{\overline F_j(x)}
\rightarrow c_{ij}\in(0,\infty)$ as $x\to\infty$,
by asymptotic proportionality only along an exploding sequence $\Theta$. Particularly, our present 
theorem applies especially when the global marginal tail ratios fail to exist.
\end{remark}
For $\epsilon\in(0,1)$ and
$\mathbf{x}=(x_1,\ldots,x_d)\in\overline{\mathbb R}_+^d:=[0,\infty]^d$, write
$x_j^{(\epsilon)}:=x_j\wedge\epsilon^{-1}$, and define
\[
G_j^\epsilon(\mathbf{x};\Theta)
:=
\epsilon^\alpha
\left\{
\sum_{i\neq j}
x_i^{-\gamma\alpha}
\bigl(w_j^i(\Theta)\bigr)^{-\gamma}
+
\bigl(x_j^{(\epsilon)}\bigr)^{-\gamma\alpha}
\right\}^{-\frac{1}{\gamma}}.
\]
To prove Theorem \ref{thm:Qexist}, we need to first establish the following conditional weak convergence result, where we also adopt the conventions that 
$\frac{\theta}{0}=\infty$, $\infty^{-\gamma\alpha}=\frac{1}{\infty}=0$, $0^{-\gamma\alpha}=\frac{1}{0}=\infty$, and 
consider the random ratios $\frac{\theta_k}{X_i}$ taking values in $\overline{\mathbb R}_+$, and $\mathbbm{1}_A$ denotes the indicator function taking 1 for any argument taking values in the event set $A$.
\begin{lemma}\label{lemma:Gdist}
Under all the conditions and assumptions specified in Theorem \ref{thm:Qexist}, fix a $j\in\{1,\dots,d\}$, then given the event set $\{X_j>\epsilon\theta_k\}$, the conditional (confined) joint distribution of
\[
\left(
\frac{\theta_k}{X_1},
\ldots,
\frac{\theta_k}{X_d}
\right)
\;\bigg|\;
\{X_j>\epsilon\theta_k\}
\]
converges weakly to $G_j^\epsilon(\cdot;\Theta)$ on $\overline{\mathbb R}_+^d$, as $k\to\infty$; this limiting distribution has probability mass concentrated only on $\mathbb{R}_+^d$ and is
absolutely continuous there, with the density function:
\begin{equation}\label{eq:density_g}
g_j^\epsilon(\mathbf{y};\Theta)
=
\epsilon^\alpha
h_j(\mathbf{y};\Theta)
\mathbbm{1}_{\{y_j<\epsilon^{-1}\}},
\qquad
\mathbf{y}\in\mathbb{R}_+^d.
\end{equation}
\end{lemma}

\begin{proof}
Fix $\mathbf{y}\in\mathbb{R}_+^d$. Without loss of generality, the conditional joint distribution given $X_j>\epsilon\theta_k$ can be expressed by: 
\begin{align}
&\quad\mathbb P\left(
\frac{\theta_k}{X_i}<y_i,\ i=1,\ldots,d
\;\middle|\;
X_j>\epsilon\theta_k
\right)
=
\frac{
\mathbb P\left(\cap_{i\neq j}
\{X_i>\frac{\theta_k}{y_i}\}, X_j>\max\left(\frac{\theta_k}{y_j}, \frac{\theta_k}{\epsilon^{-1}}\right)
\right)
}{
\overline F_j(\epsilon\theta_k)
}
\nonumber\\
&=
\frac{
\mathbb P\left(
X_i>\frac{\theta_k}{\widetilde y_i},\ i=1,\ldots,d
\right)
}{
\overline F_j(\epsilon\theta_k)
}=
\frac{
\psi^{-1}\left(
\sum_{i=1}^d
\psi\bigl(\overline F_i\left(\frac{\theta_k}{\widetilde y_i}\right)\bigr)
\right)
}{
\overline F_j(\epsilon\theta_k)
},
\label{eq:Gdist1}
\end{align}
where $\widetilde y_i:=y_i$ for $i\in\{1,\dots,d\}\backslash\{j\}$, and $\widetilde y_j:=y_j^{(\epsilon)}$; the last equality holds due to the immediate definition of the generator $\psi$ of Archimedean copula. 
By Assumption \ref{Assumefrechet} and Definition \ref{def:tailratio},  as $k\to\infty$, $\theta_k\to\infty$ and also
\begin{equation}
\label{eq:tail_ratio_limit}
\frac{
\overline F_i\left(\frac{\theta_k}{\widetilde y_i}\right)
}{
\overline F_j(\theta_k)
}
=
\frac{
\overline F_i\left(\frac{\theta_k}{\widetilde y_i}\right)
}{
\overline F_i(\theta_k)
}
\frac{
\overline F_i(\theta_k)
}{
\overline F_j(\theta_k)
}
\rightarrow
{(\widetilde y_i)}^\alpha w_j^i(\Theta),\quad\forall i=1,\dots,d.
\end{equation}
Since $\psi\in\operatorname{RV}_{-\gamma}$ at zero by Assumption \ref{AssumeArchi}, then the uniform convergence theorem for
regularly varying functions (see Theorem 1.5.2 of \cite{bgt87}) together with \eqref{eq:tail_ratio_limit} yields:
\[
\frac{
\psi\bigl(\overline F_i\left(\frac{\theta_k}{\widetilde y_i}\right)\bigr)
}{
\psi\bigl(\overline F_j(\theta_k)\bigr)
}
\rightarrow
\left(
(\widetilde y_i)^\alpha w_j^i(\Theta)
\right)^{-\gamma},\quad\forall i=1,\dots,d,
\]
and hence
\begin{equation}\label{eq:sum_ratio}
\frac{
\sum_{i=1}^d
\psi\bigl(\overline F_i\left(\frac{\theta_k}{\widetilde y_i}\right)\bigr)
}{
\psi\bigl(\overline F_j(\theta_k)\bigr)
}
\rightarrow
\sum_{i=1}^d
{(\widetilde y_i)}^{-\gamma\alpha}
\bigl(w_j^i(\Theta)\bigr)^{-\gamma}.
\end{equation}
By the asymptotic inverse theorem (Theorem 1.5.12 of \cite{bgt87}),
$\psi^{-1}\in\operatorname{RV}_{-\frac{1}{\gamma}}$ at infinity. By Assumption \ref{AssumeArchi}, the strictness of $\psi$
warrants
$\psi^{-1}\bigl(\psi(\overline F_j(\theta_k))\bigr)
=
\overline F_j(\theta_k)$,
and together with \eqref{eq:sum_ratio},
\begin{equation}
\frac{
\psi^{-1}\left(
\sum_{i=1}^d
\psi\bigl(\overline F_i\left(\frac{\theta_k}{\widetilde y_i}\right)\bigr)
\right)
}{
\overline F_j(\theta_k)
}=\frac{
\psi^{-1}\left(
\sum_{i=1}^d
\psi\bigl(\overline F_i\left(\frac{\theta_k}{\widetilde y_i}\right)\bigr)
\right)
}{
\psi^{-1}(\psi(\overline F_j(\theta_k)))
}
\rightarrow
\left\{
\sum_{i=1}^d
{(\widetilde y_i)}^{-\gamma\alpha}
\bigl(w_j^i(\Theta)\bigr)^{-\gamma}
\right\}^{-\frac{1}{\gamma}}.
\label{eq:Gdist2}
\end{equation}
On the other hand, by regular variation of $\overline F_j$,
\begin{equation}\label{eq:Gdist3}
\frac{\overline F_j(\epsilon\theta_k)}
     {\overline F_j(\theta_k)}
\rightarrow
\epsilon^{-\alpha}\text{ as }k\to\infty.
\end{equation}
Finally, combining \eqref{eq:Gdist1}, \eqref{eq:Gdist2} and \eqref{eq:Gdist3}, we see that as $k\to\infty$,
\[
\mathbb P\left(
\frac{\theta_k}{X_i}<y_i,\ i=1,\ldots,d
\;\middle|\;
X_j>\epsilon\theta_k
\right)=\frac{
\psi^{-1}\left(
\sum_{i=1}^d
\psi\bigl(\overline F_i\left(\frac{\theta_k}{\widetilde y_i}\right)\bigr)
\right)
}{
\overline F_j(\theta_k)
}\cdot\left(\frac{\overline F_j(\epsilon\theta_k)}{\overline F_j(\theta_k)}\right)^{-1}
\rightarrow
G_j^\epsilon(\mathbf{y};\Theta).
\]
To see that $G_j^\epsilon$ is a proper distribution function indeed,
recall that $h_j(\mathbf{y};\Theta)\geq0$ on $\mathbb{R}_+^d$ by  \eqref{eq:small_h}, and 
$\mathcal H_j$ in \eqref{eq:capital_H} is smooth on $\mathbb{R}_+^d$ and vanishes whenever any
one of its coordinates decreases to zero, since the corresponding summand and hence the whole summation become infinity. Obviously, by the iterative use of fundamental theorem of calculus,
\[
\mathcal H_j(\mathbf{y};\Theta)
=
\int_0^{y_1}\cdots\int_0^{y_d}
h_j(\mathbf{u};\Theta) d\mathbf{u}.
\]
Since $w_j^j(\Theta)=1$, it follows that
\begin{align*}
G_j^\epsilon(\mathbf{y};\Theta)
=
\epsilon^\alpha
\mathcal H_j
\left(
\widetilde{y}_1,\dots,\widetilde{y}_d;\Theta
\right)=
\int_0^{y_1}\cdots\int_0^{y_d}
\epsilon^\alpha h_j(\mathbf{u};\Theta)
\mathbbm{1}_{\{u_j<\epsilon^{-1}\}}
 d\mathbf{u};
\end{align*}
and in view of \eqref{eq:capital_H}, we further see that:
\[
\lim_{y_1,\ldots,y_d\to\infty}
G_j^\epsilon(\mathbf{y};\Theta)
=\lim_{y_1,\ldots,y_d\to\infty}\epsilon^\alpha
\mathcal H_j
\left(
\widetilde{y}_1,\dots,\widetilde{y}_d;\Theta
\right)=
\epsilon^\alpha
\left\{
(\epsilon^{-1})^{-\gamma\alpha}
\right\}^{-\frac{1}{\gamma}}
=
1,
\]
thus $G_j^\epsilon$ is a proper distribution function on
$\mathbb{R}_+^d$, and is absolutely continuous there with density $
g_j^\epsilon(\mathbf{y};\Theta)
=
\epsilon^\alpha
h_j(\mathbf{y};\Theta)
\mathbbm{1}_{\{y_j<\epsilon^{-1}\}}$.
\end{proof}
\begin{proof}[Proof of Theorem \ref{thm:Qexist}]
Here $j\in\{1,\ldots,d\}$ is fixed. Denote $\mathbf Y_k:=\left(\frac{\theta_k}{X_1},\ldots,\frac{\theta_k}{X_d}\right)$ and
$\mathcal A^*
:=\{\mathbf y\in\overline{\mathbb R}_+^d:
\sum_{i=1}^d y_i^{-1}>1\}$, so that
$\{S_d=\sum_{i=1}^d X_i>\theta_k\}=\{\mathbf Y_k\in\mathcal A^*\}$.
We first let $\epsilon\in\left(0,\frac{1}{d}\right)\subset(0,1)$ be arbitrary, and consider two subcases depending on whether $X_j>\epsilon\theta_k$.
\begin{enumerate}
    \item[i)] Case of $S_d>\theta_k$ and $X_j>\epsilon\theta_k$: 
By Lemma \ref{lemma:Gdist}, conditional on the event
$\{X_j>\epsilon\theta_k\}$, $\mathbf Y_k$ converges weakly to a random vector 
$\mathbf Y^\epsilon$ whose distribution is concentrated on
$\mathbb{R}_+^d$ and has the density
$g_j^\epsilon(\mathbf y;\Theta)$ given in \eqref{eq:density_g}.
Within $\mathbb{R}_+^d$, the boundary of $\mathcal A^*$ is
$\partial\mathcal A^*:=\{\mathbf x:\sum_{i=1}^d x_i^{-1}=1\}$, with a $d$-dimensional Lebesgue measure of zero, hence
$\mathbb P(\mathbf Y^\epsilon\in\partial\mathcal A^*)=0$.
Then by Portmanteau theorem
(Theorem 5.25(iv) of \cite{kallenberg2021}), we have, as $k\to\infty$:
\begin{align*}
    &\ \mathbb P(S_d>\theta_k\mid X_j>\epsilon\theta_k)=\mathbb{P}(\mathbf Y_k\in\mathcal A^*\mid X_j>\epsilon\theta_k)\to\mathbb{P}(\mathbf Y^\epsilon\in\mathcal A^*)\\
=&\ \mathbb{P}(\mathbf Y^\epsilon\in\mathcal A)=\int_{\mathcal{A}}g_j^\epsilon(\mathbf y;\Theta)d\mathbf y=\int_{\mathcal{A}}\epsilon^\alpha h_j(\mathbf y;\Theta)\mathbbm{1}_{\{y_j<\epsilon^{-1}\}}d\mathbf y=\epsilon^\alpha\int_{\mathcal A_\epsilon}h_j(\mathbf y;\Theta) d\mathbf y=:\epsilon^\alpha\mathcal L_j(\epsilon;\Theta),
\end{align*}
where $\mathcal A_\epsilon
:=\{\mathbf y\in\mathcal A:y_j<\epsilon^{-1}\}$.
Together with \eqref{eq:Gdist3}, we deduce, as $k\to\infty$, along the definite sequence $\Theta=\{\theta_k\}_{k\geq1}$,
\begin{equation}\label{eq:Q-upper-first}
\frac{\mathbb P(S_d>\theta_k,\ X_j>\epsilon\theta_k)}
     {\overline F_j(\theta_k)}=\mathbb P(S_d>\theta_k\mid X_j>\epsilon\theta_k)\cdot\frac{\overline F_j(\epsilon\theta_k)}{\overline F_j(\theta_k)}
\to
\mathcal L_j(\epsilon;\Theta).
\end{equation}
Particularly, 
\begin{equation}\label{eq:liminf_bound}
\liminf_{k\to\infty}
\frac{\mathbb P(S_d>\theta_k)}{\overline F_j(\theta_k)}\geq \liminf_{k\to\infty}
\frac{\mathbb P(S_d>\theta_k,\ X_j>\epsilon\theta_k)}{\overline F_j(\theta_k)}\geq \mathcal L_j(\epsilon;\Theta).
\end{equation}
\item[ii)] Case of $S_d>\theta_k$ and $X_j\leq\epsilon\theta_k$: Note that by pigeonhole principle, there must exist an $i\neq j$ such that $X_i>\frac{\theta_k}{d}$, hence
\begin{equation}\label{eq:Q-remainder-union}
\mathbb{P}(S_d>\theta_k,\ X_j\leq\epsilon\theta_k)
\leq
\sum_{i\neq j}
\mathbb{P}\left(
X_i>\frac{\theta_k}{d},\
X_j\leq\epsilon\theta_k
\right)=
\sum_{i\neq j}
\left[
\overline F_i\left(\frac{\theta_k}{d}\right)
-
\mathbb{P}\left(
X_i>\frac{\theta_k}{d},\
X_j>\epsilon\theta_k
\right)
\right]. 
\end{equation}
By Assumption \ref{Assumefrechet} and
Definition \ref{def:tailratio},
$\frac{\overline F_i\left(\frac{\theta_k}{d}\right)}{\overline F_j(\theta_k)}
\to d^\alpha w_j^i(\Theta)=:c_{ji}$ as 
$k\to\infty$ for any $i=1,\dots,d$.
Furthermore, by Assumption \ref{AssumeArchi}, for any $i\neq j$, the bivariate survival
copula of $(X_i,X_j)$ also has the generator $\psi$ by the very definition; also see Section 4.6 of \cite{n06}. Therefore, by applying the uniform convergence theorem and the asymptotic
inverse theorem (Theorems 1.5.2 and 1.5.12 in \cite{bgt87}, respectively) in a similar manner as obtaining \eqref{eq:Gdist2}, we have, as $k\to\infty$:
$$
\frac{\mathbb P(X_i>\theta_k/d,\ X_j>\epsilon\theta_k)}
     {\overline F_j(\theta_k)}=
\frac{
\psi^{-1}\left(
\psi\bigl(\overline F_i\left(\frac{\theta_k}{d}\right)\bigr)
+
\psi\bigl(\overline F_j(\epsilon\theta_k)\bigr)
\right)
}{
\psi^{-1}\bigl(\psi(\overline F_j(\theta_k))\bigr)
}\to
\left(c_{ji}^{-\gamma}+\epsilon^{\gamma\alpha}\right)^{-\frac{1}{\gamma}}.
$$
Consequently, \eqref{eq:Q-remainder-union} gives:
\begin{align*}
    &\ \frac{\mathbb P(S_d>\theta_k,\ X_j\leq\epsilon\theta_k)}
     {\overline F_j(\theta_k)}\leq\sum_{i\neq j}\left(\frac{\overline F_i\left(\frac{\theta_k}{d}\right)}
     {\overline F_j(\theta_k)}-\frac{\mathbb P\left(X_i>\frac{\theta_k}{d}, X_j>\epsilon\theta_k\right)}
     {\overline F_j(\theta_k)}\right)\\
     \to&\ \sum_{i\neq j}\left(c_{ji}-\left(c_{ji}^{-\gamma}+\epsilon^{\gamma\alpha}\right)^{-\frac{1}{\gamma}}\right)=:r_j(\epsilon;\Theta), \text{ as }k\to\infty,
\end{align*}
hence
\begin{equation}\label{eq:remainder-bound}
\limsup_{k\to\infty}
\frac{\mathbb P(S_d>\theta_k,\ X_j\leq\epsilon\theta_k)}
     {\overline F_j(\theta_k)}
     \leq r_j(\epsilon;\Theta)
.
\end{equation}
Combining \eqref{eq:Q-upper-first} and
\eqref{eq:remainder-bound}, we conclude that
\begin{equation}\label{eq:Q_upper_final}
\begin{aligned}
    \limsup_{k\to\infty}
\frac{\mathbb P(S_d>\theta_k)}{\overline F_j(\theta_k)}
&=\limsup_{k\to\infty}\left(\frac{\mathbb P(S_d>\theta_k,\ X_j>\epsilon\theta_k)}
     {\overline F_j(\theta_k)}+\frac{\mathbb P(S_d>\theta_k,\ X_j\leq\epsilon\theta_k)}
     {\overline F_j(\theta_k)}\right)\\
     &\leq
\mathcal L_j(\epsilon;\Theta)+r_j(\epsilon;\Theta).
\end{aligned}
\end{equation}
\end{enumerate}
Now let $\epsilon\downarrow0$. We first note that $r_j(\epsilon;\Theta)\to0$ since $c_{ji}$ is clearly positive and finite for any $i,j=1,\dots,d$ due to the very condition in the theorem statement that $w_j^i(\Theta)\in\mathbb{R}_+$. Furthermore, $\mathcal A_\epsilon\uparrow\mathcal A$, and recall \eqref{eq:small_h} that $h_j(\mathbf x;\Theta)\geq0$, hence we invoke monotone convergence theorem for Lebesgue integration to see that:
\begin{equation}\label{eq:L0}
\mathcal L_j(\epsilon;\Theta)
\uparrow
\int_{\mathcal A}h_j(\mathbf x;\Theta) d\mathbf x,\text{ as }\epsilon\downarrow0.
\end{equation}
Therefore, combining \eqref{eq:liminf_bound} and \eqref{eq:Q_upper_final}, and in light of \eqref{eq:L0}, we have:
\[
\lim_{k\to\infty}
\frac{\mathbb P(S_d>\theta_k)}{\overline F_j(\theta_k)}=
\int_{\mathcal A}h_j(\mathbf x;\Theta) d\mathbf x=:Q_j(\gamma,\alpha;\Theta).
\]
In particular, by the nonnegativity of $X_i$'s and another application of pigeonhole principle,
$\overline F_j(\theta_k)\leq\mathbb P(S_d>\theta_k)
\leq\sum_{i=1}^d\overline F_i\left(\frac{\theta_k}{d}\right)$, hence:
\begin{equation}\label{eq:elementary_bound}
1
\leq
\liminf_{k\to\infty}
\frac{\mathbb P(S_d>\theta_k)}{\overline F_j(\theta_k)}
\leq
\limsup_{k\to\infty}
\frac{\mathbb P(S_d>\theta_k)}{\overline F_j(\theta_k)}
\leq
\sum_{i=1}^d c_{ji}
<\infty,
\end{equation}
henceforth, we conclude that
$1\leq Q_j(\gamma,\alpha;\Theta)<\infty$.
\end{proof}

The following result extends Theorem 2.3 of \cite{enw08} to asymptotically proportional, non-identically distributed marginals and particularly it is strengthened by establishing strict monotonicity for $\alpha\neq1$. Indeed, when the marginals are identical, $w_j^i(\Theta)=1$ for all $i,j$, and $Q_j(\gamma,\alpha;\Theta)$ reduces to the limiting constant considered in \cite{enw08}. As noted in the remark following Theorem 2.3 of \cite{enw08}, their argument does not establish strict monotonicity when $d\geq3$.  We here give two different proofs: the first one is a novel purely analytic proof, while the other is based on a convex-order comparison and a coordinate-wise replacement argument.

\begin{theorem}[Strict monotonicity in $\gamma$ of $Q_j(\gamma,\alpha;\Theta)$]\label{thm:qin}
Under the assumptions and conditions specified in Theorem \ref{thm:Qexist}, for any fixed
$j\in\{1,\ldots,d\}$, $Q_j(\gamma,\alpha;\Theta)$ is:
\begin{enumerate}[label=(\roman*)]
    \item  strictly decreasing in $\gamma$, if $0<\alpha<1$; or
    \item identically equal to the constant $\sum_{i=1}^d w_j^i(\Theta)$, if $\alpha=1$; or
    \item strictly increasing in $\gamma$, if $\alpha>1$.
\end{enumerate}
\end{theorem}
To facilitate our further discussion, we first establish the following useful probabilistic representation of $Q_j(\gamma,\alpha;\Theta)$:

\begin{lemma}
\label{lemma:Q-exponential}
Under the assumptions and conditions specified in Theorem \ref{thm:Qexist},
\begin{equation}\label{eq:Q-V-rho}
Q_j(\gamma,\alpha;\Theta)
=\frac{1}{\Gamma\left(\frac{1}{\gamma}+1\right)}
\mathbb E\left[
\left(
\sum_{i=1}^d \left(w_j^i(\Theta)\right)^{\frac{1}{\alpha}} E_i^{\frac{1}{\gamma\alpha}}
\right)^\alpha
\right],
\end{equation}
where $E_1,\ldots,E_d$ are independent standard exponential
random variables.
\end{lemma}

\begin{proof}
We first introduce a change of variables,
$t_i:=x_i^{-\gamma\alpha}(w_j^i(\Theta))^{-\gamma}$ for $i=1,\dots,d$, then $x_i=(w_j^i(\Theta))^{-\frac{1}{\alpha}}t_i^{-\frac{1}{\gamma\alpha}}$, so that the integral constraint specified in $\mathcal{A}$ becomes:
$$\sum_{i=1}^d\frac{1}{x_i}=\sum_{i=1}^d(w_j^i(\Theta))^{\frac{1}{\alpha}}t_i^{\frac{1}{\gamma\alpha}}>1;$$
meanwhile, the Jacobian of the transformation is:
$$\left|\det\left(\frac{\partial(x_1,\dots,x_d)}{\partial(t_1,\dots,t_d)}\right)\right|=\left|\prod_{i=1}^d\left(-\frac{1}{\gamma\alpha}\right)t_i^{-\frac{1}{\gamma\alpha}-1}(w_j^i(\Theta))^{-\frac{1}{\alpha}}\right|= \left(\frac{1}{\gamma\alpha}\right)^d \left(\prod_{i=1}^d t_i\right)^{-\frac{1}{\gamma\alpha}-1}\left(\prod_{i=1}^d w_j^i(\Theta)\right)^{-\frac{1}{\alpha}}.$$
Therefore, \eqref{qjint} becomes:
\begin{equation}\label{eq:Q-t-representation}
\begin{aligned}
    &\ Q_j(\gamma,\alpha;\Theta)\\
=&\ \int_{\sum_{i=1}^d(w_j^i(\Theta))^{\frac{1}{\alpha}}t_i^{\frac{1}{\gamma\alpha}}>1}\alpha^d
\prod_{r=0}^{d-1}(1+r\gamma)
\left\{
\prod_{i=1}^d
\bigl(w_j^i(\Theta)\bigr)^{-\gamma}
\left(t_i\bigl(w_j^i(\Theta)\bigr)^{\gamma}\right)^{-\frac{1}{\gamma\alpha}(-\gamma\alpha-1)}
\right\}
\left(
\sum_{i=1}^d
t_i\right)^{-\frac{1}{\gamma}-d}\\
&\qquad\cdot \left(\frac{1}{\gamma\alpha}\right)^d \left(\prod_{i=1}^d t_i\right)^{-\frac{1}{\gamma\alpha}-1}\left(\prod_{i=1}^d w_j^i(\Theta)\right)^{-\frac{1}{\alpha}}dt_1\cdots dt_d\\    
=&\ \prod_{r=0}^{d-1}\left(\frac{1}{\gamma}+r\right)\int_{\sum_{i=1}^d(w_j^i(\Theta))^{\frac{1}{\alpha}}t_i^{\frac{1}{\gamma\alpha}}>1}
\left(\sum_{i=1}^d t_i\right)^{-\frac{1}{\gamma}-d}
dt_1\cdots dt_d\\
=&\ 
\frac{\Gamma\left(\frac{1}{\gamma}+d\right)}{\Gamma\left(\frac{1}{\gamma}\right)}
\int_{\mathbb{R}_+^d}
\mathbbm{1}_{\left\{\sum_{i=1}^d(w_j^i(\Theta))^{\frac{1}{\alpha}}t_i^{\frac{1}{\gamma\alpha}}>1\right\}}
\left(\sum_{i=1}^d t_i\right)^{-\frac{1}{\gamma}-d}
dt_1\cdots dt_d.
\end{aligned}
\end{equation}
We next apply another change of variables:
\begin{equation}\label{eq:var_transformation}
(t_1,\dots,t_d)\mapsto(u_1,\dots,u_{d-1},r):=\left(\frac{t_1}{\sum_{i=1}^d t_i},\dots,\frac{t_{d-1}}{\sum_{i=1}^d t_i},{\sum_{i=1}^d t_i}\right),
\end{equation}
whose Jacobian is clearly $r^{d-1}$; for the sake of convenience, we also denote $u_d:=1-\sum_{i=1}^{d-1}u_i=\frac{t_d}{\sum_{i=1}^d t_i}$. Henceforth, \eqref{eq:Q-t-representation} can be rewritten as:
\begin{equation}\label{eq:Q-dirichlet-integral}
\begin{aligned}
Q_j(\gamma,\alpha;\Theta)
&=
\frac{\Gamma\left(\frac{1}{\gamma}+d\right)}{\Gamma\left(\frac{1}{\gamma}\right)}\int_0^1\int_0^{1-u_1}\cdots\int_0^{1-u_1-\dots-u_{d-2}}\int_{(\sum_{i=1}^d(w_j^i(\Theta))^{\frac{1}{\alpha}}u_i^{\frac{1}{\gamma\alpha}})^{-\gamma\alpha}}^\infty
r^{-\frac{1}{\gamma}-d}\cdot r^{d-1}drdu_{d-1}\cdots du_1\\
&=
\frac{\Gamma\left(\frac{1}{\gamma}+d\right)}{\Gamma\left(\frac{1}{\gamma}+1\right)}
\int_0^1\int_0^{1-u_1}\cdots\int_0^{1-u_1-\dots-u_{d-2}}
\left(\sum_{i=1}^d(w_j^i(\Theta))^{\frac{1}{\alpha}}u_i^{\frac{1}{\gamma\alpha}}\right)^\alpha du_{d-1}\cdots du_1.
\end{aligned}
\end{equation}
Next, we consider the expectation in the RHS of \eqref{eq:Q-V-rho}, which can be written as: 
$$\mathbb E\left[
\left(
\sum_{i=1}^d \left(w_j^i(\Theta)\right)^{\frac{1}{\alpha}} E_i^{\frac{1}{\gamma\alpha}}
\right)^\alpha
\right]=\int_{\mathbb{R}_+^d}\left(
\sum_{i=1}^d \left(w_j^i(\Theta)\right)^{\frac{1}{\alpha}} e_i^{\frac{1}{\gamma\alpha}}
\right)^\alpha e^{-\sum_{i=1}^de_i}de_1\cdots de_d;$$
inspired by R\'enyi representation (see \cite{renyi1953theory} and Chapter 2 of \cite{david2004order}), we adopt the following change of variable analogous to \eqref{eq:var_transformation}:
$$(e_1,\dots,e_d)\mapsto(\tilde{u}_1,\dots,\tilde{u}_{d-1},\tilde{r}):=\left(\frac{e_1}{\sum_{i=1}^d e_i},\dots,\frac{e_{d-1}}{\sum_{i=1}^d e_i},{\sum_{i=1}^d e_i}\right),$$
whose Jacobian is $\tilde{r}^{d-1}$; similarly, we also denote $\tilde{u}_d:=1-\sum_{i=1}^{d-1}\tilde{u}_i=\frac{e_d}{\sum_{i=1}^d e_i}$. Then we have:
\begin{equation}\label{eq:Q-dirichlet-integral-RHS}
\begin{aligned}
&\quad \mathbb E\left[
\left(
\sum_{i=1}^d \left(w_j^i(\Theta)\right)^{\frac{1}{\alpha}} E_i^{\frac{1}{\gamma\alpha}}
\right)^\alpha
\right]\\
&=\int_0^1\int_0^{1-\tilde{u}_1}\cdots\int_0^{1-\tilde{u}_1-\dots-\tilde{u}_{d-2}}\int_0^{\infty}\left(
\sum_{i=1}^d \left(w_j^i(\Theta)\right)^{\frac{1}{\alpha}} (\tilde{u}_i\tilde{r})^{\frac{1}{\gamma\alpha}}
\right)^\alpha e^{-\tilde{r}}\tilde{r}^{d-1}d\tilde{r}d\tilde{u}_{d-1}\cdots d\tilde{u}_1\\
&=\int_0^1\int_0^{1-\tilde{u}_1}\cdots\int_0^{1-\tilde{u}_1-\dots-\tilde{u}_{d-2}}\left(
\sum_{i=1}^d \left(w_j^i(\Theta)\right)^{\frac{1}{\alpha}} \tilde{u}_i^{\frac{1}{\gamma\alpha}}
\right)^\alpha \int_0^{\infty} \tilde{r}^{\frac{1}{\gamma}+d-1}e^{-\tilde{r}}d\tilde{r} d\tilde{u}_{d-1}\cdots d\tilde{u}_1\\
&=
\Gamma\left(\frac{1}{\gamma}+d\right)
\int_0^1\int_0^{1-\tilde{u}_1}\cdots\int_0^{1-\tilde{u}_1-\dots-\tilde{u}_{d-2}}
\left(\sum_{i=1}^d(w_j^i(\Theta))^{\frac{1}{\alpha}}\tilde{u}_i^{\frac{1}{\gamma\alpha}}\right)^\alpha d\tilde{u}_{d-1}\cdots d\tilde{u}_1.
\end{aligned}
\end{equation}
which is exactly \eqref{eq:Q-dirichlet-integral} up to the scaling factor of $\frac{1}{\Gamma\left(\frac{1}{\gamma}+1\right)}$, hence \eqref{eq:Q-V-rho} is established.
\end{proof}

We now present the first novel purely analytic proof of Theorem \ref{thm:qin}:

\begin{proof}[First Proof of Theorem \ref{thm:qin}]
Fix $j\in\{1,\ldots,d\}$, and also fix the sequence $\Theta=\{\theta_k\}_{k=1}^{\infty}$ with $\theta_k\to\infty$ that satisfies
$
w_j^i(\Theta)\in \mathbb{R}_+$ for any $i,j\in\{1,\ldots,d\}$ as specified in Theorem \ref{thm:Qexist}. The case of $\alpha=1$ is trivial, since $\mathbb{E}\left(E_i^{\frac{1}{\gamma}}\right)=\int_0^{\infty}t^{\frac{1}{\gamma}}e^{-t}dt=\Gamma\left(\frac{1}{\gamma}+1\right)$ for $i=1,\dots,d$, henceforth: 
\begin{equation}\label{eq:alpha_equals_1}
Q_j(\gamma,1;\Theta)=\frac{1}{\Gamma\left(\frac{1}{\gamma}+1\right)}
\mathbb E
\left(
\sum_{i=1}^d w_j^i(\Theta) E_i^{\frac{1}{\gamma}}
\right)=\frac{1}{\Gamma\left(\frac{1}{\gamma}+1\right)}\sum_{i=1}^dw_j^i(\Theta)\Gamma\left(\frac{1}{\gamma}+1\right)=\sum_{i=1}^dw_j^i(\Theta).
\end{equation}
It then suffices to study the case of $\alpha\neq1$; note that standard measure-theoretical arguments in connection with dominated convergence theorem allow the interchange of the order between differentiation and expectation-taking, and so:
\begin{equation}\label{eq:Q-analytic-derivative}
\begin{aligned}
\frac{\partial}{\partial\gamma}Q_j(\gamma,\alpha;\Theta)
&=\frac{-\Gamma'\left(\frac{1}{\gamma}+1\right)\cdot\left(-\frac{1}{\gamma^2}\right)}{\left(\Gamma\left(\frac{1}{\gamma}+1\right)\right)^2}\mathbb E\left[
\left(
\sum_{i=1}^d \left(w_j^i(\Theta)\right)^{\frac{1}{\alpha}} E_i^{\frac{1}{\gamma\alpha}}
\right)^\alpha
\right]\\
&\qquad+\frac{1}{\Gamma\left(\frac{1}{\gamma}+1\right)}\mathbb E\left[
\alpha \left(
\sum_{i=1}^d \left(w_j^i(\Theta)\right)^{\frac{1}{\alpha}} E_i^{\frac{1}{\gamma\alpha}}
\right)^{\alpha-1}\sum_{i=1}^d \left(w_j^i(\Theta)\right)^{\frac{1}{\alpha}} E_i^{\frac{1}{\gamma\alpha}}(\ln E_i)\left(-\frac{1}{\gamma^2\alpha}\right)
\right]
\\
&=-\frac{1}{\gamma^2}
\sum_{i=1}^d\frac{\left(w_j^i(\Theta)\right)^{\frac{1}{\alpha}}}{\Gamma\left(\frac{1}{\gamma}+1\right)}\mathbb E\left[
\left(
\sum_{l=1}^d \left(w_j^l(\Theta)\right)^{\frac{1}{\alpha}} E_l^{\frac{1}{\gamma\alpha}}
\right)^{\alpha-1}E_i^{\frac{1}{\gamma\alpha}}\left(\ln E_i-\Psi\left(\frac{1}{\gamma}+1\right)\right)
\right],
\end{aligned}
\end{equation}
where $\Psi(x):=\frac{\Gamma'(x)}{\Gamma(x)}$ is the digamma function.

For any $i=1,\dots,d$, we rewrite the summand as follows by an application of tower property on the random sub-sum $\sum_{l\neq i} \left(w_j^l(\Theta)\right)^{\frac{1}{\alpha}} E_l^{\frac{1}{\gamma\alpha}}$, which is independent of $E_i$ by construction:
\begin{equation}\label{eq:tower_property}
\begin{aligned}
    &\ \frac{\left(w_j^i(\Theta)\right)^{\frac{1}{\alpha}}}{\Gamma\left(\frac{1}{\gamma}+1\right)}\mathbb E\left[
\left(
\sum_{l=1}^d \left(w_j^l(\Theta)\right)^{\frac{1}{\alpha}} E_l^{\frac{1}{\gamma\alpha}}
\right)^{\alpha-1}E_i^{\frac{1}{\gamma\alpha}}\left(\ln E_i-\Psi\left(\frac{1}{\gamma}+1\right)\right)
\right]\\
=&\ \frac{\left(w_j^i(\Theta)\right)^{\frac{1}{\alpha}}}{\Gamma\left(\frac{1}{\gamma}+1\right)}\mathbb E\left[\left.\mathbb E\left[
\left(
\sum_{l=1}^d \left(w_j^l(\Theta)\right)^{\frac{1}{\alpha}} E_l^{\frac{1}{\gamma\alpha}}
\right)^{\alpha-1}E_i^{\frac{1}{\gamma\alpha}}\left(\ln E_i-\Psi\left(\frac{1}{\gamma}+1\right)\right)
\right]\right|\sum_{l\neq i} \left(w_j^l(\Theta)\right)^{\frac{1}{\alpha}} E_l^{\frac{1}{\gamma\alpha}}\right]\\
=&\ \frac{\left(w_j^i(\Theta)\right)^{\frac{1}{\alpha}}}{\Gamma\left(\frac{1}{\gamma}+1\right)}\mathbb E\left[\int_0^{\infty}
\left(\left(w_j^i(\Theta)\right)^{\frac{1}{\alpha}} t^{\frac{1}{\gamma\alpha}}+\sum_{l\neq i} \left(w_j^l(\Theta)\right)^{\frac{1}{\alpha}} E_l^{\frac{1}{\gamma\alpha}}
\right)^{\alpha-1}t^{\frac{1}{\gamma\alpha}}\left(\ln t-\Psi\left(\frac{1}{\gamma}+1\right)\right)e^{-t}dt
\right]\\
=:&\ \left(w_j^i(\Theta)\right)^{\frac{1}{\alpha}}\mathbb{E}\left[\int_0^{\infty}h_{E_{-i}}(t)g_\gamma(t)dt\right],
\end{aligned}
\end{equation}
where we denote
$$h_{E_{-i}}(t):=\left(\left(w_j^i(\Theta)\right)^{\frac{1}{\alpha}} +\sum_{l\neq i} \left(w_j^l(\Theta)\right)^{\frac{1}{\alpha}} \left(\frac{E_l}{t}\right)^{\frac{1}{\gamma\alpha}}
\right)^{\alpha-1},\quad g_\gamma(t):=\left(\ln t-\Psi\left(\frac{1}{\gamma}+1\right)\right)\frac{t^{\frac{1}{\gamma}}e^{-t}}{\Gamma\left(\frac{1}{\gamma}+1\right)}.$$
It is clear that $\sum_{l\neq i} \left(w_j^l(\Theta)\right)^{\frac{1}{\alpha}} \left(\frac{E_l}{t}\right)^{\frac{1}{\gamma\alpha}}$ is strictly decreasing in $t$, hence $h_{E_{-i}}(t)$ is strictly increasing if $0<\alpha<1$, but strictly decreasing if $\alpha>1$. On the other hand, note that $\frac{t^{\frac{1}{\gamma}}e^{-t}}{\Gamma\left(\frac{1}{\gamma}+1\right)}$ is the density of $\text{Gamma}\left(\frac{1}{\gamma}+1,1\right)$, and we have:
$$g_\gamma(t)\left\{\begin{aligned}
    &<0,\text{ if } t<t_\gamma:=e^{\Psi\left(\frac{1}{\gamma}+1\right)};\\
    &=0,\text{ if } t=t_\gamma;\\
    &>0,\text{ if } t>t_\gamma.
\end{aligned}\right.$$
furthermore, by the definition of digamma function, $\Psi\left(\frac{1}{\gamma}+1\right)=\int_0^{\infty} \frac{t^{\frac{1}{\gamma}}e^{-t}}{\Gamma\left(\frac{1}{\gamma}+1\right)}\ln tdt$, and as the integral of density equals 1, $\Psi\left(\frac{1}{\gamma}+1\right)\left(\int_0^{\infty}\frac{t^{\frac{1}{\gamma}}e^{-t}}{\Gamma\left(\frac{1}{\gamma}+1\right)}dt\right)=\int_0^{\infty} \frac{t^{\frac{1}{\gamma}}e^{-t}}{\Gamma\left(\frac{1}{\gamma}+1\right)}\ln tdt$, hence $\int_0^{\infty}g_\gamma(t)dt=0$; therefore, the last expectation in \eqref{eq:tower_property} can simply be rewritten as:
\begin{equation}\label{eq:expectation_sign}
\mathbb{E}\left[\int_0^{\infty}h_{E_{-i}}(t)g_\gamma(t)dt\right]=\mathbb{E}\left[\int_0^{\infty}(h_{E_{-i}}(t)-h_{E_{-i}}(t_\gamma))g_\gamma(t)dt\right].
\end{equation}
\begin{enumerate}
    \item[i)] If $0<\alpha<1$, then $h_{E_{-i}}(t)-h_{E_{-i}}(t_\gamma)<0$ for $t<t_\gamma$, but $>0$ for $t>t_\gamma$, henceforth $(h_{E_{-i}}(t)-h_{E_{-i}}(t_\gamma))g_\gamma(t)\geq0$ for all $t>0$, with equality if and only if $t=t_\gamma$. Therefore, \eqref{eq:expectation_sign} and subsequently \eqref{eq:tower_property} are positive for all $i=1,\dots,d$, thus \eqref{eq:Q-analytic-derivative} is negative, i.e., $\frac{\partial}{\partial\gamma}Q_j(\gamma,\alpha;\Theta)<0$, and so $Q_j(\gamma,\alpha;\Theta)$ is strictly decreasing in $\gamma$. 
    \item[ii)] If $\alpha>1$, then $h_{E_{-i}}(t)-h_{E_{-i}}(t_\gamma)>0$ for $t<t_\gamma$, but $<0$ for $t>t_\gamma$, henceforth $(h_{E_{-i}}(t)-h_{E_{-i}}(t_\gamma))g_\gamma(t)\leq0$ for all $t>0$, with equality if and only if $t=t_\gamma$. Therefore, \eqref{eq:expectation_sign} and subsequently \eqref{eq:tower_property} are negative for all $i=1,\dots,d$, thus \eqref{eq:Q-analytic-derivative} is positive, i.e., $\frac{\partial}{\partial\gamma}Q_j(\gamma,\alpha;\Theta)>0$, and so $Q_j(\gamma,\alpha;\Theta)$ is strictly increasing in $\gamma$.
\end{enumerate}
Combining the three cases of $\alpha$, the proof is complete. 
\end{proof}

About the second proof, based on Lemma \ref{lemma:Q-exponential}, we shall apply a convex-order comparison and a coordinate-wise replacement argument on that probabilistic representation. To this end, we first recall the definition and some properties of convex ordering; also see Section 3.A of
\cite{shaked2007stochastic}. For random variables $X$ and $Y$, we
say $X$ is less than $Y$ in the convex order, denoted by $X\leq_{\mathrm{cx}}Y$, if $\mathbb E[\phi(X)]\leq\mathbb E[\phi(Y)]$ for every convex function $\phi$ whenever the expectations exist. 

\begin{proposition}[Quantile characterization of convex order; Theorem 3.A.5 of \cite{shaked2007stochastic}]
\label{prop:convex-order-quantile}
Given that $X$ and $Y$ have equal finite means and quantile functions $q_X$ and
$q_Y$, respectively, then
\[
X\leq_{\mathrm{cx}}Y
\quad\Longleftrightarrow\quad
\int_0^u q_X(v) dv
\geq
\int_0^u q_Y(v) dv,
\qquad 0\leq u\leq1.
\]
\end{proposition}

\begin{proposition}[Theorem 3.A.43 of \cite{shaked2007stochastic}]
\label{prop:convex-order-strict1}
Given that $X\leq_{\mathrm{cx}}Y$, and $\mathbb E[\phi(X)]=\mathbb E[\phi(Y)]$ for some strictly convex function $\phi$ whenever the expectations are finite,
then $X$ and $Y$ have the same distribution. 
\end{proposition}

Proposition \ref{prop:convex-order-strict1} has an immediate variant as follows, which forms the backbone of our second proof of Theorem \ref{thm:qin}.

\begin{corollary}\label{cor:convex-order-strict2}
Given that $X\leq_{\mathrm{cx}}Y$, and $X$ and $Y$
have different distributions, then $\mathbb E[\phi(X)]<\mathbb E[\phi(Y)]$ (resp. $\mathbb E[\phi(X)]>\mathbb E[\phi(Y)]$)
for every strictly convex (resp. strictly concave) $\phi$ whenever the expectations are finite.    
\end{corollary}

\begin{proof}
The convex case follows from Proposition \ref{prop:convex-order-strict1},
since equality would imply that $X$ and $Y$ have the same distribution.
The concave case follows by applying the convex case to $-\phi$.
\end{proof}

We are now ready to state the second proof based on convex ordering of Theorem \ref{thm:qin}.

\begin{proof}[Second Proof of Theorem \ref{thm:qin}]
Fix $j\in\{1,\ldots,d\}$. The case of $\alpha=1$ is the same as that in the first proof, hence we omit it here, and only focus on the case of $\alpha\neq1$. Rewrite \eqref{eq:Q-V-rho} as:
\begin{equation}\label{eq:Q-V-rho-rewritten}
    Q_j(\gamma,\alpha;\Theta)
=
\mathbb E\left[
\left(
\sum_{i=1}^d \left(w_j^i(\Theta)\right)^{\frac{1}{\alpha}} \left(\frac{E_i^{\frac{1}{\gamma}}}{\Gamma\left(\frac{1}{\gamma}+1\right)}\right)^{\frac{1}{\alpha}}
\right)^\alpha
\right]=:\mathbb{E}\left(\phi(G_{\gamma,1},\dots,G_{\gamma,d})\right),
\end{equation}
where $G_{\gamma,i}:=\frac{E_i^{\frac{1}{\gamma}}}{\Gamma\left(\frac{1}{\gamma}+1\right)}$ for $i=1,\dots,d$, and $\phi(v_1,\dots,v_d):=\left(\sum_{l=1}^d\left(w_j^l(\Theta)\right)^{\frac{1}{\alpha}}v_l^{\frac{1}{\alpha}}\right)^\alpha$. Clearly, $G_{\gamma,i}$'s are independent copies of $G_{\gamma}:=\frac{E^{\frac{1}{\gamma}}}{\Gamma\left(\frac{1}{\gamma}+1\right)}$, with $E\sim\text{Exp}(1)$.

Furthermore, by the definition of quantile function, we can solve, for $0<u<1$:
$$F_{G_\gamma}(q_\gamma(u))=\mathbb{P}\left(\frac{E^{\frac{1}{\gamma}}}{\Gamma\left(\frac{1}{\gamma}+1\right)}\leq q_\gamma(u)\right)=1-e^{-\left(q_\gamma(u)\Gamma\left(\frac{1}{\gamma}+1\right)\right)^\gamma}=u\implies q_\gamma(u)=\frac{(-\ln(1-u))^{\frac{1}{\gamma}}}{\Gamma\left(\frac{1}{\gamma}+1\right)}.$$

For $\gamma_1>\gamma_2>0$,
\[
\frac{q_{\gamma_1}(u)}{q_{\gamma_2}(u)}
=
\frac{\Gamma\left(\frac{1}{\gamma_2}+1\right)}{\Gamma\left(\frac{1}{\gamma_1}+1\right)}
(-\ln(1-u))^{\frac{1}{\gamma_1}-\frac{1}{\gamma_2}},
\]
which is strictly decreasing in $u\in(0,1)$, tending to $\infty$
as $u\downarrow0$ and to $0$ as $u\uparrow1$, hence there exists a unique $u_0\in(0,1)$ such that $\frac{q_{\gamma_1}(u_0)}{q_{\gamma_2}(u_0)}=1$; and
$q_{\gamma_1}(u)>q_{\gamma_2}(u)$ for $0<u<u_0$, but
$q_{\gamma_1}(u)<q_{\gamma_2}(u)$ for $u_0<u<1$. Also note that $\int_0^1q_{\gamma}(v)dv=1$ for all $\gamma>0$, since $\int_0^1q_{\gamma}(v)dv=\mathbb{E}(q_{\gamma}(U))=\mathbb{E}(G_\gamma)$ where $U\sim U(0,1)$, where the second equality holds by the property of probability integral transform, and also $\mathbb{E}(G_\gamma)=1$ as $\mathbb{E}\left(E^{\frac{1}{\gamma}}\right)=\Gamma\left(\frac{1}{\gamma}+1\right)$. Obviously, $\int_0^u(q_{\gamma_1}(v)-q_{\gamma_2}(v))dv=0$ for $u=0,1$, and
$\int_0^u(q_{\gamma_1}(v)-q_{\gamma_2}(v))dv>0$ for $0<u\leq u_0$; whereas for $u_0<u<1$,
we have:
\[
\int_0^u
(q_{\gamma_1}(v)-q_{\gamma_2}(v))dv
=\int_0^1
(q_{\gamma_1}(v)-q_{\gamma_2}(v))dv-\int_u^1
(q_{\gamma_1}(v)-q_{\gamma_2}(v))dv=
\int_u^1
(q_{\gamma_2}(v)-q_{\gamma_1}(v))dv
>0;
\]
altogether yields $\int_0^u(q_{\gamma_1}(v)-q_{\gamma_2}(v))dv\geq0$ for any $u\in[0,1]$. By Proposition \ref{prop:convex-order-quantile}, 
\begin{equation}\label{eq:G_convex_ordering}
    G_{\gamma_1}\leq_{\mathrm{cx}}G_{\gamma_2}\text{ for any }\gamma_1>\gamma_2>0.
\end{equation}

Besides, for any $i=1,\dots,d$,
\begin{equation}\label{eq:phi_complexity}
\frac{\partial^2\phi}{\partial v_i^2}=\left(w_j^i(\Theta)\right)^{\frac{1}{\alpha}}\left(\frac1\alpha-1\right)
\left(\sum_{l\neq i}\left(w_j^l(\Theta)\right)^{\frac{1}{\alpha}}v_l^{\frac{1}{\alpha}}\right) v_i^{\frac{1}{\alpha}-2}
\left(\sum_{l=1}^d\left(w_j^l(\Theta)\right)^{\frac{1}{\alpha}}v_l^{\frac{1}{\alpha}}\right)^{\alpha-2}\left\{\begin{aligned}
    &>0,\text{ for }0<\alpha<1;\\
    &<0,\text{ for }\alpha>1.
\end{aligned}\right.
\end{equation}
Therefore, $\phi$ is strictly convex (resp. strictly concave) in $v_i$ for any $i=1,\dots,d$ when $0<\alpha<1$ (resp. $\alpha>1$).

For any $\gamma_1>\gamma_2>0$, and for any $i=1,\dots,d$, denote $\mathbf{G}_{\gamma_2}^{(i)}:=(G_{\gamma_2,1},\cdots,G_{\gamma_2,i},G_{\gamma_1,i+1},\cdots,G_{\gamma_1,d})$; as a custom, we denote $\mathbf{G}_{\gamma_2}^{(0)}=\mathbf{G}_{\gamma_1}:=(G_{\gamma_1,1},\cdots,G_{\gamma_1,d})$, and $\mathbf{G}_{\gamma_2}:=\mathbf{G}_{\gamma_2}^{(d)}$. In particular, note that the only difference between $\mathbf{G}_{\gamma_2}^{(i-1)}$ and $\mathbf{G}_{\gamma_2}^{(i)}$  lies in the $i$-th entry, namely $G_{\gamma_1,i}$ in the former and $G_{\gamma_2,i}$ in the latter; in light of \eqref{eq:G_convex_ordering}, $G_{\gamma_1,i}\leq_{\rm cx}G_{\gamma_2,i}$, and they clearly have different distributions. By further denoting $\phi_i(v_i):=\phi(v_1,\dots,v_d)|_{v_l,\forall l\neq i}$, we have the following:
\begin{enumerate}
    \item[1)] If $0<\alpha<1$, then $\phi_i(v_i)$ is strictly convex in $v_i$ by \eqref{eq:phi_complexity}. By Corollary \ref{cor:convex-order-strict2},
    $$\mathbb{E}[\phi(\mathbf{G}_{\gamma_2}^{(i-1)})]=\mathbb{E}[\phi_i(G_{\gamma_1,i})]<\mathbb{E}[\phi_i(G_{\gamma_2,i})]=\mathbb{E}[\phi(\mathbf{G}_{\gamma_2}^{(i)})];$$
    applying the same argument for $i=1,\dots,d$ then yields:
    $$\mathbb{E}[\phi(\mathbf{G}_{\gamma_1})]=\mathbb{E}[\phi(\mathbf{G}_{\gamma_2}^{(0)})]<\mathbb{E}[\phi(\mathbf{G}_{\gamma_2}^{(1)})]<\cdots<\mathbb{E}[\phi(\mathbf{G}_{\gamma_2}^{(d)})]=\mathbb{E}[\phi(\mathbf{G}_{\gamma_2})].$$
    In light of \eqref{eq:Q-V-rho-rewritten}, $Q_j(\gamma_1,\alpha;\Theta)=\mathbb{E}[\phi(\mathbf{G}_{\gamma_1})]<\mathbb{E}[\phi(\mathbf{G}_{\gamma_2})]=Q_j(\gamma_2,\alpha;\Theta)$ for any $\gamma_1>\gamma_2>0$, thus $Q_j(\gamma,\alpha;\Theta)$ is strictly decreasing in $\gamma$.
    \item[2)] If $\alpha>1$, then $\phi_i(v_i)$ is strictly concave in $v_i$ by \eqref{eq:phi_complexity}. By Corollary \ref{cor:convex-order-strict2},
    $$\mathbb{E}[\phi(\mathbf{G}_{\gamma_2}^{(i-1)})]=\mathbb{E}[\phi_i(G_{\gamma_1,i})]>\mathbb{E}[\phi_i(G_{\gamma_2,i})]=\mathbb{E}[\phi(\mathbf{G}_{\gamma_2}^{(i)})];$$
    applying the same argument for $i=1,\dots,d$ then yields: $$\mathbb{E}[\phi(\mathbf{G}_{\gamma_1})]=\mathbb{E}[\phi(\mathbf{G}_{\gamma_2}^{(0)})]>\mathbb{E}[\phi(\mathbf{G}_{\gamma_2}^{(1)})]>\cdots>\mathbb{E}[\phi(\mathbf{G}_{\gamma_2}^{(d)})]=\mathbb{E}[\phi(\mathbf{G}_{\gamma_2})].$$
    In light of \eqref{eq:Q-V-rho-rewritten}, $Q_j(\gamma_1,\alpha;\Theta)=\mathbb{E}[\phi(\mathbf{G}_{\gamma_1})]>\mathbb{E}[\phi(\mathbf{G}_{\gamma_2})]=Q_j(\gamma_2,\alpha;\Theta)$ for any $\gamma_1>\gamma_2>0$, thus $Q_j(\gamma,\alpha;\Theta)$ is strictly increasing in $\gamma$.
\end{enumerate}
The desired claim then follows.
\end{proof}

The monotonicity result of Theorem \ref{thm:qin} describes how the limiting tail constant varies
with the regular-variation index $\gamma$, which governs the strength of
upper-tail dependence. The following comparison between $Q_j(\gamma,\alpha;\Theta)$ and its comonotonic benchmark is crucial in establishing our main results in the sequel:

\begin{proposition}
\label{prop:Q-comonotonic}
Under the assumptions and conditions specified in Theorem \ref{thm:Qexist}, define
\[
Q_{j,\mathrm{com}}(\alpha;\Theta)
:=
\left(
\sum_{i=1}^d
\bigl(w_j^i(\Theta)\bigr)^{\frac{1}{\alpha}}
\right)^\alpha.
\]
Then, for every finite $\gamma>0$ and every
$j\in\{1,\ldots,d\}$,
\[Q_j(\gamma,\alpha;\Theta)\left\{
\begin{aligned}
&>
Q_{j,\mathrm{com}}(\alpha;\Theta),
&& 0<\alpha<1;\\
&=
Q_{j,\mathrm{com}}(1;\Theta),
&& \alpha=1;\\
&<
Q_{j,\mathrm{com}}(\alpha;\Theta),
&& \alpha>1.
\end{aligned}\right.
\]
\end{proposition}
\begin{remark}
    As we shall establish in Corollary \ref{cor:Q-comonotonic-limit}, $Q_{j,\mathrm{com}}(\alpha;\Theta)$ is the comonotonic limit of $Q_j(\gamma,\alpha;\Theta)$ as $\gamma\uparrow\infty$, justifying the use of subscript ``com''.
\end{remark}
\begin{proof}
Fix $j\in\{1,\ldots,d\}$. Based on \eqref{eq:Q-V-rho-rewritten}, we have:
\[
\frac{Q_j(\gamma,\alpha;\Theta)}{Q_{j,\mathrm{com}}(\alpha;\Theta)}
=
\mathbb E\left[
\left(
\sum_{i=1}^d \kappa_j^iG_{\gamma,i}^{\frac{1}{\alpha}}
\right)^\alpha
\right],
\]
where $\kappa_j^i:=\frac{\left(w_j^i(\Theta)\right)^{\frac{1}{\alpha}}}{\sum_{l=1}^d\left(w_j^l(\Theta)\right)^{\frac{1}{\alpha}}}$ for $i=1,\dots,d$; clearly $\kappa_j^i>0$ for all $i=1,\dots,d$, and $\sum_{i=1}^d\kappa_j^i=1$.
\begin{enumerate}
    \item[i)] If $0<\alpha<1$, by strict concavity of $x\mapsto x^\alpha$ and Jensen's inequality with weights $\kappa_j^i$, $i=1,\dots,d$, we have:
\[
\left(
\sum_{i=1}^d \kappa_j^iG_{\gamma,i}^{\frac{1}{\alpha}}
\right)^\alpha
>
\sum_{i=1}^d \kappa_j^iG_{\gamma,i},
\]
the inequality is strict since the $G_{\gamma,i}$'s are not all equal almost surely, as the event set $\{G_{\gamma,1}=\cdots=G_{\gamma,d}\}$ has a zero measure.
By recalling that $\mathbb{E}(G_{\gamma,i})=1$ for any $i=1,\dots,d$ in light of the second proof of Theorem \ref{thm:qin}, we have:
$$\frac{Q_j(\gamma,\alpha;\Theta)}{Q_{j,\mathrm{com}}(\alpha;\Theta)}=\mathbb{E}\left(
\sum_{i=1}^d \kappa_j^iG_{\gamma,i}^{\frac{1}{\alpha}}
\right)^\alpha
>\mathbb{E}\left(\sum_{i=1}^d \kappa_j^iG_{\gamma,i}\right)=\sum_{i=1}^d \kappa_j^i\mathbb{E}\left(G_{\gamma,i}\right)=\sum_{i=1}^d \kappa_j^i=1,$$
hence $Q_j(\gamma,\alpha;\Theta)>Q_{j,\mathrm{com}}(\alpha;\Theta)$.
\item[ii)] If $\alpha=1$, by \eqref{eq:alpha_equals_1}, we immediately conclude that $Q_j(\gamma,1;\Theta)=\sum_{i=1}^d
w_j^i(\Theta)=Q_{j,\mathrm{com}}(1;\Theta)$.
\item[iii)] If $\alpha>1$, then using the likewise arguments as in Case i), with strict convexity of $x\mapsto x^\alpha$, we have:
\[
\left(
\sum_{i=1}^d \kappa_j^iG_{\gamma,i}^{\frac{1}{\alpha}}
\right)^\alpha
<
\sum_{i=1}^d \kappa_j^iG_{\gamma,i};
\]
taking expectations on both sides gives $\frac{Q_j(\gamma,\alpha;\Theta)}{Q_{j,\mathrm{com}}(\alpha;\Theta)}=\mathbb{E}\left(
\sum_{i=1}^d \kappa_j^iG_{\gamma,i}^{\frac{1}{\alpha}}
\right)^\alpha
<\mathbb{E}\left(\sum_{i=1}^d \kappa_j^iG_{\gamma,i}\right)=1$, hence
$Q_j(\gamma,\alpha;\Theta)<Q_{j,\mathrm{com}}(\alpha;\Theta)$.
\end{enumerate}
\end{proof}

\begin{corollary}
\label{cor:Q-comonotonic-limit}
Under the assumptions and conditions specified in Theorem \ref{thm:Qexist}, for every $j\in\{1,\cdots,d\}$ and any $\alpha>0$,
$Q_j(\gamma,\alpha;\Theta)\rightarrow Q_{j,\mathrm{com}}(\alpha;\Theta)$ as
$\gamma\uparrow\infty$.
\end{corollary}

\begin{proof}
Fix $j\in\{1,\cdots,d\}$. As $\gamma\uparrow\infty$, $\frac{1}{\gamma}\to0$, hence $E_i^{\frac{1}{\gamma}}\overset{a.s.}{\to}1$ for any $i=1,\dots,d$, and $\Gamma\left(\frac{1}{\gamma}+1\right)\to\Gamma(1)=1$. Therefore, recalling $G_{\gamma,i}=\frac{E_i^{\frac{1}{\gamma}}}{\Gamma\left(\frac{1}{\gamma}+1\right)}$ as defined immediately after \eqref{eq:Q-V-rho-rewritten},  $G_{\gamma,i}\overset{a.s.}{\to}1$ for all $i=1,\dots,d$, and
$$\left(\sum_{i=1}^d\left(w_j^i(\Theta)\right)^{\frac{1}{\alpha}}G_{\gamma,i}^{\frac{1}{\alpha}}\right)^\alpha\overset{a.s.}{\to}\left(\sum_{i=1}^d\left(w_j^i(\Theta)\right)^{\frac{1}{\alpha}}\right)^\alpha=Q_{j,\mathrm{com}}(\alpha;\Theta),\text{ as }\gamma\uparrow\infty.$$
Therefore, we have:
$$Q_j(\gamma,\alpha;\Theta)=\mathbb{E}\left[\left(\sum_{i=1}^d\left(w_j^i(\Theta)\right)^{\frac{1}{\alpha}}G_{\gamma,i}^{\frac{1}{\alpha}}\right)^\alpha\right]\to\mathbb{E}\left[\left(\sum_{i=1}^d\left(w_j^i(\Theta)\right)^{\frac{1}{\alpha}}\right)^\alpha\right]
=Q_{j,\mathrm{com}}(\alpha;\Theta),\text{ as }\gamma\uparrow\infty,$$
as long as limit-taking and expectation-taking can be interchanged; indeed, to show that it is viable, we shall establish a uniform bound for $\left(\sum_{i=1}^d\left(w_j^i(\Theta)\right)^{\frac{1}{\alpha}}G_{\gamma,i}^{\frac{1}{\alpha}}\right)^\alpha$ that is
independent of $\gamma$, by then the dominated convergence theorem can warrant the aforementioned interchange. To this end, we first note the following inequality: if $c_i>0$ for all $i=1,\dots,d$, then
\begin{equation}\label{eq:power-sum-inequality}
    \left(\sum_{i=1}^dc_i^{\frac{1}{\alpha}}\right)^\alpha\leq d^{(\alpha-1)_+}\sum_{i=1}^d c_i,\text{ where }(\alpha-1)_+:=\max(\alpha-1,0).\footnote{To see this, we observe that when $0<\alpha\leq1$, $\frac{1}{\alpha}\geq1$, and
$\left(\sum_{i=1}^dc_i^{\frac{1}{\alpha}}\right)^\alpha=\left(\sum_{l=1}^d c_l\right)\left(\sum_{i=1}^d\left(\frac{c_i}{\sum_{l=1}^d c_l}\right)^{\frac{1}{\alpha}}\right)^\alpha\leq\left(\sum_{l=1}^d c_l\right)\left(\sum_{i=1}^d\frac{c_i}{\sum_{l=1}^d c_l}\right)^\alpha=\sum_{l=1}^d c_l$;
whereas for $\alpha>1$, the function $x\mapsto x^{\frac{1}{\alpha}}$ is concave as $0<\frac{1}{\alpha}<1$, then by Jensen's inequality,
$\frac{1}{d}\sum_{i=1}^d c_i^{\frac{1}{\alpha}}\leq\left(\frac{1}{d}\sum_{i=1}^dc_i\right)^{\frac{1}{\alpha}}$, implying $\left(\sum_{i=1}^d c_i^{\frac{1}{\alpha}}\right)^\alpha\leq d^{\alpha-1}\sum_{i=1}^dc_i$.}
\end{equation}
Furthermore, taking $\gamma>1$, note that $\frac{1}{\gamma}+1\in(1,2)$, and 
$G_{\gamma,i}\leq \frac{E_i^{\frac{1}{\gamma}}}{m_\Gamma} \leq \frac{1+E_i}{m_\Gamma}$, where $m_\Gamma:=\inf_{x\in(1,2)}\Gamma(x)\in\mathbb{R}_+$, then together with \eqref{eq:power-sum-inequality}, we have:
\[
0\leq \left(\sum_{i=1}^d\left(w_j^i(\Theta)\right)^{\frac{1}{\alpha}}G_{\gamma,i}^{\frac{1}{\alpha}}\right)^\alpha\leq
d^{(\alpha-1)_+}
\sum_{i=1}^d w_j^i(\Theta)G_{\gamma,i}\leq
\frac{1}{m_\Gamma} d^{(\alpha-1)_+}
\sum_{i=1}^d w_j^i(\Theta)(1+E_i),
\]
where the rightmost expression is integrable and independent of $\gamma$.
\end{proof}

\begin{theorem}[Asymptotic V@R super/subadditivity for risks with a common Fr\'echet index]
\label{thm:main}
Under Assumptions \ref{AssumeArchi} and \ref{Assumefrechet}, recalling $D_p$ in \eqref{eq:diversification-ratio}, V@R is:
\begin{enumerate}[label=(\roman*)]
    \item asymptotically superadditive with $\liminf_{p\uparrow 1}D_p\geq 1$, for $0<\alpha<1$; or
    \item asymptotically additive with
    $\lim_{p\uparrow 1}D_p= 1$, for $\alpha=1$; or 
    \item asymptotically subadditive with
    $\limsup_{p\uparrow 1}D_p\leq 1$, for $\alpha>1$.
\end{enumerate}
\end{theorem}
\begin{remark}
    We only require a common Fr\'echet index among $X_i$'s, but not  tail equivalence among them. 
\end{remark}
\begin{proof}[Proof of Theorem \ref{thm:main}]
Let $\alpha$ be the common Fr\'echet index. We use a subsequential argument. Let $\{p_n\}_{n\geq1}$ be an arbitrary sequence with
$p_n\uparrow1$, and set $x_{i,n}:=\V@R_{p_n}(X_i)$, $A_n:=\sum_{i=1}^d x_{i,n}$, $r_{i,n}:=\frac{x_{i,n}}{A_n}\in[0,1]$.
Since each $X_i\in\operatorname{MDA}(\Phi_\alpha)$ has infinite right
endpoint, $x_{i,n}\rightarrow\infty$ for all $i=1,\dots,d$, and hence $A_n\rightarrow\infty$.
After passing to a subsequence, we may assume that $r_{i,n}\rightarrow r_i$ for $i=1,\ldots,d$,
for some $r_i\in[0,1]$ with $\sum_{i=1}^dr_i=1$. 
By Lemma \ref{lemma:tail-at-var}, we have:
\begin{equation}\label{eq:tail-at-var}
\frac{\overline F_i(x_{i,n})}{1-p_n}
\rightarrow1,\text{ as }p_n\uparrow1.
\end{equation}
Define
$\mathcal I:=\{k:r_k>0\}\neq\varnothing$; indeed, since $\sum_{k=1}^dr_k=1$, $r_k$'s cannot all equal 0. Let $i\in\mathcal I$ and fix $t>0$. Since
$r_{i,n}\to r_i>0$, the uniform convergence theorem for regularly varying
functions (Theorem 1.5.2 of \cite{bgt87}) together with \eqref{eq:tail-at-var} gives:
\begin{align}
\frac{\overline F_i(A_nt)}{1-p_n}
=
\frac{\overline F_i(A_nt)}
     {\overline F_i(x_{i,n})}
\frac{\overline F_i(x_{i,n})}{1-p_n}=
\frac{\overline F_i(A_nt)}
     {\overline F_i(A_nr_{i,n})}
\frac{\overline F_i(x_{i,n})}{1-p_n}
\rightarrow
\left(\frac{t}{r_i}\right)^{-\alpha}
=
r_i^\alpha t^{-\alpha},\text{ as }p_n\uparrow1.
\label{eq:tail-var-normalization}
\end{align}
Particularly, for $i,j\in\mathcal I$, and with $t=1$,
\begin{equation}\label{eq:w-r}
\frac{\overline F_i(A_n)}
     {\overline F_j(A_n)}
\rightarrow
\left(\frac{r_i}{r_j}\right)^\alpha,\text{ as }p_n\uparrow1.
\end{equation}
Now, take $\Theta:=\{A_n\}_{n\geq1}$, so that $w_j^i(\Theta)=\left(\frac{r_i}{r_j}\right)^\alpha$ for $i,j\in\mathcal I$ by Definition \ref{def:tailratio}; indeed, when $|\mathcal I|=1$, let say $\mathcal I=\{i\}$, then $r_i=1$ and $w_i^i(\Theta)=1$; while for $|\mathcal I|\geq2$,  $r_i,r_j\in(0,1)$ for any $i,j\in\mathcal I$ by the very construction of $\mathcal I$, hence $w_j^i(\Theta)\in\mathbb{R}_+$. 

Consider the non-trivial case of $|\mathcal I|\geq2$, and write
$S_{\mathcal I}:=\sum_{k\in\mathcal I}X_k$.
Note that by the very definition of Archimedean copula, every marginal subvector of an Archimedean copula is again Archimedean
with the same generator; also see Section 4.6 of \cite{n06}. For any fixed $t>0$, let $\Theta_t:=\{tA_n\}$. For
$i,j\in\mathcal I$,
\[
\frac{\overline F_i(tA_n)}
     {\overline F_j(tA_n)}
=
\frac{\overline F_i(tA_n)}
     {\overline F_i(A_n)}
\frac{\overline F_i(A_n)}
     {\overline F_j(A_n)}
\frac{\overline F_j(A_n)}
     {\overline F_j(tA_n)}
\rightarrow
t^{-\alpha}w_j^i(\Theta)t^\alpha
=
w_j^i(\Theta),\text{ as }p_n\uparrow1,
\]
thus $w_j^i(\Theta_t)=w_j^i(\Theta)\in\mathbb{R}_+$ for all $i,j\in\mathcal I$.
Since the integral density representation \eqref{eq:small_h}-\eqref{qjint} depends on
the sequence $\Theta$ only through these limiting tail ratios, $\lim_{n\to\infty}\frac {\mathbb{P}\left(S_{\mathcal I}>tA_n\right)}{\overline{F}_j(tA_n)}=:Q_j^{\mathcal I}(\gamma,\alpha;\Theta_t)=Q_j^{\mathcal I}(\gamma,\alpha;\Theta):=\lim_{k\to\infty}\frac {\mathbb{P}\left(S_{\mathcal I}>\theta_k\right)}{\overline{F}_j(\theta_k)}$; together with \eqref{eq:tail-var-normalization}, we have:
\begin{equation}\label{eq:SI-tail}
\frac{\mathbb P(S_{\mathcal I}>tA_n)}
     {1-p_n}=\frac{\mathbb P(S_{\mathcal I}>tA_n)}
     {\overline{F}_j(A_nt)}\frac{\overline{F}_j(A_nt)}
     {1-p_n}
\rightarrow
Q_j^{\mathcal I}(\gamma,\alpha;\Theta)r_j^\alpha t^{-\alpha}=:C_\gamma t^{-\alpha},\text{ as }n\to\infty;
\end{equation}
in particular, the limit only depends on the set $\mathcal I$, but not on the specific choice of $j\in\mathcal I$; this is also the case for $C_\gamma=Q_j^{\mathcal I}(\gamma,\alpha;\Theta)r_j^\alpha$ that is independent of the choice of $j\in\mathcal I$. 
Recalling the definition of $w_j^i(\Theta)$ and $\mathcal I$, we have $Q^{\mathcal I}_{j,\rm{com}}(\alpha;\Theta):=\left(
\sum_{i\in\mathcal I}
(w_j^i(\Theta))^{\frac{1}{\alpha}}
\right)^\alpha
=
\left(
\sum_{i\in\mathcal I}\frac{r_i}{r_j}
\right)^\alpha
=
r_j^{-\alpha}$.
Therefore, in light of Proposition \ref{prop:Q-comonotonic}, when $|\mathcal I|\geq2$,
\begin{equation}\label{eq:C-strict}
C_\gamma=\frac{Q_j^{\mathcal I}(\gamma,\alpha;\Theta)}{Q^{\mathcal I}_{j,\rm{com}}(\alpha;\Theta)}\frac{Q^{\mathcal I}_{j,\rm{com}}(\alpha;\Theta)}{r_j^{-\alpha}}=\frac{Q_j^{\mathcal I}(\gamma,\alpha;\Theta)}{Q^{\mathcal I}_{j,\rm{com}}(\alpha;\Theta)}\left\{\begin{aligned}
 &> 1, \text{ for }0<\alpha<1;\\
&= 1, \text{ for }\alpha=1;\\
&< 1, \text{ for }\alpha>1.
\end{aligned}\right.
\end{equation}
For the trivial case $|\mathcal I|=1$, let say $\mathcal I=\{i\}$ without loss of generality, then $r_i=1$ and
\eqref{eq:tail-var-normalization} directly gives $\frac{\mathbb P(X_i>tA_n)}{1-p_n}\rightarrow t^{-\alpha}$, and since $S_{\mathcal I}=X_i$, \eqref{eq:SI-tail} remains valid with $C_\gamma=1$.

We next show that the risks outside $\mathcal I$ are asymptotically
negligible. If $i\notin\mathcal I$, then
$\frac{x_{i,n}}{A_n}=r_{i,n}\to r_i=0$. Fix $c>0$, and choose $\delta\in(0,\alpha)$ and $K>1$. By Potter's bound
(Theorem 1.5.6 of \cite{bgt87}), there exists $x_0>0$ such that, whenever
$x,y\geq x_0$,
\[
\frac{\overline F_i(y)}{\overline F_i(x)}
\leq
K\max\left\{
\left(\frac{y}{x}\right)^{-\alpha+\delta},
\left(\frac{y}{x}\right)^{-\alpha-\delta}
\right\}.
\]
Since $x_{i,n}\to\infty$ and $\frac{cA_n}{x_{i,n}} \to \infty$, we have $x_{i,n},cA_n\geq x_0$ and
$\frac{cA_n}{x_{i,n}}\geq1$ for all sufficiently large $n$, by then $\left(\frac{cA_n}{x_{i,n}}\right)^{-\alpha+\delta}\geq\left(\frac{cA_n}{x_{i,n}}\right)^{-\alpha-\delta}$, hence:
\begin{equation}\label{eq:risks_outside_I}
\frac{\overline F_i(cA_n)}{1-p_n}
=
\frac{\overline F_i(cA_n)}
     {\overline F_i(x_{i,n})}
\frac{\overline F_i(x_{i,n})}{1-p_n}
\leq
K \left(\frac{cA_n}{x_{i,n}}\right)^{-\alpha+\delta}
\frac{\overline F_i(x_{i,n})}{1-p_n}
\rightarrow0,
\end{equation}
where the last step follows from
\eqref{eq:tail-at-var} and $-\alpha+\delta<0$.
Assume first that $\mathcal I^c\neq\varnothing$, and let
$m:=|\mathcal I^c|$. Then for every $\epsilon>0$, based on \eqref{eq:risks_outside_I} by taking $c=\frac{\epsilon}{m}$, Boole's inequality gives:
\[
\frac{
\mathbb P\left(
\sum_{i\notin\mathcal I}X_i>\epsilon A_n
\right)
}{
1-p_n
}
\leq
\sum_{i\notin\mathcal I}
\frac{
\mathbb P\left(
X_i>\frac{\epsilon}{m}A_n
\right)
}{
1-p_n
}
\rightarrow0;
\]
on the other hand, if $\mathcal I^c=\varnothing$, clearly $\mathbb P\left(
\sum_{i\notin\mathcal I}X_i>\epsilon A_n
\right)\equiv0$ for every $\epsilon>0$. Combining both cases, we have:
\begin{equation}\label{eq:negligible-sum}
\frac{
\mathbb P\left(
\sum_{i\notin\mathcal I}X_i>\epsilon A_n
\right)
}{
1-p_n
}
\rightarrow0.
\end{equation}

Fix $t>0$ and $\epsilon\in(0,t)$. Since $\mathbb P(S_{\mathcal I}>tA_n)\leq\mathbb P(S_d>tA_n)$ and
\[
\mathbb P(S_d>tA_n)
\leq
\mathbb P(S_{\mathcal I}>(t-\epsilon)A_n)
+
\mathbb P\left(
\sum_{i\notin\mathcal I}X_i>\epsilon A_n
\right),
\]
\eqref{eq:SI-tail} and \eqref{eq:negligible-sum} together imply that:
\[
C_\gamma t^{-\alpha}
\leq
\liminf_{n\to\infty}
\frac{\mathbb P(S_d>tA_n)}{1-p_n}
\leq
\limsup_{n\to\infty}
\frac{\mathbb P(S_d>tA_n)}{1-p_n}
\leq
C_\gamma(t-\epsilon)^{-\alpha},
\]
and letting $\epsilon\downarrow0$ gives:
\begin{equation}\label{eq:S-tail}
\frac{\mathbb P(S_d>tA_n)}
     {1-p_n}
\rightarrow
C_\gamma t^{-\alpha},
\qquad t>0.
\end{equation}
If $t<C_\gamma^{\frac{1}{\alpha}}$, then $C_\gamma t^{-\alpha}>1$, so
\eqref{eq:S-tail} eventually implies
$\mathbb P(S_d>tA_n)>1-p_n$, and hence
$\V@R_{p_n}(S_d)>tA_n$. Therefore,
$\liminf_{n\to\infty}
\frac{\V@R_{p_n}(S_d)}{A_n}
\geq t$
for every $t<C_\gamma^{\frac{1}{\alpha}}$, which gives:
\begin{equation}\label{eq:quantile-liminf}
\liminf_{n\to\infty}
\frac{\V@R_{p_n}(S_d)}{A_n}
\geq C_\gamma^{\frac{1}{\alpha}}.
\end{equation}

Similarly, if $t>C_\gamma^{\frac{1}{\alpha}}$, then
$C_\gamma t^{-\alpha}<1$, so \eqref{eq:S-tail} eventually implies
$\mathbb P(S_d>tA_n)<1-p_n$, and hence
$\V@R_{p_n}(S_d)\leq tA_n$. Thus,
\begin{equation}\label{eq:quantile-limsup}
\limsup_{n\to\infty}
\frac{\V@R_{p_n}(S_d)}{A_n}
\leq C_\gamma^{\frac{1}{\alpha}}.
\end{equation}
In view of \eqref{eq:quantile-liminf} and
\eqref{eq:quantile-limsup},
\begin{equation}\label{eq:quantile-inversion}
\frac{\V@R_{p_n}(S_d)}{A_n}
\rightarrow
C_\gamma^{\frac{1}{\alpha}},\text{ as }n\to\infty.
\end{equation}
Since
$A_n=\sum_{i=1}^d\V@R_{p_n}(X_i)$,
\eqref{eq:quantile-inversion} gives
\begin{equation}\label{eq:D-limit}
D_{p_n}
=
\frac{\V@R_{p_n}(S_d)}
     {\sum_{i=1}^d\V@R_{p_n}(X_i)}
\rightarrow
C_\gamma^{\frac{1}{\alpha}},\text{ as }n\to\infty.
\end{equation}
By \eqref{eq:C-strict} and recalling $C_\gamma=1$ when $|\mathcal I|=1$, along the extracted
subsequence, the limit in \eqref{eq:D-limit} satisfies
\begin{equation}\label{eq:subsequence_limit}
C_\gamma^{\frac{1}{\alpha}}\left\{\begin{aligned}
 &\geq 1, \text{ for } 0<\alpha<1;\\
&= 1,\text{ for } \alpha=1;\\
&\leq 1, \text{ for }  \alpha>1.\\
\end{aligned}\right.
\end{equation}
Since the original sequence $\{p_n\}$ is arbitrary, this yields the
statement of Theorem \ref{thm:main}; indeed,
\begin{enumerate}
    \item[i)] If $0<\alpha<1$ and $\lim\inf_{p\uparrow1}D_p<1$, there would exist $\epsilon>0$ and a sequence
$\{p_n\}\uparrow1$ such that
$D_{p_n}\leq1-\epsilon$
for every $n$. On the other hand, by \eqref{eq:D-limit} and \eqref{eq:subsequence_limit}, there exists a subsequence $\{p_{n_i}\}$ such that $D_{p_{n_i}}$ converges to a limit $\geq1$, which results in a contradiction. Therefore, $\lim\inf_{p\uparrow1}D_p\geq1$.
\item[ii)] If $\alpha>1$ and $\lim\sup_{p\uparrow1}D_p>1$, there would exist $\epsilon>0$ and a sequence
$\{p_n\}\uparrow1$ such that
$D_{p_n}\geq1+\epsilon$
for every $n$. On the other hand, by \eqref{eq:D-limit} and \eqref{eq:subsequence_limit}, there exists a subsequence $\{p_{n_i}\}$ such that $D_{p_{n_i}}$  converges to a limit $\leq1$, which results in a contradiction. Therefore, $\lim\sup_{p\uparrow1}D_p\leq1$.
\item[iii)] When $\alpha=1$, if $\lim\inf_{p\uparrow1}D_p<1$, there would exist $\epsilon>0$ and a sequence
$\{p_n\}\uparrow1$ such that
$D_{p_n}\leq1-\epsilon$
for every $n$, yet by \eqref{eq:D-limit} and \eqref{eq:subsequence_limit},  every
sequence $\{p_n\}\uparrow1$ has a further subsequence $\{p_{n_i}\}$ such that $D_{p_{n_i}}$
converges to 1, which leads to a contradiction, so that $\lim\inf_{p\uparrow1}D_p\geq1$. On the other hand, by the same arguments, $\lim\sup_{p\uparrow1}D_p>1$ also leads to a contradiction, so $\lim\sup_{p\uparrow1}D_p\leq1$. Altogether leads to $\lim_{p\uparrow1}D_p=1$. \qedhere
\end{enumerate}  
\end{proof}

The following theorem shows that without assuming tail equivalence among risks, our obtained result in Theorem \ref{thm:main} is sharp, in the sense that it cannot be further improven unless stronger assumptions are imposed.  
\begin{theorem}
\label{thm:limsup-liminf-sharp}
Let $\alpha>0$, $\alpha\neq1$, and let $\gamma\in(0,\infty)$.
There exist two marginal distributions
$F_1,F_2\in\operatorname{MDA}(\Phi_\alpha)$, coupled by an
Archimedean survival copula satisfying Assumption \ref{AssumeArchi},
such that the diversification ratio $D_p$ as defined in \eqref{eq:diversification-ratio} does not converge as
$p\uparrow1$. More precisely, there exist sequences
$\{p_n^{(1)}\}\uparrow1$ and $\{p_n^{(2)}\}\uparrow1$ such that as $n\to\infty$,
\[
D_{p_n^{(1)}}\rightarrow d_{\alpha,\gamma},
\qquad
D_{p_n^{(2)}}\rightarrow1,
\]
where $d_{\alpha,\gamma}>1$ if $0<\alpha<1$,
but $<1$ if $\alpha>1$.
Thus, for any $\alpha\neq1$, $D_p$ has no limit as $p\uparrow1$.
\end{theorem}

\begin{proof}
    See Appendix \ref{sect:proof_of_sharpness}.
\end{proof}

The following result establishes
the corresponding closure property of Fr\'echet MDA under aggregation under Assumptions \ref{AssumeArchi} and
\ref{Assumefrechet} without requiring global tail equivalence of the marginal
risks.

\begin{corollary}\label{cor:sum-frechet}
Under Assumptions \ref{AssumeArchi} and
\ref{Assumefrechet}, 
$S_d=\sum_{i=1}^d X_i\in\operatorname{MDA}(\Phi_\alpha)$.
\end{corollary}

\begin{proof}
    As stated in Assumption \ref{Assumefrechet}, $X_i\in\operatorname{MDA}(\Phi_\alpha)$ for $i=1,\dots,d$, and tail equivalence is not required. It suffices to show that for every $t>0$,
$\frac{\mathbb P(S_d>tx)}
     {\mathbb P(S_d>x)}
\rightarrow t^{-\alpha}$ as $x\to\infty$.
Let $\{x_n\}$ be an arbitrary sequence with $x_n\to\infty$ and set
$L_n:=\sum_{i=1}^d\overline F_i(x_n)$. By compactness of the unit simplex,
after passing to a subsequence, we may assume that
$\frac{\overline F_i(x_n)}{L_n}\rightarrow\lambda_i$ for $i=1,\ldots,d$,
for some $\lambda_i\geq0$ such that $\sum_{i=1}^d\lambda_i=1$. Define
$\mathcal I:=\{k:\lambda_k>0\}$ and fix $j\in\mathcal I$, then
$\frac{\overline F_i(x_n)}
     {\overline F_j(x_n)}
\rightarrow\frac{\lambda_i}{\lambda_j}\mathbbm{1}_{i\in\mathcal I}$.
Take $\Theta:=\{x_n\}_{n\geq1}$ and again denote
$S_{\mathcal I}:=\sum_{i\in\mathcal I}X_i$. For every fixed $c>0$ and
$i,\ell\in\mathcal I$, regular variation gives
$\frac{\overline F_i(cx_n)}
     {\overline F_\ell(cx_n)}
\rightarrow
\frac{\lambda_i}{\lambda_\ell}$, thus $\Theta_c:=\{cx_n\}_{n\geq1}$ has the same limiting tail ratios as $\Theta$.
Since the integral density representation \eqref{eq:small_h}-\eqref{qjint} depends on
the sequence $\Theta$ only through these limiting tail ratios,
\begin{equation}\label{eq:SI_frechet}
\lim_{n\to\infty}\frac {\mathbb{P}\left(S_{\mathcal I}>cx_n\right)}{\overline{F}_j(cx_n)}=:Q_j^{\mathcal I}(\gamma,\alpha;\Theta_c)=Q_j^{\mathcal I}(\gamma,\alpha;\Theta):=\lim_{n\to\infty}\frac {\mathbb{P}\left(S_{\mathcal I}>x_n\right)}{\overline{F}_j(x_n)},
\end{equation}
where $Q_j^{\mathcal I}(\gamma,\alpha;\Theta)\in\mathbb{R}_+$ is independent of $c$.
For $i\notin\mathcal I$ and fixed $a,b>0$,
$\frac{\overline F_i(ax_n)}
     {\overline F_j(bx_n)}
=
\frac{\overline F_i(ax_n)}
     {\overline F_i(x_n)}
\frac{\overline F_i(x_n)}
     {\overline F_j(x_n)}
\frac{\overline F_j(x_n)}
     {\overline F_j(bx_n)}
\rightarrow0$.
Therefore, by Boole's inequality, for every $\epsilon,c>0$,
\begin{equation}\label{eq:neg_frechet}
\frac{
\mathbb P\left(
\sum_{i\notin\mathcal I}X_i>\epsilon cx_n
\right)
}{
\overline F_j(cx_n)
}
\rightarrow0;
\end{equation}
note that the same limit holds trivially if $\mathcal I^c=\varnothing$, as $\mathbb P\left(
\sum_{i\notin\mathcal I}X_i>\epsilon cx_n
\right)\equiv0$ for any $\epsilon>0$.

Now fix $c>0$ and $\epsilon\in(0,1)$. Since
\[
\mathbb P(S_{\mathcal I}>cx_n)
\leq
\mathbb P(S_d>cx_n)
\leq
\mathbb P(S_{\mathcal I}>(1-\epsilon)cx_n)
+
\mathbb P\left(
\sum_{i\notin\mathcal I}X_i>\epsilon cx_n
\right),
\]
together with \eqref{eq:SI_frechet}, \eqref{eq:neg_frechet}, and regular
variation, we have:
\[
Q_j^{\mathcal I}(\gamma,\alpha;\Theta)
\leq
\liminf_{n\to\infty}
\frac{\mathbb P(S_d>cx_n)}
     {\overline F_j(cx_n)}
\leq
\limsup_{n\to\infty}
\frac{\mathbb P(S_d>cx_n)}
     {\overline F_j(cx_n)}
\leq
Q_j^{\mathcal I}(\gamma,\alpha;\Theta)(1-\epsilon)^{-\alpha}.
\]
Letting $\epsilon\downarrow0$ gives:
\begin{equation}\label{eq:S_frechet_subseq}
\frac{\mathbb P(S_d>cx_n)}
     {\overline F_j(cx_n)}
\rightarrow Q_j^{\mathcal I}(\gamma,\alpha;\Theta),
\qquad c>0.
\end{equation}
Consequently,
\[
\frac{\mathbb P(S_d>tx_n)}
     {\mathbb P(S_d>x_n)}
=
\frac{\mathbb P(S_d>tx_n)}
     {\overline F_j(tx_n)}
\frac{\overline F_j(tx_n)}
     {\overline F_j(x_n)}
\frac{\overline F_j(x_n)}
     {\mathbb P(S_d>x_n)}
\rightarrow
Q_j^{\mathcal I}(\gamma,\alpha;\Theta)t^{-\alpha}(Q_j^{\mathcal I}(\gamma,\alpha;\Theta))^{-1}
=
t^{-\alpha}.
\]
Finally, if the limit of $\frac{\mathbb P(S_d>tx)}
     {\mathbb P(S_d>x)}
\rightarrow t^{-\alpha}$ as $x\to\infty$ failed for some
$t>0$, then there would exist $\eta>0$ and a sequence $\{x_n\}\to\infty$
such that
$\left|
\frac{\mathbb P(S_d>tx_n)}
     {\mathbb P(S_d>x_n)}
-t^{-\alpha}
\right|
\geq\eta$
for every $n$, and there would exist a subsequence $\{x_{n_k}\}$ such that $\frac{\mathbb P(S_d>tx_{n_k})}
     {\mathbb P(S_d>x_{n_k})}\to t^{-\alpha}$, resulting in a contradiction, hence the limit of $\frac{\mathbb P(S_d>tx)}
     {\mathbb P(S_d>x)}
\rightarrow t^{-\alpha}$ as $x\to\infty$ must hold for any
$t>0$. Therefore,
$\overline F_{S_d}\in\operatorname{RV}_{-\alpha}$, and hence
$S_d\in\operatorname{MDA}(\Phi_\alpha)$.
\end{proof}

\subsubsection{Dominating Fr\'echet Risks}
\label{sec:frechet-other}

We now consider portfolios containing risks with different tail indices.
If a risk belongs to $\operatorname{MDA}(\Phi_\alpha)$, then risks in a
Fr\'echet MDA with a larger tail index $\alpha'>\alpha$, risks in the Gumbel MDA, and risks
with finite right endpoints are asymptotically negligible relative to it.
See, for example, Lemma 3.1 of \cite{cltyy18}. We distinguish between the
cases in which there is a single dominating Fr\'echet risk and several
dominating Fr\'echet risks.

The case of a single dominating Fr\'echet risk is already well understood.
In this case, the dependence structure is asymptotically irrelevant; see, for
example, Lemmas 3.2, 3.3, and 3.5 of \cite{cltyy18}, or
Proposition 3.5(i) of \cite{zlyxs23}. We record the result and its proof below for
completeness and for comparison with the multiple-dominating-risk case.

\begin{proposition}\label{prop:single-dominating}
Let $j^*\in\{1,\ldots,d\}$. Suppose that $X_{j^*}\in\operatorname{MDA}(\Phi_\alpha)$ for some
$\alpha>0$ and
$\frac{\overline F_i(x)}
     {\overline F_{j^*}(x)}
\rightarrow0$ as $x\to\infty$ whenever $i\neq j^*$, then V@R is asymptotically additive with
$\lim_{p\uparrow 1} D_p = 1$.
\end{proposition}

\begin{proof}
By Lemma 3.5 of \cite{cltyy18},
$
\frac{\mathbb P\left(\sum_{i=1}^dX_i>x\right)}{
\mathbb P(X_{j^*}>x)}
\rightarrow1$, hence Lemma 3.3 of \cite{cltyy18} gives
$\frac{\V@R_p(\sum_{i=1}^dX_i)}{\V@R_p(X_{j^*})}\to1$.
Moreover, by Lemma 3.2 of \cite{cltyy18},
$\frac{\V@R_p(X_i)}{\V@R_p(X_{j^*})}\to0$ for every $i\neq j^*$.
Therefore,
\[
D_p
=
\frac{
\V@R_p(S_d)
}{
\V@R_p(X_{j^*})
}
\left(
1+\sum_{i\neq j^*}
\frac{\V@R_p(X_i)}
     {\V@R_p(X_{j^*})}
\right)^{-1}
\rightarrow1.
\]
\end{proof}
In particular, the tail-dominance condition in
Proposition \ref{prop:single-dominating} holds if for every $i\neq j^*$, $X_i$: i)
belongs to $\operatorname{MDA}(\Phi_{\alpha_i})$ with
$\alpha_i>\alpha$; or ii) belongs to the Gumbel MDA; or iii) has a finite right
endpoint; see Lemma 3.1 of \cite{cltyy18}.

We next consider the case in which more than one risks share the smallest
Fr\'echet tail index. Importantly, these risks are not required to be tail-equivalent. The remaining risks can be taken as asymptotically negligible
relative to each of these several dominating risks.

\begin{theorem}[V@R aggregation with multiple dominating Fr\'echet risks]
\label{thm:multiple-dominating}
Suppose that, for some $2\leq m\leq d$,
$X_1,\ldots,X_m\in\operatorname{MDA}(\Phi_\alpha), \alpha>0$,
and that the survival copula of $(X_1,\ldots,X_m)$ is Archimedean with a
strict generator $\psi$ that is regularly varying at zero with index
$-\gamma<0$. Suppose further that
$\frac{\overline F_i(x)}{\overline F_j(x)}\rightarrow0$,
$x\to\infty$, $j\in\{1,\ldots,m\}$, $i\in\{m+1,\ldots,d\}$.
Then V@R is:
\begin{enumerate}[label=(\roman*)]
    \item asymptotically superadditive with
    $\liminf_{p\uparrow1} D_p \geq 1$, if $0<\alpha<1$; or
    \item asymptotically additive with
    $\lim_{p\uparrow1} D_p = 1$, if $\alpha=1$; or
    \item asymptotically subadditive with
    $\limsup_{p\uparrow1} D_p \leq 1$, if $\alpha>1$.
\end{enumerate}
\end{theorem}
\begin{remark}\label{remark:frechet}
    This setting is closely related to Proposition 3.5(ii) of \cite{zlyxs23}, which also reduces the asymptotic diversification behavior of a portfolio to a subvector of dominating regularly varying risks, under an assumption of stable tail dependence function (STDF) on the copula of the dominant subvector. However, the authors of \cite{zlyxs23} assumed that the dominating risks are tail-equivalent, whereas our result only requires the dominating risks to share a common Fr\'echet index, while their relative tail ratios $\frac{\overline{F}_i(\theta_k)}{\overline{F}_j(\theta_k)}$ need not converge. Besides, we need not impose a specific dependence structure on $X_{m+1},\dots,X_d$. 
\end{remark}
\begin{proof}
If $m=d$, the result follows directly from
Theorem \ref{thm:main}, so it suffices to show the statement when $m<d$. Fix $j\in\{1,\ldots,m\}$.
By applying Corollary \ref{cor:sum-frechet} to the subvector
$(X_1,\ldots,X_m)$, we have $S_m=\sum_{i=1}^m X_i\in\operatorname{MDA}(\Phi_\alpha)$. For every $i>m$, clearly $\overline{F}_j(x)\leq \mathbb{P}(S_m>x)$, thus
$\frac{\mathbb P(X_i>x)}
     {\mathbb P(S_m>x)}
\leq
\frac{\overline F_i(x)}
     {\overline F_j(x)}
\rightarrow0$.
Hence, by Lemma 3.5 of \cite{cltyy18},
$\frac{\mathbb P(S_d>x)}{\mathbb P(S_m>x)}\rightarrow1$ as $x\to\infty$.
Since $S_m\in\operatorname{MDA}(\Phi_\alpha)$, Lemma 3.3 of
\cite{cltyy18} gives
$\frac{\V@R_p(S_d)}{\V@R_p(S_m)}\rightarrow1$,
$p\uparrow1$.
Moreover, by Lemma 3.2 of \cite{cltyy18}, for every $i>m$,
$\frac{\V@R_p(X_i)}{\V@R_p(X_j)}\rightarrow0$.
Since $d$ is fixed,
$\sum_{i=m+1}^d\V@R_p(X_i)
=o\left(\sum_{i=1}^m\V@R_p(X_i)\right)$.
Consequently,
\[
D_p
=
\frac{\V@R_p(S_d)}
     {\sum_{i=1}^d\V@R_p(X_i)}
\sim
\frac{\V@R_p(S_m)}
     {\sum_{i=1}^m\V@R_p(X_i)},
\text{ as } p\uparrow1.
\]
The result now follows by applying Theorem \ref{thm:main} to the
subvector $(X_1,\ldots,X_m)$.
\end{proof}
As in the single-dominating-risk case, the tail-dominance condition is
automatically satisfied if for each $i>m$, $X_i$: i) belongs to
$\operatorname{MDA}(\Phi_{\alpha_i})$ with $\alpha_i>\alpha$; or ii) belongs
to the Gumbel MDA; or iii) has a finite right endpoint.

\subsection{Portfolios with Gumbel Risks}
\label{sec:gumbel-infinite}

\subsubsection{Infinite Right Endpoints}

We next consider the case in which all risks belong to the Gumbel MDA and
have infinite right endpoints. 
For identically distributed risks in the Gumbel MDA with an infinite right endpoint, \cite{ccyyz23} showed that V@R is asymptotically subadditive. Although their result is developed in the context of extreme-value dependence, their proof does not use any assumption on the dependence structure and therefore holds under arbitrary dependence. 
The following theorem extends this result to heterogeneous marginal risks: the risks need not be identically distributed, or even tail equivalent, provided that each belongs to the Gumbel MDA with an infinite right endpoint.

\begin{theorem}[Asymptotic subadditivity for Gumbel risks under arbitrary dependence]\label{thm:gumbel-subadditivity}
Suppose that $X_1,\ldots,X_d\geq0$ and that
$X_i\in\operatorname{MDA}(\Lambda)$ with infinite right endpoint for every
$i\in\{1,\ldots,d\}$. Then
$   \limsup_{p\uparrow1} D_p \leq 1$, i.e. V@R is asymptotically subadditive for non-heavy-tail risks.
\end{theorem}

\begin{proof}
For $p\in(0,1)$, define
$q_p:=1-\frac{1-p}{d}$.
We first show that, for every $i$,
\begin{equation}\label{eq:gumbel-quantile-ratio}
\frac{\V@R_{q_p}(X_i)}{\V@R_p(X_i)}
\rightarrow1,
\qquad p\uparrow1;
\end{equation}
indeed, write
$x_p:=\V@R_p(X_i).
$
Since $q_p\geq p$, we have $\V@R_{q_p}(X_i)\geq x_p$, while the infinite right endpoint
implies $x_p\to\infty$.
Fix $\epsilon>0$. If $\V@R_{q_p}(X_i)>(1+\epsilon)x_p$, then, by the definition of
V@R,
\[
\overline F_i((1+\epsilon)x_p)
>
1-q_p
=
\frac{1-p}{d},
\]
whereas
$
\overline F_i(x_p)\leq1-p.
$
Consequently,
\[
\frac{
\overline F_i((1+\epsilon)x_p)
}{
\overline F_i(x_p)
}
>\frac{\frac{1-p}{d}}{1-p}=
\frac1d,
\]
contradicting the rapid variation of $\overline F_i$, since $x_p\to\infty$.
Thus $\V@R_{q_p}(X_i)\leq(1+\epsilon)x_p$ for all sufficiently large $p$. Since
$\V@R_{q_p}(X_i)\geq x_p$ and $\epsilon>0$ is arbitrary,
\eqref{eq:gumbel-quantile-ratio} follows.

On the other hand, if
$\sum_{i=1}^dX_i>\sum_{i=1}^d \V@R_{q_p}(X_i)$, then $X_i> \V@R_{q_p}(X_i)$ for at least one $i\in\{1,\dots,d\}$, hence by Boole's inequality,
\begin{align*}
\mathbb P\left(
\sum_{i=1}^d X_i>\sum_{i=1}^d \V@R_{q_p}(X_i)
\right)
\leq
\sum_{i=1}^d
\mathbb P(X_i>\V@R_{q_p}(X_i))
\leq
d(1-q_p)
=
1-p.
\end{align*}
Therefore,
\[
\V@R_p\left(\sum_{i=1}^d X_i\right)
\leq 
\sum_{i=1}^d \V@R_{q_p}(X_i).
\]
It follows that
\[
\frac{
\V@R_p\left(\sum_{i=1}^d X_i\right)
}{
\sum_{i=1}^d\V@R_p(X_i)
}
\leq
\frac{
\sum_{i=1}^d\V@R_{q_p}(X_i)
}{
\sum_{i=1}^d\V@R_p(X_i)
}
\leq
\max_{1\leq i\leq d}
\frac{\V@R_{q_p}(X_i)}{\V@R_p(X_i)},
\]
where the last inequality holds due to a standard application of the mediant inequality. Since $d$ is fixed, \eqref{eq:gumbel-quantile-ratio} implies that the
right-hand side converges to $1$, and the result follows.
\end{proof}

\subsubsection{Positive Upper-Tail Dependence}
In the previous subsection, we showed that V@R is asymptotically subadditive for Gumbel risks with infinite right endpoints under arbitrary dependence. Under additional assumptions on the dependence structure, stronger asymptotic additivity results are known. In particular, under upper tail comonotonicity, \cite{cltyy18} established asymptotic additivity for a broad class of heterogeneous risks, including risks in the Gumbel MDA, without requiring tail equivalence. Their dependence condition, however, is substantially stronger than positive diagonal upper-tail dependence and, in particular, implies that the diagonal upper-tail dependence coefficient equals one. More recently, Proposition 4.1 of \cite{zlyxs23} established asymptotic additivity under a positive diagonal upper-tail dependence coefficient, but for continuous log-regularly varying marginal distributions sharing a common log-regular-variation index.

The following theorem retains the weaker dependence condition of \cite{zlyxs23} while allowing heterogeneous marginals throughout the Gumbel MDA with infinite right endpoints. In particular, no common log-regular-variation index, tail-equivalence condition, or continuity assumption is required.

\begin{theorem}[Asymptotic additivity under positive upper-tail dependence for infinite Gumbel tail]
\label{thm:gumbel-tail-dependent}
Suppose that $X_1,\ldots,X_d\geq0$, and
$
X_i\in\operatorname{MDA}(\Lambda)
$
with infinite right endpoint for every $i\in\{1,\ldots,d\}$. Suppose that there exists a survival copula $\widehat C$ of
$(X_1,\ldots,X_d)$ for which the diagonal upper-tail dependence coefficient
\[
\lambda_\Delta
:=
\lim_{u\downarrow0}
\frac{\widehat C(u,\ldots,u)}{u}
\]
exists and is strictly positive. Then V@R is asymptotically additive with
$\lim_{p \uparrow 1} D_p = 1$.
\end{theorem}

\begin{proof}
Since $\widehat C$ is a copula,
$
\widehat C(u,\ldots,u)\leq u,
u\in[0,1],
$
so $0<\lambda_\Delta\leq1$.
Fix $c>\lambda_\Delta^{-1}$, and set 
$
q_p:=1-c(1-p)
$ for $p$ sufficiently close to 1.
For each $i=1,\dots,d$, write
$
x_{i,p}:=\V@R_p(X_i),
y_{i,p}:=\V@R_{q_p}(X_i).
$
Clearly $q_p<p$, and we have
$
y_{i,p}\leq x_{i,p},
 i=1,\ldots,d.
$
By Lemma \ref{lemma:tail-at-var}, as $p\uparrow1$,
\begin{equation}\label{eq:gumbel-tail-at-var}
\frac{\overline F_i(x_{i,p})}{1-p}
\rightarrow1,
\quad
\frac{\overline F_i(y_{i,p})}{c(1-p)}
\rightarrow1.
\end{equation}
Consequently,
\begin{equation}\label{eq:gumbel-tail-ratio-xy}
\frac{\overline F_i(y_{i,p})}
     {\overline F_i(x_{i,p})}
\rightarrow c.
\end{equation}

We claim that
\begin{equation}\label{eq:gumbel-quantile-equivalence}
\frac{y_{i,p}}{x_{i,p}}
\rightarrow1,
\qquad p\uparrow1.
\end{equation}
Suppose that \eqref{eq:gumbel-quantile-equivalence} does not hold. Since $y_{i,p}\leq x_{i,p}$, there would exist
$\epsilon\in(0,1)$ and a sequence $\{p_n\}\uparrow1$ such that
$
y_{i,p_n}\leq(1-\epsilon)x_{i,p_n}
$
for every $n$. By monotonicity of the survival function,
\[
\frac{\overline F_i(y_{i,p_n})}
     {\overline F_i(x_{i,p_n})}
\geq
\frac{\overline F_i((1-\epsilon)x_{i,p_n})}
     {\overline F_i(x_{i,p_n})}\to\infty,
\]
since $\overline F_i$ is rapidly varying and
$x_{i,p_n}\to\infty$, yet this
contradicts \eqref{eq:gumbel-tail-ratio-xy}. Therefore,
\eqref{eq:gumbel-quantile-equivalence} holds.

Set
$
u_p:=c(1-p),
u_{i,p}:=\overline F_i(y_{i,p}).
$
By \eqref{eq:gumbel-tail-at-var}, $\frac{u_{i,p}}{u_p}\to1$ as $p\uparrow1$ for every $i$.
Using the Lipschitz property of copulae,
\[
\left|
\widehat C(u_{1,p},\ldots,u_{d,p})
-
\widehat C(u_p,\ldots,u_p)
\right|
\leq
\sum_{i=1}^d|u_{i,p}-u_p|=u_p\sum_{i=1}^d\left|\frac{u_{i,p}}{u_p}-1\right|=(1-p)c\sum_{i=1}^d\left|\frac{u_{i,p}}{u_p}-1\right|
=
o(1-p).
\]
Therefore, since $1-p=\frac{u_p}{c}$,
$$
\frac{
\mathbb P\left(
X_i>y_{i,p},\ i=1,\ldots,d
\right)
}{
1-p
}
=
\frac{
\widehat C(u_{1,p},\ldots,u_{d,p})
}{
1-p
}=
\frac{
\widehat C(u_p,\ldots,u_p)
}{
1-p
}
+o(1)
\rightarrow
c\lambda_\Delta
>1.
$$
Hence, for all sufficiently large $p$,
$
\mathbb P\left(
X_i>y_{i,p},\ i=1,\ldots,d
\right)>1-p.
$
Since
$\{X_i>y_{i,p},\ i=1,\ldots,d\}
\subset
\left\{
\sum_{i=1}^dX_i>
\sum_{i=1}^dy_{i,p}
\right\}$,
we obtain $\mathbb P\left(
\sum_{i=1}^dX_i>
\sum_{i=1}^dy_{i,p}
\right)>1-p$ for all sufficiently large $p$.

By the definition of V@R,
$\V@R_p(S_d)>
\sum_{i=1}^dy_{i,p}.
$
Consequently,
\begin{equation}\label{eq:gumbel-lower}
\liminf_{p\uparrow1}
\frac{
\V@R_p\left(\sum_{i=1}^dX_i\right)
}{
\sum_{i=1}^d\V@R_p(X_i)
}
\geq
\liminf_{p\uparrow1}
\frac{\sum_{i=1}^dy_{i,p}}
     {\sum_{i=1}^dx_{i,p}}
=1,
\end{equation}
where the last equality follows from
\eqref{eq:gumbel-quantile-equivalence}. 
On the other hand, by Theorem \ref{thm:gumbel-subadditivity},
\begin{equation}\label{eq:gumbel-upper}
\limsup_{p\uparrow1}
\frac{
\V@R_p\left(\sum_{i=1}^dX_i\right)
}{
\sum_{i=1}^d\V@R_p(X_i)
}
\leq1.
\end{equation}
The desired result then follows by combining \eqref{eq:gumbel-lower} and \eqref{eq:gumbel-upper}.
\end{proof}
\begin{corollary}[Infinite Gumbel risks under Archimedean survival dependence]
\label{cor:gumbel-archimedean}
Suppose that Assumption \ref{AssumeArchi} holds and that
$X_i\in\operatorname{MDA}(\Lambda)$
with infinite right endpoint for every $i\in\{1,\ldots,d\}$. Then
$\lim_{p \uparrow 1} D_p = 1$.
\end{corollary}
\begin{proof}
By Theorem 3.1 of \cite{charpentier2009tails},
$\lim_{u\downarrow0}\frac{\widehat C(u,\ldots,u)}{u}
=d^{-1/\gamma}>0$
under Assumption \ref{AssumeArchi}, hence
Theorem \ref{thm:gumbel-tail-dependent} applies.
\end{proof}

\subsubsection{Finite-Endpoint Additions}
\label{sec:gumbel-mixed}

We next allow some of the risks to have finite right endpoints. These risks
are asymptotically negligible relative to those in the Gumbel MDA with
infinite right endpoints.

\begin{theorem}[Asymptotic subadditivity with finite-endpoint additions]\label{thm:gumbel-mixed}
Suppose that $X_1,\ldots,X_d\geq0$. Suppose that
$X_i\in\operatorname{MDA}(\Lambda)$ with infinite right endpoint for every
$i\in\mathcal I$, where $\mathcal I$ is a nonempty subset of $\{1,\ldots,d\}$, while $X_i$ has a finite right endpoint whenever
$i\notin\mathcal I$, then V@R is asymptotically subadditive with
$\limsup_{p \uparrow 1} D_p \leq 1$.
\end{theorem}

\begin{proof}
Denote $S_{\mathcal I}=\sum_{i\in\mathcal I}X_i$ as before, and also
$b:=\sum_{i\notin\mathcal I}x_{F_i}<\infty$. Since
$\sum_{i\notin\mathcal I}X_i\leq b$ almost surely, monotonicity and
translation invariance of V@R give:
\[
\V@R_p\left(\sum_{i=1}^dX_i\right)
\leq
\V@R_p(S_{\mathcal I})+b.
\]
Moreover, since $\mathcal I\neq\varnothing$, and for every $i\in\mathcal I$, $X_i$ has an infinite right endpoint, then
$\sum_{i\in\mathcal I}\V@R_p(X_i)\to\infty$ as $p\uparrow1$. Hence,
\[
\limsup_{p\uparrow1}
\frac{\V@R_p\left(\sum_{i=1}^dX_i\right)}
     {\sum_{i=1}^d\V@R_p(X_i)}
\leq
\limsup_{p\uparrow1}
\frac{\V@R_p(S_{\mathcal I})+b}
     {\sum_{i\in\mathcal I}\V@R_p(X_i)}
\leq1,
\]
where the last inequality follows from
Theorem \ref{thm:gumbel-subadditivity}.
\end{proof}

\subsection{Portfolios of Risks with Finite Right Endpoints}
\label{sec:finite-endpoints}

We finally consider the case in which all marginal risks have finite right endpoints when the survival copula is Archimedean. The result follows directly from Proposition 4.8 of \cite{cltyy18}.

\begin{theorem}[Aggregation of risks with finite right endpoints]
\label{thm:finite-endpoints}
Suppose that the survival copula of $(X_1,\ldots,X_d)$ is Archimedean
with a strict generator $\psi$, and suppose that
$x_{F_i}<\infty$ for every $i\in\{1,\ldots,d\}$, then
\[
\lim_{p\uparrow1}
\V@R_p\left(\sum_{i=1}^d X_i\right)
=
\lim_{p\uparrow1}
\sum_{i=1}^d\V@R_p(X_i)
=
\sum_{i=1}^d x_{F_i}.
\]
\end{theorem}

\begin{proof}
For every $\epsilon\in(0,1)$,
$\widehat C(\epsilon,\ldots,\epsilon)
=
\psi^{-1}\left(d\cdot\psi(\epsilon)\right)>0$, and the strictness of $\psi$ gives $\psi(\epsilon)<\infty$ and
$\psi^{-1}(t)>0$ for every finite $t$. The result now follows from
Proposition 4.8 of \cite{cltyy18}.
\end{proof}

\section{Numerical Illustrations}\label{sec:num}

In this section, we provide pictorial illustrations that share the feeling about the convergence rate for our claimed asymptotic results established in
Section \ref{sect:var_aa}. The
quantity displayed in Figures \ref{fig:num-frechet}--\ref{fig:num-gumbel} is
primarily the diversification ratio $D_p$.
Values below, equal to, and above 1 correspond to the respective cases of
subadditivity, additivity, and superadditivity of V@R at the confidence level $p$. Throughout the simulation studies, unless stated otherwise, we assume the dependence structure of a survival Clayton
copula with parameter $\gamma>0$ and generator
$\psi_\gamma(u)=u^{-\gamma}-1$.  We use its frailty representation (also see Chapter 4 of \cite{j97} and Chapter 4 of \cite{n06}): we let
\[
U_i=\left(1+\frac{E_i}{V}\right)^{-\frac{1}{\gamma}},
\qquad
X_i=F_i^{\leftarrow}(1-U_i),
\qquad i=1,\ldots,d,
\]
where $V\sim\operatorname{Gamma}\left(\frac{1}{\gamma},1\right)$, $E_1,\ldots,E_d\stackrel{\mathrm{iid}}{\sim}\operatorname{Exp}(1)$ and are independent of $V$, then $(U_1,\dots,U_d)$ has the Clayton copula; indeed, this construction
assigns the Clayton copula to the survival probabilities of the marginal
risks.  Each experiment uses Monte Carlo sampling with size $10^7$.  The confidence levels
range from $0.95$ to $0.9995$ over 50 values that are equally spaced on the
$-\log_{10}(1-p)$ scale.  At the largest confidence level, the simulation
therefore contains approximately $5{,}000$ observations above the target
quantile.  The marginal V@Rs in the denominator of
\eqref{eq:diversification-ratio} are evaluated from their exact quantile
functions without smoothing; only the V@R of the aggregate loss is estimated empirically, using
the sample generalized inverse.  

\subsection{A common Fr\'echet index without tail equivalence}

Our first experiment demonstrates that
Theorem \ref{thm:main} does not require tail equivalence. For
$\alpha>0$ and $
|a|<\frac1\alpha$,
define a distribution through its upper quantile function by
\begin{equation}\label{eq:num-oscillatory-quantile}
F_{\alpha,a}^{\leftarrow}(1-u)
=
u^{-\frac{1}{\alpha}}
\exp\left\{
a\sin\bigl(\ln(1-\ln u)\bigr)
\right\}-1,
\qquad 0<u<1;
\end{equation}
indeed, to see that it is a valid upper quantile function, we let
$t=-\ln u$, then
\begin{align*}
\frac{d}{dt}
\ln\left\{
F_{\alpha,a}^{\leftarrow}(1-e^{-t})+1
\right\}
=
\frac1\alpha
+
\frac{
a\cos\bigl(\ln(1+t)\bigr)
}{
1+t
}\geq
\frac1\alpha-|a|
>0.
\end{align*}
Hence the right-hand side of
\eqref{eq:num-oscillatory-quantile} is increasing as $u\downarrow0$. Moreover, the factor
\[
L_a(t)
=
\exp\left\{
a\sin\bigl(\ln(1+\ln t)\bigr)
\right\},
\qquad t\geq1,
\]
is clearly slowly varying, and so
$F_{\alpha,a}^{\leftarrow}(1-u)$ is regularly varying at zero with index
$-\frac{1}{\alpha}$, and the corresponding distribution belongs to
$\operatorname{MDA}(\Phi_\alpha)$. On the other hand, for $a\neq b\in\mathbb{R}$,
\[
\frac{F_{\alpha,a}^{\leftarrow}(1-u)+1}
     {F_{\alpha,b}^{\leftarrow}(1-u)+1}
=
\exp\left\{
(a-b)\sin\bigl(\ln(1-\ln u)\bigr)
\right\},
\]
which has no limit as $u\downarrow0$. As a result, the corresponding marginal distributions are not
tail-equivalent.

\begin{figure}[t]
    \centering
    \includegraphics[width=\textwidth]{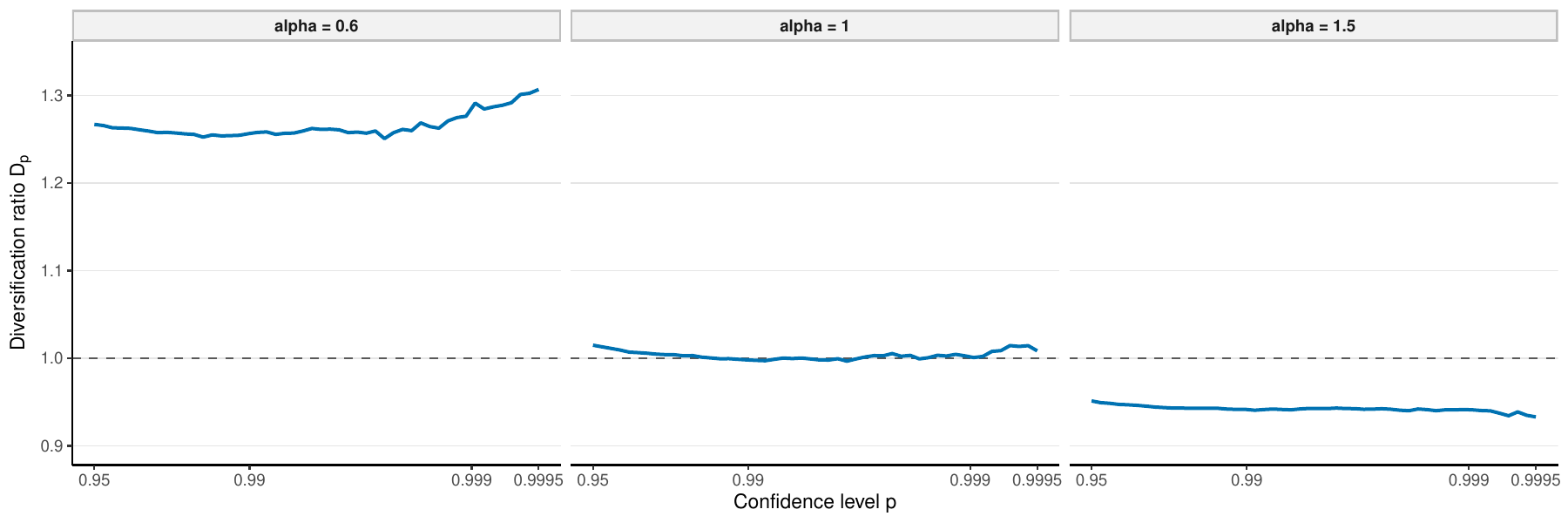}
    \caption{Diversification ratio $D_p$ for three non-tail-equivalent
    Fr\'echet marginals with a common tail index.  The marginal upper quantiles
    are given by \eqref{eq:num-oscillatory-quantile} with
    $a\in\{-0.35,0,0.35\}$, and the survival copula is Clayton with
    $\gamma=1$.  The dashed horizontal line is $D_p=1$.}
    \label{fig:num-frechet}
\end{figure}

We take three marginals of the form \eqref{eq:num-oscillatory-quantile} with
$a=-0.35,0,0.35$, model their dependence structure by a survival Clayton copula with $\gamma=1$,
and consider $\alpha=0.6,1$, and $1.5$.  Figure \ref{fig:num-frechet} shows
the three regimes predicted by Theorem \ref{thm:main}.  For $\alpha=0.6$,
the diversification ratio remains above 1 at the far tail, displaying
superadditivity.  For $\alpha=1.5$, it remains below 1, displaying
subadditivity.  At the boundary $\alpha=1$, the ratio approaches 1 despite
the persistent oscillation in the relative marginal tail scales.  

\subsection{Dependence strength and the number of dominating risks}

\begin{figure}[t]
    \centering
    \includegraphics[width=1\textwidth]{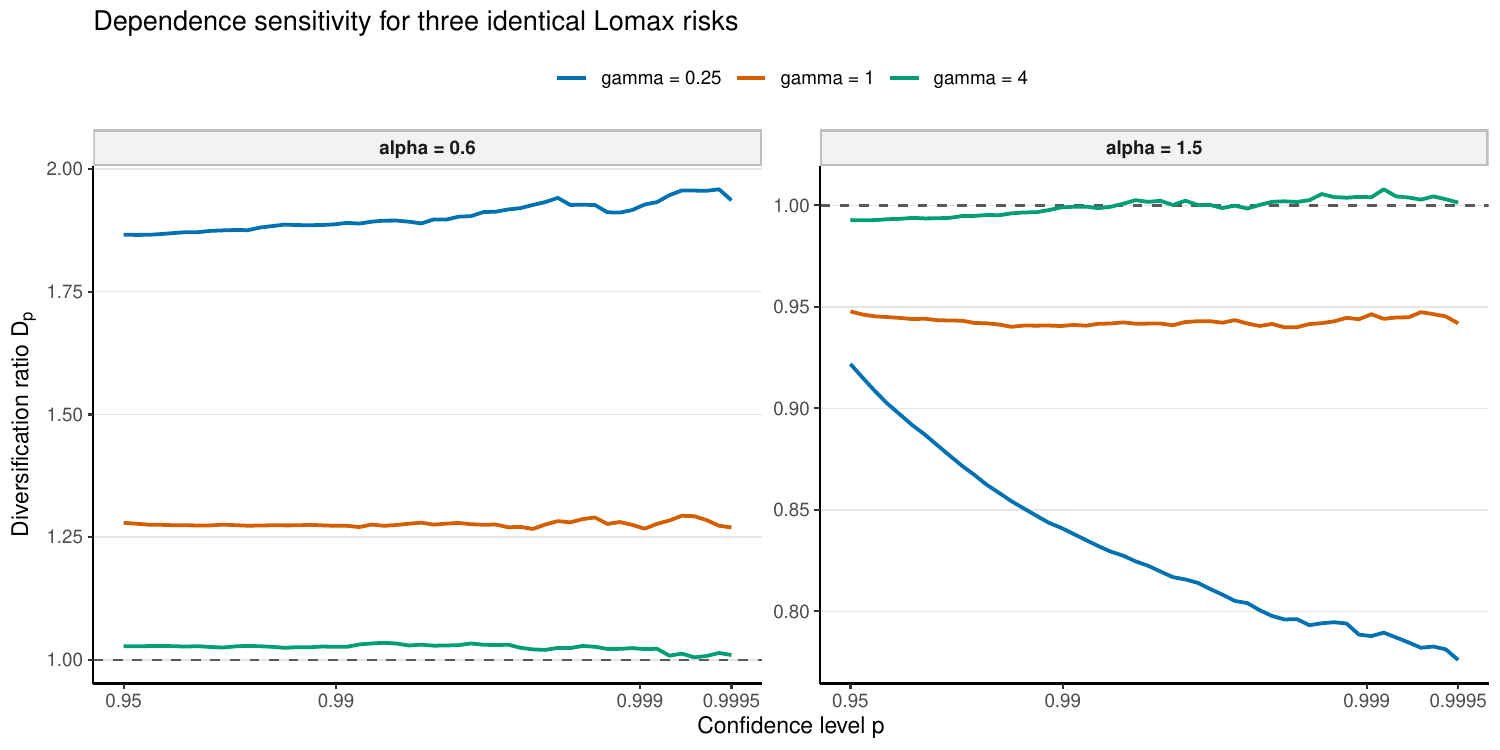}
    \caption{Effect of dependence strength on the diversification ratio.
    The portfolio consists of three identically distributed Lomax risks under
    a survival Clayton copula with $\gamma=0.25,1$, and $4$. The two panels
    correspond to $\alpha=0.6$ and $\alpha=1.5$, respectively. The dashed
    horizontal line is $D_p=1$.}
    \label{fig:num-dependence}
\end{figure}

We next discuss the effects of two theoretically distinct aspects caused by:
i) the regular-variation index of the Archimedean generator; and ii) the number of
dominating risks at the far tail. Throughout this subsection, $L_\alpha$
denotes a Lomax (generalized Pareto) risk parametrized by $\alpha$, with
\[
\mathbb P(L_\alpha>x)=(1+x)^{-\alpha},
\qquad x\geq0.
\]
We first examine the effect of the strength of dependence modelling among individual risks. In
Figure \ref{fig:num-dependence}, the portfolio consists of three identically
distributed copies of $L_\alpha$, with survival Clayton parameters
$\gamma=0.25,1$, and $4$. When $\alpha=0.6$, the diversification ratios are
above 1 and decrease toward the comonotonic benchmark 1 as $\gamma$
increases. When $\alpha=1.5$, the ratios are below 1 and increase toward
1. This provides a finite-level illustration of the monotonicity in
Theorem \ref{thm:qin}, namely a stronger Clayton tail dependence reduces the magnitude
of the diversification effect.

\begin{figure}[t]
    \centering
    \includegraphics[width=0.55\textwidth]{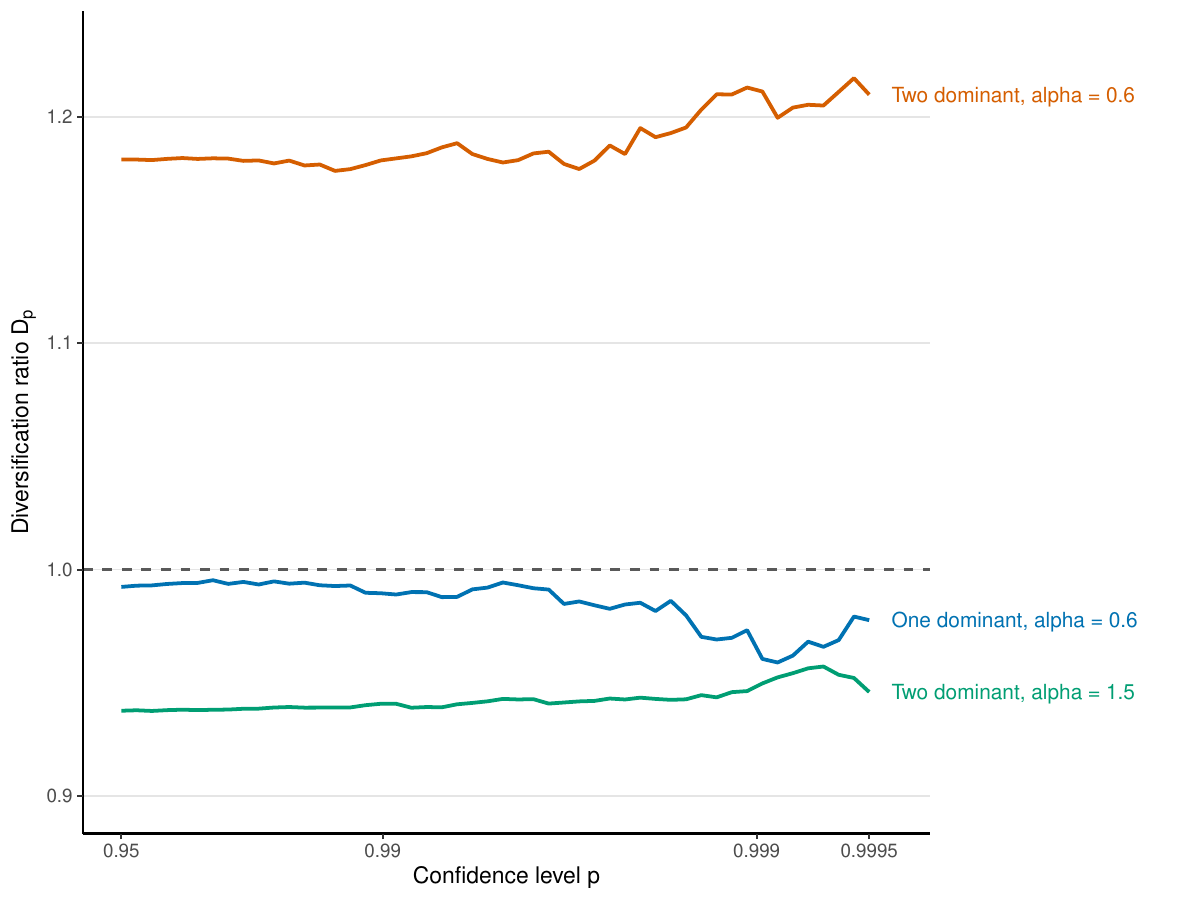}
    \caption{Effect of the number of dominating risks on the diversification
    ratio under a survival Clayton copula with $\gamma=1$.
    Portfolio $\mathcal P_1$ has one dominating Fr\'echet risk, whereas
    $\mathcal P_2$ and $\mathcal P_3$ each have two dominating Fr\'echet risks
    with common tail indices $0.6$ and $1.5$, respectively. The dashed
    horizontal line is $D_p=1$.}
    \label{fig:num-dominance}
\end{figure}

We next examine the role of the number of dominating risks.
Figure \ref{fig:num-dominance} uses a survival Clayton copula with
$\gamma=1$ and compares the following three portfolios:
\begin{align*}
\mathcal P_1
&=L_{0.6}+L_{2.5}+E_{0.25}+U_{0.5},\\
\mathcal P_2
&=L_{0.6}+0.7L_{0.6}+L_{2.5}+E_{0.25}+U_{0.5},\\
\mathcal P_3
&=L_{1.5}+0.7L_{1.5}+L_{3.5}+E_{0.25}+U_{0.5},
\end{align*}
where $E_s$ is an exponential random variable with mean $s$ and $U_b$ is another uniform random variable on $[0,b]$.
The first portfolio $\mathcal{P}_1$ has a single dominating Fr\'echet risk $L_{0.6}$. Its
diversification ratio approaches 1 even though the dominating risk $L_{0.6}$ has
infinite mean, in agreement with
Proposition \ref{prop:single-dominating}. The second portfolio $\mathcal{P}_2$ has two dominating
risks with common index $0.6<1$ and exhibits asymptotic superadditivity,
whereas the third $\mathcal{P}_3$ has two dominating risks with common index $1.5>1$ and
exhibits asymptotic subadditivity; these latter two cases illustrate
Theorem \ref{thm:multiple-dominating}. Together, the three curves show why the
number of asymptotically relevant risks matters; indeed, the presence of an
infinite-mean Fr\'echet risk alone does not produce a strict superadditivity
effect when all remaining risks are asymptotically negligible.

\subsection{Gumbel risks and finite right endpoints}
\begin{figure}[t]
    \centering

    \begin{subfigure}[t]{0.49\textwidth}
        \centering
        \includegraphics[width=\textwidth]{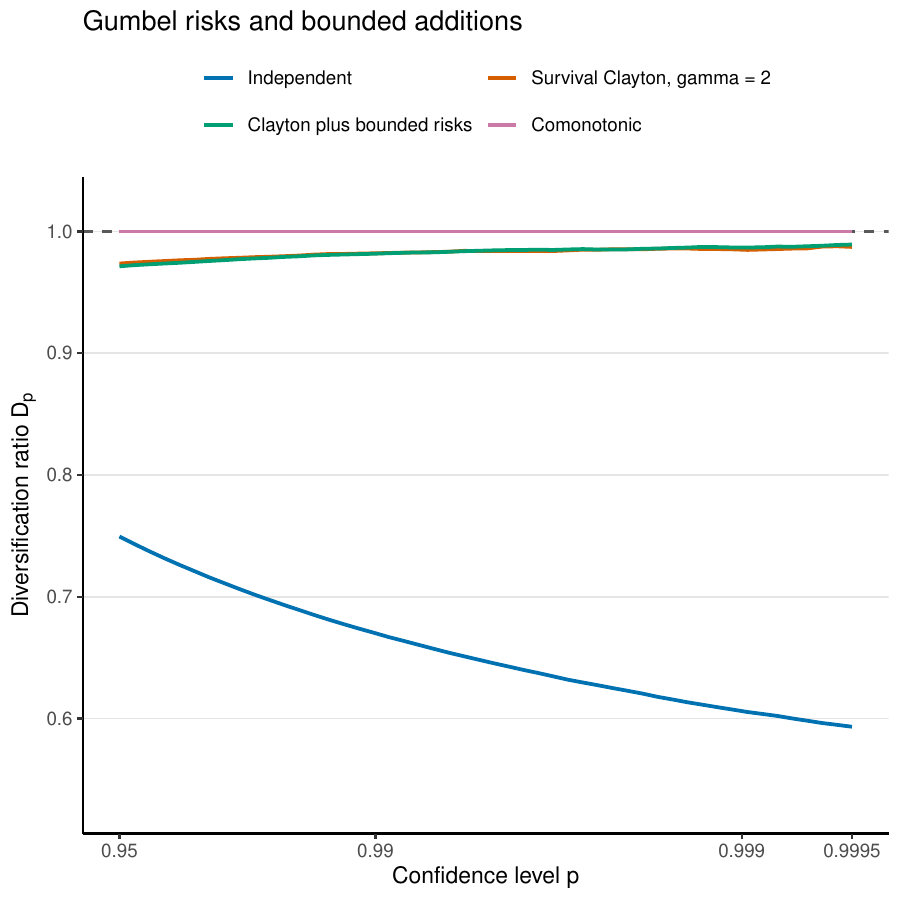}
        \caption{Gumbel risks and bounded additions.}
    \end{subfigure}
    \hfill
    \begin{subfigure}[t]{0.49\textwidth}
        \centering
        \includegraphics[width=\textwidth]{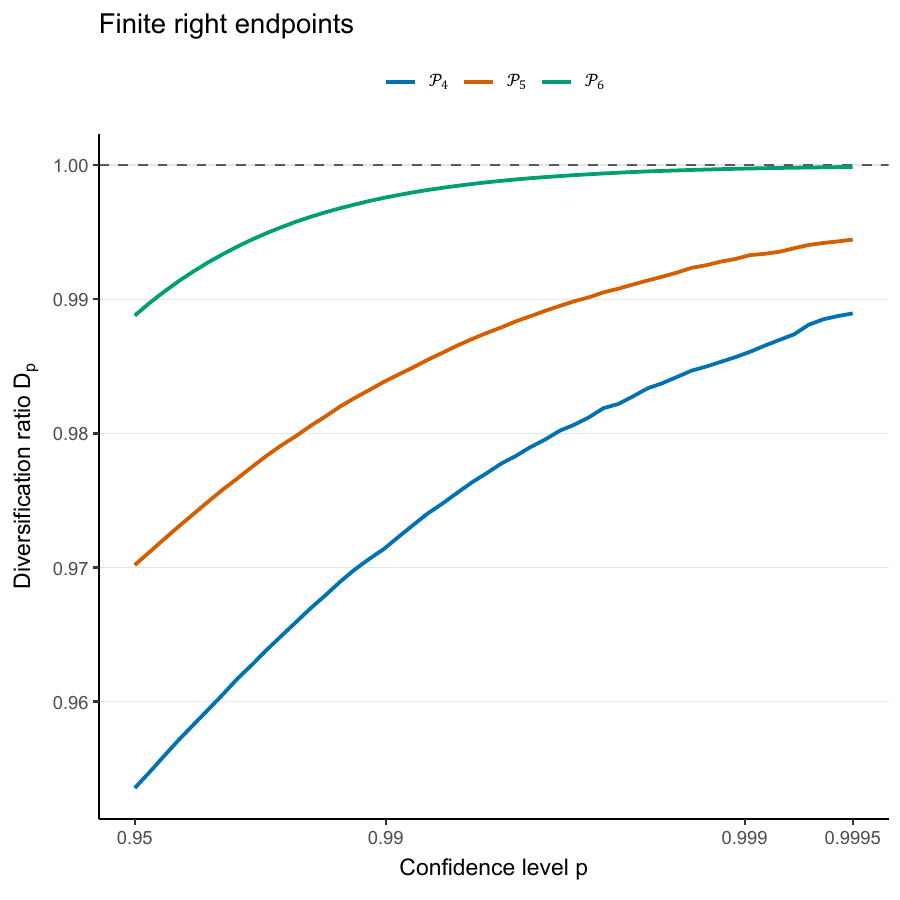}
        \caption{Risks with finite right endpoints.}
    \end{subfigure}

    \caption{Gumbel risks and finite right endpoints. Panel (a) compares independence, survival Clayton dependence with $\gamma=2$,
and comonotonicity for the same heterogeneous Gumbel portfolio, and also shows
the effect of adding two bounded risks under the same survival Clayton
dependence. Panel (b) shows the
    diversification ratios for three portfolios whose marginal risks all have
    finite right endpoints. Dashed horizontal lines indicate $D_p=1$.}
    \label{fig:num-gumbel}
\end{figure}

For Figure \ref{fig:num-gumbel}(a), consider the following unbounded Gumbel portfolio:
\[
S_G=E_1+W_{2,2}+Z,
\]
where $E_1$ is an exponential random variable with mean 1, $W_{2,2}$ is a random variable with tail distribution
$\mathbb P(W_{2,2}>x)=\exp\left\{-\left(\frac{x}{2}\right)^2\right\}$, and
$Z$ is another lognormal random variable with log-scale parameters $(0,0.6)$. We compare three
dependence structures for these same marginal risks, namely i) independence; ii) a survival
Clayton copula with $\gamma=2$; and iii) comonotonicity, together with another portfolio iv) obtained by adding two bounded risks to $S_G$ under ii), namely $U_1+2B$, where $U_1$ is a uniform random variable on $[0,1]$ and the Beta random variable
$B\sim\operatorname{Beta}(2,4)$.

The independent curve remains below 1 over the displayed confidence levels,
consistent with the asymptotic subadditivity result in
Theorem \ref{thm:gumbel-subadditivity}. In contrast, the diversification ratio
under survival Clayton dependence moves close to 1 as the confidence level
increases, illustrating the asymptotic additivity established in
Corollary \ref{cor:gumbel-archimedean}. Here $\gamma=2$ is used so that this
convergence is visible over the displayed range while the dependence remains
distinct from comonotonicity. The comonotonic ratio is identically 1.
Finally, the survival-Clayton curves with and without the two bounded risks
become close at high confidence levels, illustrating the asymptotic minority of the addition of bounded random variables underlying
Theorem \ref{thm:gumbel-mixed}.

Figure \ref{fig:num-gumbel}(b) considers three portfolios consisting entirely
of risks with finite right endpoints. In each case, the marginal endpoints are
$1$, $2$, and $3$, respectively, so that their sum is 6. Specifically, we consider the following portfolios:
\begin{align*}
\mathcal P_4
&=B_{1,5}+2B_{2,4}+3B_{3,3},\\
\mathcal P_5
&=B_{2,5}+2B_{3,3}+3B_{5,2},\\
\mathcal P_6
&=B_{2,1}+2B_{3,1}+3B_{5,1},
\end{align*}
where the Beta random variables $B_{a,b}\sim\operatorname{Beta}(a,b)$. All three portfolios are coupled
by a survival Clayton copula with $\gamma=1$. Although Theorem \ref{thm:finite-endpoints} is stated in terms of convergence
to the common endpoint sum rather than directly in ratio form, it immediately
implies $D_p\to1$, since both the numerator and denominator converge to 6.
The three examples illustrate this convergence for different marginal
behavior near their finite right endpoints.

\section{Discussion}
\label{sec:discussion}

This article develops asymptotic theory for V@R
aggregation under marginal tail heterogeneity. The main results for the case of Fr\'echet MDA
 show that global tail equivalence is not required for the familiar asymptotic superadditive or subadditive results depending on whether $\alpha<1$ or $\alpha>1$, respectively. The key argument is a subsequential one of the relative
marginal tail scales. Despite the possible dependence of the aggregate-tail constant on the composition of the limiting tail, its
asymptotic sub/superadditivity still solely depends only on whether $\alpha$ is above, equal, or below 1. This sequence-independent sign
property yields the global limsup and liminf asymptotic classification for V@R. The Archimedean survival-dependence framework also permits a detailed analysis of how dependence affects aggregation. For fixed positive limiting marginal-tail ratios, the explicit
representation of the aggregate-tail constant $Q_j(\gamma,\alpha;\Theta)$ yields strict monotonicity
in the generator parameter $\gamma$ when $\alpha\ne1$ and convergence
to the corresponding comonotonic limit $Q_{j,{\rm com}}(\alpha;\Theta)$ as $\gamma\to\infty$.

The tail-index classification can possibly be extended to more general dependence structures without imposing global tail equivalence. To motivate the idea, inspired by \cite{zlyxs23}, let $C$ denote the ordinary copula of
$(X_1,\ldots,X_d)$ with Lipschitz property, and suppose that for $\mathbf z\in[0,\infty)^d$, the STDF (see Remark \ref{remark:frechet}) $\ell(\mathbf z)
:=
\lim_{p\uparrow1}
\frac{1-C(1-(1-p)z_1,\ldots,1-(1-p)z_d)}{1-p}$ exists, then it admits a spectral representation (see \cite{segers2012max}), namely $\ell(\mathbf z)
=
\mathbb E\left[\max_{i\in\{1,\dots,d\}}(z_iW_i)\right]$
 for some non-negative random variables $W_i$ with $\mathbb{E}[W_i]=1$,  for $i=1,\dots,d$.
Assuming regularly varying marginal tails with a common index $-\alpha<0$, we anticipate the following result:
\begin{equation}\label{eq:discussion_Dp}
D_p
=\frac{K_\alpha\left(\V@R_p(X_1),\ldots,\V@R_p(X_d)\right)}
       {\sum_{i=1}^d \V@R_p(X_i)}
+o(1),\qquad p\uparrow1.    
\end{equation}
where $K_\alpha(\mathbf z)
:=
\left\{
\mathbb E\left[
\left(\sum_{i=1}^d z_iW_i^{\frac{1}{\alpha}}\right)^\alpha
\right]
\right\}^{\frac{1}{\alpha}}$ for $\mathbf z\in[0,\infty)^d$; this should not require the convergence of the ratios $\frac{\overline F_i(x)}{\overline F_j(x)}$ as $x\to\infty$.

Some approach to obtain \eqref{eq:discussion_Dp} can be understood as follows. Recall that $x_{i,p}=\V@R_p(X_i)$, and denote
$\mathbf Y_p:=(Y_{1,p},\ldots,Y_{d,p})^\top$, where $Y_{i,p}:=\frac{X_i}{x_{i,p}}$, $i=1,\dots,d$. Fix $\mathbf y=(y_1,\ldots,y_d)^\top\in\mathbb{R}_+^d$, and let $U_{\mathbf y}:=\{\mathbf x\in[0,\infty)^d:x_i>y_i\text{ for some }i\}$, then
$\mathbb P(\mathbf Y_p\in U_{\mathbf y})
=1-C\left(F_1(x_{1,p}y_1),\ldots,F_d(x_{d,p}y_d)\right)$.
Note that $\frac{\overline F_i(x_{i,p})}{1-p}\to1$ as $p\uparrow1$ by Lemma \ref{lemma:tail-at-var}, hence $z_{i,p}:=\frac{\overline F_i(x_{i,p}y_i)}{1-p}
=\frac{\overline F_i(x_{i,p}y_i)}{\overline F_i(x_{i,p})}
 \frac{\overline F_i(x_{i,p})}{1-p}
\to y_i^{-\alpha}$ for $i=1,\dots,d$.
Using the Lipschitz property of $C$, and noting that $F_i(x_{i,p}y_i)=1-(1-p)z_{i,p}$, we obtain:
\begin{align*}
    &\ \left|\frac{\mathbb P(\mathbf Y_p\in U_{\mathbf y})}{1-p}-\frac{1-C(1-(1-p)y_1^{-\alpha},\dots,1-(1-p)y_d^{-\alpha})}{1-p}\right|\\
    =&\ \left|\frac{1-C(1-(1-p)z_{1,p},\dots,1-(1-p)z_{d,p})}{1-p}-\frac{1-C(1-(1-p)y_1^{-\alpha},\dots,1-(1-p)y_d^{-\alpha})}{1-p}\right|\\
    \leq&\  \sum_{i=1}^d |z_{i,p}-y_i^{-\alpha}|\to0,\text{ as }p\uparrow1;
\end{align*}
together with the definition of STDF $\ell(\cdot)$ and a standard application of triangular inequality, we anticipate that $\frac{\mathbb P(\mathbf Y_p\in U_{\mathbf y})}{1-p}
\to
\ell(y_1^{-\alpha},\ldots,y_d^{-\alpha})
=\mathbb E\left[\max_{i\in\{1,\dots,d\}}\frac{W_i}{y_i^\alpha}\right]$.

Next, define a new measure $\nu(B):=\mathbb E\left[\int_0^\infty
\mathbbm{1}_{\{(sW_1^{\frac{1}{\alpha}},\ldots,sW_d^{\frac{1}{\alpha}})\in B\}} \alpha s^{-\alpha-1} ds\right]$ for any Borel set $B\subseteq[0,\infty)^d$. In particular, when $B=U_{\mathbf y}$, the indicator equals 1 when $sW_i^{\frac{1}{\alpha}}>y_i$ or $s^{\alpha}\frac{W_i}{y_i^{\alpha}}>1$ for some $i=1,\dots,d$, i.e., when $s^\alpha\max_{i\in\{1,\dots,d\}}\frac{W_i}{y_i^\alpha}>1$, hence
\begin{equation}\label{eq:nu_u}
\nu(U_{\mathbf y})=\mathbb E\left[\int_{\left(\max_{i\in\{1,\dots,d\}}\frac{W_i}{y_i^\alpha}\right)^{-\frac{1}{\alpha}}}^\infty
\alpha s^{-\alpha-1} ds\right]=\mathbb E\left[-s^{-\alpha}\Bigg|_{\left(\max_{i\in\{1,\dots,d\}}\frac{W_i}{y_i^\alpha}\right)^{-\frac{1}{\alpha}}}^\infty\right]=\mathbb E\left[\max_{i\in\{1,\dots,d\}}\frac{W_i}{y_i^\alpha}\right],
\end{equation}
from which we also have $\frac{\mathbb P(\mathbf Y_p\in U_{\mathbf y})}{1-p}\to\nu(U_{\mathbf y})$.

We then fix $t>0$, and let $\mathbf 0_d$, $\mathbf 1_d$ denote the respective $d$-dimensional vectors of 0's and 1's, then
$U_{t\cdot\mathbf 1_d}:=\{\mathbf x\in[0,\infty)^d:x_i>t\text{ for some }i\}=\{\mathbf x\in[0,\infty)^d:\max_{i=1,\dots,d} x_i>t\}$.
Note that $\{\mathbf Y_p\leq\mathbf y\}\cap\{\mathbf Y_p\in U_{t\cdot\mathbf 1_d}\}=\{\mathbf Y_p\in[\mathbf 0_d,\mathbf y]\cap U_{t\cdot\mathbf 1_d}\}=\{\mathbf Y_p\in U_{\mathbf y\wedge t\cdot\mathbf 1_d}\}\setminus\{\mathbf Y_p\in U_{\mathbf y}\}$,
where $U_{\mathbf y\wedge t\cdot\mathbf 1_d}:=\{\mathbf x\in[0,\infty)^d:x_i>(y_i\wedge t)\text{ for some }i\}$ analogously, hence by using the convergence obtained above, we have the following limiting joint conditional distribution of $\mathbf Y_p$ as $p\uparrow1$:
\begin{align*}
\mathbb P(\mathbf Y_p\le\mathbf y\mid\mathbf Y_p\in U_{t\cdot\mathbf 1_d})&=\frac{\mathbb P(\{\mathbf Y_p\leq\mathbf y\}\cap\{\mathbf Y_p\in U_{t\cdot\mathbf 1_d}\})}
      {\mathbb P(\mathbf Y_p\in U_{t\cdot\mathbf 1_d})}
=\frac{\frac{\mathbb P(\mathbf Y_p\in U_{\mathbf y\wedge t\cdot\mathbf 1_d})}{1-p}
       -\frac{\mathbb P(\mathbf Y_p\in U_{\mathbf y})}{1-p}}
      {\frac{\mathbb P(\mathbf Y_p\in U_{t\cdot\mathbf 1_d})}{1-p}}\\
      &\to\frac{\nu(U_{\mathbf y\wedge t\cdot\mathbf 1_d})-\nu(U_{\mathbf y})}{\nu(U_{t\cdot\mathbf 1_d})}=\frac{\nu(U_{\mathbf y\wedge t\cdot\mathbf 1_d}\setminus U_{\mathbf y})}{\nu(U_{t\cdot\mathbf 1_d})}=\frac{\nu([\mathbf 0_d,\mathbf y]\cap U_{t\cdot\mathbf 1_d})}{\nu(U_{t\cdot\mathbf 1_d})};
\end{align*}
to establish the second last equality, we note that $U_{\mathbf y}\subseteq U_{\mathbf y\wedge t\cdot\mathbf 1_d}$ by the definition of $U_{\cdot}$, then the equality is immediate by the very construction of $\nu(\cdot)$. We anticipate similar arguments that lead to the following limiting joint law for $\mathbf Y_p$ given $U_{t\cdot\mathbf 1_d}$ as $p\uparrow1$, for any Borel set $B\subseteq[0,\infty)^d$:
\begin{equation}\label{eq:cond}
\mathbb P(\mathbf Y_p\in B\mid\mathbf Y_p\in U_{t\cdot\mathbf 1_d})\to\frac{\nu(B\cap U_{t\cdot\mathbf 1_d})}{\nu(U_{t\cdot\mathbf 1_d})}=:P_t(B).
\end{equation}
Consider any sequence $p_n\uparrow1$ along
which $\mathbf r_n:=(r_{1,n},\ldots,r_{d,n})^\top$, with $r_{i,n}:=\frac{x_{i,p_n}}{\sum_{j=1}^d x_{j,p_n}}$ for  $i=1,\dots,d$, converges to
some $\mathbf r:=(r_1,\dots,r_d)^\top$; clearly all entries of $\mathbf r_n$ are nonnegative and sum up to 1, which is also the case for $\mathbf r$. Note that $$\mathbb{P}\left(S_d>t\sum_{j=1}^dx_{j,p_n}\right)=\mathbb{P}\left(\sum_{i=1}^d\frac{X_i}{x_{i,p_n}}\frac{x_{i,p_n}}{\sum_{j=1}^dx_{j,p_n}}>t\right)=\mathbb{P}\left(\sum_{i=1}^dr_{i,n}Y_{i,p_n}=\mathbf r_n^\top\mathbf Y_{p_n}>t\right),$$
hence as $p_n\uparrow1$, $\mathbf r_n\to\mathbf r$ such that $\{\mathbf y\in[0,\infty)^d:\mathbf r_n^\top\mathbf y>t\}\to \{\mathbf y \in [0,\infty)^d:\mathbf r^\top\mathbf y>t\}=:H_t(\mathbf r)$, and in light of \eqref{eq:cond} together with an application of Theorem 2.7(v) of \cite{vdv98} and then the continuous mapping theorem
(see Theorem 2.3(i) of \cite{vdv98}), we foresee that:
\[
\begin{aligned}
\frac{\mathbb{P}\left(S_d>t\sum_{j=1}^dx_{j,p_n}\right)}{1-p_n}
&=\mathbb P\left(\mathbf r_n^\top\mathbf Y_{p_n}>t
                   \mid\mathbf Y_{p_n}\in U_{t\cdot\mathbf 1_d}\right)
  \frac{\mathbb P(\mathbf Y_{p_n}\in U_{t\cdot\mathbf 1_d})}{1-p_n}\\
  &\to P_t(H_t(\mathbf r)) \nu(U_{t\cdot\mathbf 1_d})=\frac{\nu(H_t(\mathbf r)\cap U_{t\cdot\mathbf 1_d})}{\nu(U_{t\cdot\mathbf 1_d})}\nu(U_{t\cdot\mathbf 1_d})
 =\nu(H_t(\mathbf r)),
\end{aligned}
\]
where the last equality holds since $\sum_{i=1}^d r_ix_i>t$ implies $\max_{i\in\{1,\dots,d\}}x_i>t$ by pigeonhole principle, indicating that $H_t(\mathbf r)\subseteq U_{t\cdot\mathbf 1_d}$. Similar to \eqref{eq:nu_u}, with the indicator taking 1 if $\sum_{i=1}^dr_isW_i^{\frac{1}{\alpha}}>t$, we have:
\[
\nu(H_t(\mathbf r))=\mathbb E\left[-s^{-\alpha}\bigg|_{t\left(\sum_{i=1}^dr_iW_i^{\frac{1}{\alpha}}\right)^{-1}}^\infty\right]
=t^{-\alpha}\mathbb E\left[\left(\sum_i r_iW_i^{1/\alpha}\right)^\alpha\right]
=t^{-\alpha}\left(K_\alpha(\mathbf r)\right)^\alpha.
\]
By a similar argument as in the proof of Theorem \ref{thm:main}, we then have $D_{p_n}\to K_\alpha(\mathbf r)$, hence $D_{p_n}-K_\alpha(\mathbf r_n)\to0$ by continuity of $K_\alpha$. Meanwhile, by homogeneity, $K_\alpha(\mathbf r_n)=\frac{K_\alpha(x_{1,p_n},\ldots,x_{d,p_n})}{\sum_{j=1}^dx_{j,p_n}}$. \eqref{eq:discussion_Dp} then follows since every sequence $p_n\uparrow1$ has a subsequence along which the proportions converge, based on which we perceive the following:

i) For $\alpha>1$ and any nonzero $\mathbf z\in[0,\infty)^d$, applying Jensen's inequality to $x\mapsto x^\alpha$ with weights $\frac{z_i}{\sum_{j=1}^d z_j}$, $i=1,\dots,d$ gives:
\[
\left(\frac{K_\alpha(\mathbf z)}{\sum_{j=1}^d z_j}\right)^\alpha
=\mathbb E\left[\left(\sum_{i=1}^d\frac{z_i}{\sum_{j=1}^d z_j}
                         W_i^{1/\alpha}\right)^\alpha\right]
\le \sum_{i=1}^d\frac{z_i}{\sum_{j=1}^d z_j}\mathbb E [W_i]=1.
\]
Taking $z_i=\V@R_p(X_i)$ in this bound, \eqref{eq:discussion_Dp}
yields $\limsup_{p\uparrow1}D_p\le1$.

ii) For $0<\alpha<1$, concavity
reverses the inequality, giving $\liminf_{p\uparrow1}D_p\ge1$.

iii) Finally, for $\alpha=1$, $K_1(\mathbf z)=\sum_{i=1}^d z_i$, and hence
$D_p\to1$ as $p\uparrow1$.

Other future interesting works may study second-order approximations for the
diversification ratio, or to obtain similar asymptotic results for other tail risk measures (for Expected Shortfall, which is naturally subadditive as a coherent risk measure, only asymptotic additivity would be of interest).

\section{Acknowledgements}
This work has undergone substantial changes in terms of idea and methodology from a preliminary version. We thank Jia Zeng for her contributions in preparing that preliminary version, which also constitutes part of the PhD thesis \cite{zeng2022contributions} at the University of Hong Kong and King's College London. Yongzhao Chen acknowledges the financial support from HKGRF-UGC/FDS14/P03/23 with the project title ``InsurTech: Risk Classification and Premium Calibration with Data Analytics''. Ka Chun Cheung was partially supported by a grant from the Research Grants Council of the Hong Kong Special Administrative Region, China (Project No.\ 17324516). Hok Kan Ling gratefully acknowledges the support by NSERC Grant RGPIN/03124-2021. Qihe Tang acknowledges support from the Australian Government through the Australian Research Council’s Discovery Projects funding scheme (DP200101859, DP220100090, and DP250100078). Phillip Yam acknowledges the financial support from 
HKGRF-14300123 with the project title ``Well-posedness of Some Poisson-driven Mean Field
Learning Models and their Applications'', and HKGRF-14300025 with the project title ``A Generic Theory for Stochastic Control against Fractional Brownian Motions''. The work described in this article was also supported by
a grant from the Germany/Hong Kong Joint Research Scheme sponsored by the Research Grants
Council of Hong Kong and the German Academic Exchange Service of Germany (Reference No.
G-CUHK411/23).

\bibliographystyle{plain}
\bibliography{references_revised}

\appendix

\section{Summary for Asymptotic Behavior of $D_p$}

\begin{table}[!htbp]
\centering
\caption{High-confidence behavior of the diversification ratio $D_p$.}
\label{tab:main-results}
\small
\renewcommand{\arraystretch}{1.28}
\setlength{\tabcolsep}{6pt}

\begin{tabularx}{\textwidth}{
    @{}>{\raggedright\arraybackslash}p{0.34\textwidth}
    >{\raggedright\arraybackslash}X
    >{\centering\arraybackslash}p{0.27\textwidth}@{}
}
\toprule
\textbf{Marginal configuration}
&
\textbf{Dependence}
&
\textbf{As $p\uparrow1$}
\\
\midrule

\multicolumn{3}{@{}l}{
\textit{Fr\'echet: $I$ denotes the set of dominating risks}
} \\[2pt]

$|I|\geq2$, with common tail index $\alpha$
&
Regularly varying Archimedean survival dependence on $I$
&
$\displaystyle
\begin{cases}
\liminf D_p\geq1, & 0<\alpha<1;\\
D_p\to1,          & \alpha=1;\\
\limsup D_p\leq1, & \alpha>1.
\end{cases}$
\\[6pt]

$|I|=1$
&
Arbitrary
&
$\displaystyle D_p\to1$
\\

\addlinespace[4pt]
\midrule

\multicolumn{3}{@{}l}{
\textit{Gumbel: infinite right endpoints}
} \\[2pt]

Gumbel portfolio
&
Arbitrary
&
$\displaystyle \limsup D_p\leq1$
\\

Gumbel portfolio
&
Positive diagonal upper-tail dependence
&
$\displaystyle D_p\to1$
\\

Gumbel portfolio with finite-endpoint additions
&
Arbitrary
&
$\displaystyle \limsup D_p\leq1$
\\

\addlinespace[4pt]
\midrule

\multicolumn{3}{@{}l}{
\textit{Finite right endpoints}
} \\[2pt]

$x_{F_i}<\infty$ for all $i$
&
Strict Archimedean survival dependence
&
$\displaystyle D_p\to1$
\\

\bottomrule
\end{tabularx}
\end{table}

\section{Proof of Theorem \ref{thm:limsup-liminf-sharp}}\label{sect:proof_of_sharpness}

\begin{proof}
Consider the following two marginal survival functions for all sufficiently
large $x$:
\[
\overline F_1(x)=x^{-\alpha},
\quad
\overline F_2(x)=x^{-\alpha}L(x),
\quad
\text{ where }L(x)
=
\exp\left\{
\sqrt{\ln\ln x} 
\sin\bigl(\ln\ln\ln x\bigr)
\right\}.
\]

Obviously, $\overline{F}_2$ defines an admissible
regularly varying tail; indeed, write
$h(u):=\sqrt{u}\sin(\ln u)$, then
$\ln L(x)=h(\ln\ln x)$. Since
$h'(u)
=
\frac{\sin(\ln u)}{2\sqrt{u}}
+
\sqrt{u}\cos(\ln u)\frac{1}{u}
=
O(u^{-\frac{1}{2}})$ as $u\to\infty$,
and for every fixed $t>0$,  $\ln\ln(tx)-\ln\ln x\to0$ as $x\to\infty$, thus the mean value theorem gives:
\[
\ln\frac{L(tx)}{L(x)}
=
h(\ln\ln(tx))-h(\ln\ln x)
\rightarrow0;
\]
hence
$\frac{L(tx)}{L(x)}\to1$ as $x\to\infty$, and therefore $L$ is slowly varying. Also note that $x^{-\alpha}L(x)$ is continuous in $x$, and for all sufficiently large $x$,
\[
x\frac{d}{dx}\ln L(x)
=
\frac{h'(\ln\ln x)}{\ln x}
\rightarrow0
\implies
x\frac{d}{dx}\ln\bigl(x^{-\alpha}L(x)\bigr)
=
-\alpha+o(1)<0,
\]
hence $x^{-\alpha}L(x)$ is eventually
strictly decreasing. With a proper modification upon smaller values of $x$,
$\overline F_1$ and $\overline F_2$ may therefore be taken to be continuous
survival functions. In particular,
$\overline F_1,\overline F_2\in\operatorname{RV}_{-\alpha}$, so both marginals
belong to $\operatorname{MDA}(\Phi_\alpha)$.

We next couple $X_1$ and $X_2$ by a survival Clayton copula with the generator
$\psi_\gamma(u)=u^{-\gamma}-1$ for
$\gamma>0$,
and as before, we write
$S_2=X_1+X_2$. We now construct two sequences of confidence levels as follows:
\begin{enumerate}
\item[i)] Let
$x_n^{(1)}
:=
\exp\left\{\exp\left(e^{2\pi n}\right)\right\}$,
$p_n^{(1)}
:=
1-\bigl(x_n^{(1)}\bigr)^{-\alpha}$.
Then $\ln\ln\ln x_n^{(1)}=2\pi n$, and hence
$L(x_n^{(1)})=1$. For all sufficiently large $n$,
$\overline F_1(x_n^{(1)})
=
\overline F_2(x_n^{(1)})
=
1-p_n^{(1)}$.
Since the marginal tails are continuous and strictly decreasing in this
region,
$\V@R_{p_n^{(1)}}(X_1)
=
\V@R_{p_n^{(1)}}(X_2)
=
x_n^{(1)}$, therefore
\begin{equation}
\V@R_{p_n^{(1)}}(X_1)
+
\V@R_{p_n^{(1)}}(X_2)
=
2x_n^{(1)}.
\label{eq:sharp-denom-c}
\end{equation}

Let $\Theta^{(1)}:=\{x_n^{(1)}\}_{n\geq1}$. Note that
$\frac{\overline F_2(x_n^{(1)})}{\overline F_1(x_n^{(1)})}=L(x_n^{(1)})=1$, and $w_1^1(\Theta^{(1)})=w_1^2(\Theta^{(1)})=1$; furthermore, for every fixed $t>0$, by its slowly  varying nature,
$L(tx_n^{(1)})\sim L(x_n^{(1)})=1$, and hence
$\frac{\overline F_2(tx_n^{(1)})}
     {\overline F_1(tx_n^{(1)})}=\frac{\overline F_2(tx_n^{(1)})}
     {\overline F_2(x_n^{(1)})}\frac{\overline F_2(x_n^{(1)})}
     {\overline F_1(x_n^{(1)})}\frac{\overline F_1(x_n^{(1)})}
    {\overline F_1(tx_n^{(1)})}
\rightarrow t^{-\alpha}\cdot 1\cdot t^\alpha=1$ as $n\to\infty$;
thus Theorem \ref{thm:Qexist} applies to
$\Theta_t^{(1)}:=\{tx_n^{(1)}\}_{n\geq1}$, with $w_j^i(\Theta_t^{(1)})=1$ for $i,j\in\{1,2\}$. Following similar arguments as in the proof of Theorem \ref{thm:main}, since the integral density representation \eqref{eq:small_h}-\eqref{qjint} depends on
the sequence $\Theta_t^{(1)}$ only through these limiting tail ratios, $\lim_{n\to\infty}\frac {\mathbb{P}\left(S_2>tx_n^{(1)}\right)}{\overline{F}_1(tx_n^{(1)})}=:Q_1(\gamma,\alpha;\Theta_t^{(1)})=Q_1(\gamma,\alpha;\Theta^{(1)}):=\lim_{n\to\infty}\frac {\mathbb{P}\left(S_2>x_n^{(1)}\right)}{\overline{F}_1(x_n^{(1)})}$.
Since
$\frac{\overline F_1(tx_n^{(1)})}{1-p_n^{(1)}}=t^{-\alpha}$, we have that as $n\to\infty$:
\begin{equation}
\frac{\mathbb P(S_2>tx_n^{(1)})}
     {1-p_n^{(1)}}
\rightarrow
Q_1(\gamma,\alpha;\Theta^{(1)})t^{-\alpha},
\qquad t>0.
\label{eq:sharp-tail-c}
\end{equation}

Analogous to the arguments leading to \eqref{eq:quantile-liminf} and
\eqref{eq:quantile-limsup}, if
$t<(Q_1(\gamma,\alpha;\Theta^{(1)}))^{\frac{1}{\alpha}}$, then $Q_1(\gamma,\alpha;\Theta^{(1)})t^{-\alpha}>1$, so
\eqref{eq:sharp-tail-c} implies, for all sufficiently large $n$,
$\mathbb P(S_2>tx_n^{(1)})>1-p_n^{(1)}$, and therefore
$\V@R_{p_n^{(1)}}(S_2)>tx_n^{(1)}$. Similarly, if
$t>(Q_1(\gamma,\alpha;\Theta^{(1)}))^{\frac{1}{\alpha}}$, then $Q_1(\gamma,\alpha;\Theta^{(1)})t^{-\alpha}<1$, and hence
eventually $\V@R_{p_n^{(1)}}(S_2)\leq tx_n^{(1)}$. Letting
$t\uparrow (Q_1(\gamma,\alpha;\Theta^{(1)}))^{\frac{1}{\alpha}}$ and $t\downarrow (Q_1(\gamma,\alpha;\Theta^{(1)}))^{\frac{1}{\alpha}}$ yields
$\frac{\V@R_{p_n^{(1)}}(S_2)}
     {x_n^{(1)}}
\rightarrow
(Q_1(\gamma,\alpha;\Theta^{(1)}))^{\frac{1}{\alpha}}$,
then together with \eqref{eq:sharp-denom-c},
\begin{equation}
D_{p_n^{(1)}}=\frac{\V@R_{p_n^{(1)}}(S_2)}{\V@R_{p_n^{(1)}}(X_1)+\V@R_{p_n^{(1)}}(X_2)}=\frac{V@R_{p_n^{(1)}}(S_2)}{2x_n^{(1)}}
\rightarrow
\frac{(Q_1(\gamma,\alpha;\Theta^{(1)}))^{\frac{1}{\alpha}}}{2}=:d_{\alpha,\gamma}.
\label{eq:sharp-D-c}
\end{equation}
Recall that
$w_1^1(\Theta^{(1)})=w_1^2(\Theta^{(1)})=1$ as mentioned above, so the corresponding
comonotonic limit for $Q_1(\gamma,\alpha;\Theta^{(1)})$, as defined in Proposition \ref{prop:Q-comonotonic}, is:
$$Q_{1,\mathrm{com}}(\alpha;\Theta^{(1)})=\left(\sum_{i=1}^2\left(w_1^i(\Theta^{(1)})\right)^{\frac{1}{\alpha}}\right)^\alpha=(1+1)^\alpha=2^\alpha.$$
By Proposition \ref{prop:Q-comonotonic}, the comparison is strict for
finite $\gamma$, and therefore $Q_1(\gamma,\alpha;\Theta^{(1)})>2^\alpha$ (resp. $<2^\alpha$) if $0<\alpha<1$ (resp. $\alpha>1$).
Consequently,
\begin{equation}\label{eq:sharp-c-sub-super}
d_{\alpha,\gamma}\left\{\begin{aligned}
    &>1,\text{ for }0<\alpha<1.\\
&<1,\text{ for } \alpha>1.
\end{aligned}\right.
\end{equation}

\item[ii)] Let
$x_n^{(2)}
:=
\exp\left\{
\exp\left(e^{\frac{3\pi}{2}+2\pi n}\right)
\right\}$ and
$p_n^{(2)}
:=
1-\bigl(x_n^{(2)}\bigr)^{-\alpha}$,
then $\ln\ln\ln x_n^{(2)}=\frac{3\pi}{2}+2\pi n$, so
$L(x_n^{(2)})
=
\exp\left\{
-\sqrt{\ln\ln x_n^{(2)}}
\right\}
\rightarrow0$.
Since $\overline F_1(x_n^{(2)})=1-p_n^{(2)}$, and recalling the modification of $\overline{F}_1$ and $\overline{F}_2$ on smaller values of $x$ mentioned in establishing their regularly varying nature, we have
$\V@R_{p_n^{(2)}}(X_1)=x_n^{(2)}$ for all sufficiently large $n$. Since $p_n^{(2)}\uparrow1$, we have
$\V@R_{p_n^{(2)}}(X_2)\to\infty$, and the continuity and eventual strict monotonicity of
$\overline F_2$ give
$\overline F_2(\V@R_{p_n^{(2)}}(X_2))
=
1-p_n^{(2)}
=
\bigl(x_n^{(2)}\bigr)^{-\alpha}$, whereas
$\overline F_2(x_n^{(2)})
=
\bigl(x_n^{(2)}\bigr)^{-\alpha}L(x_n^{(2)})$; hence
$\frac{\overline F_2(\V@R_{p_n^{(2)}}(X_2))}
     {\overline F_2(x_n^{(2)})}
=
\frac{1}{L(x_n^{(2)})}
\rightarrow\infty$.
Since $L(x_n^{(2)})<1$ eventually, we also have $\V@R_{p_n^{(2)}}(X_2)<x_n^{(2)}$.
Fix $t>0$ and
$\epsilon\in(0,t)$. By a similar argument as that leading to \eqref{eq:S-tail}, since $X_1,X_2\geq0$,
we have
$\mathbb P(X_1>tx_n^{(2)})
\leq
\mathbb P(S_2>tx_n^{(2)})
\leq
\mathbb P\bigl(X_1>(t-\epsilon)x_n^{(2)}\bigr)
+
\mathbb P\bigl(X_2>\epsilon x_n^{(2)}\bigr)$;
dividing by 
$1-p_n^{(2)}=(x_n^{(2)})^{-\alpha}$ gives:
\begin{align}
  \frac{\mathbb P(S_2>tx_n^{(2)})}{1-p_n^{(2)}}
&\geq \frac{\mathbb P(X_1>tx_n^{(2)})}{(x_n^{(2)})^{-\alpha}};\label{eq:sequence2_lower}\\
\frac{\mathbb P(S_2>tx_n^{(2)})}{1-p_n^{(2)}}
&\leq
\frac{\mathbb P\bigl(X_1>(t-\epsilon)x_n^{(2)}\bigr)}{(x_n^{(2)})^{-\alpha}}
+
\frac{\mathbb P\bigl(X_2>\epsilon x_n^{(2)}\bigr)}{(x_n^{(2)})^{-\alpha}}.\label{eq:sequence2_upper}
\end{align}
Note that as $n\to\infty$, 
$\frac{\mathbb P(X_1>tx_n^{(2)})}{(x_n^{(2)})^{-\alpha}}=\frac{\mathbb P(X_1>tx_n^{(2)})}{\mathbb P(X_1>x_n^{(2)})}\to t^{-\alpha}$ and likewise
$\frac{\mathbb P\bigl(X_1>(t-\epsilon)x_n^{(2)}\bigr)}{(x_n^{(2)})^{-\alpha}}=\frac{\mathbb P\bigl(X_1>(t-\epsilon)x_n^{(2)}\bigr)}{\mathbb P(X_1>x_n^{(2)})}\to (t-\epsilon)^{-\alpha}$, while
$$
\frac{\mathbb P(X_2>\epsilon x_n^{(2)})}
     {1-p_n^{(2)}}
=\frac{(\epsilon x_n^{(2)})^{-\alpha}L(\epsilon x_n^{(2)})}{(x_n^{(2)})^{-\alpha}}=
\epsilon^{-\alpha}L(\epsilon x_n^{(2)})=
\epsilon^{-\alpha}
\frac{L(\epsilon x_n^{(2)})}
     {L(x_n^{(2)})}
L(x_n^{(2)})
\rightarrow0,
$$
by the slowly varying nature of $L$ and also $L(x_n^{(2)})\to0$. Applying the respective limits to \eqref{eq:sequence2_lower} and \eqref{eq:sequence2_upper} yields:
\[
t^{-\alpha}
\leq
\liminf_{n\to\infty}
\frac{\mathbb P(S_2>tx_n^{(2)})}
     {1-p_n^{(2)}}
\leq
\limsup_{n\to\infty}
\frac{\mathbb P(S_2>tx_n^{(2)})}
     {1-p_n^{(2)}}
\leq
(t-\epsilon)^{-\alpha},
\]
and letting $\epsilon\downarrow0$ gives
\begin{equation}
\frac{\mathbb P(S_2>tx_n^{(2)})}
     {1-p_n^{(2)}}
\rightarrow
t^{-\alpha},
\qquad t>0.
\label{eq:sharp-tail-d}
\end{equation}

Similar to the argument before \eqref{eq:sharp-D-c} in Case i), if
$t<1$, then $t^{-\alpha}>1$, so
\eqref{eq:sharp-tail-d} implies, for all sufficiently large $n$,
$\mathbb P(S_2>tx_n^{(2)})>1-p_n^{(2)}$, and therefore
$\V@R_{p_n^{(2)}}(S_2)>tx_n^{(2)}$; and if
$t>1$, then $t^{-\alpha}<1$, and hence
eventually $\V@R_{p_n^{(2)}}(S_2)\leq tx_n^{(2)}$. Letting
$t\uparrow 1$ and $t\downarrow 1$ yields
$\frac{\V@R_{p_n^{(2)}}(S_2)}{x_n^{(2)}}\to1$. On the other hand, by the very construction of $X_1$,
\[
\frac{
\V@R_{p_n^{(2)}}(X_1)+\V@R_{p_n^{(2)}}(X_2)
}{
x_n^{(2)}
}
=
1+\frac{\V@R_{p_n^{(2)}}(X_2)}{x_n^{(2)}}
\rightarrow1,\text{ as }n\to\infty;
\]
in particular, to see that
$\frac{\V@R_{p_n^{(2)}}(X_2)}{x_n^{(2)}}\rightarrow0$,
we assume otherwise, then there would exist a subsequence $\{p_{n_k}^{(2)}\}$ and
$c\in(0,1]$ such that $\frac{\V@R_{p_{n_k}^{(2)}}(X_2)}{x_{n_k}^{(2)}}\to c$. Since
$\overline F_2\in\operatorname{RV}_{-\alpha}$, the uniform convergence
theorem for regularly varying functions (Theorem 1.5.2 of \cite{bgt87}) would then give
$\frac{\overline F_2(\V@R_{p_{n_k}^{(2)}}(X_2))}
     {\overline F_2(x_{n_k}^{(2)})}
\rightarrow
c^{-\alpha}<\infty$; yet recalling that $L(x_n^{(2)})\to0$ and $\overline F_2(\V@R_{p_n^{(2)}}(X_2))
=
(x_n^{(2)})^{-\alpha}$ as mentioned above, $\frac{\overline F_2(\V@R_{p_{n_k}^{(2)}}(X_2))}
     {\overline F_2(x_{n_k}^{(2)})}=\frac{(x_{n_k}^{(2)})^{-\alpha}}{(x_{n_k}^{(2)})^{-\alpha}L(x_{n_k}^{(2)})}=\frac{1}{L(x_{n_k}^{(2)})}\to\infty$, which leads to a contradiction. Finally, we can conclude that:
\begin{equation}
D_{p_n^{(2)}}=\frac{\V@R_{p_n^{(2)}}(S_2)}{\V@R_{p_n^{(2)}}(X_1)+\V@R_{p_n^{(2)}}(X_2)}=\frac{\frac{\V@R_{p_n^{(2)}}(S_2)}{x_n^{(2)}}}{\frac{\V@R_{p_n^{(2)}}(X_1)+\V@R_{p_n^{(2)}}(X_2)}{x_n^{(2)}}}\rightarrow\frac{1}{1}=1,\text{ as }n\to\infty.
\label{eq:sharp-D-d}
\end{equation}
\end{enumerate}
Clearly, the limit of \eqref{eq:sharp-D-c}, \eqref{eq:sharp-c-sub-super} in Case i) is not equal to that of \eqref{eq:sharp-D-d} in Case ii) whenever $\alpha\neq1$, hence $D_p$ does
not converge as $p\uparrow1$ in this specific example.

\end{proof}

\end{document}